\documentclass{amsart}
\usepackage[a4paper, total={6in, 8in}]{geometry} 	%Paper size
\usepackage[utf8]{inputenc}
\usepackage[english]{babel}
\usepackage{mypreamble}
\newcommand{\Rgrmod}{R\mathrm{Mod}}
\newcommand{\Zgrmod}{\Z\mathrm{Mod}}
\newcommand{\GrCob}{\mathrm{GCob}}
\newcommand{\bolddet}{\mathbf{det}}
\renewcommand{\Cob}{\mathrm{Cob}}
\newcommand{\fatCob}{\mathrm{fCob}}
\newcommand{\fatGrCob}{\mathrm{fGCob}}

\newcommand{\pGrCob}{\mathrm{pGCob}}
\newcommand{\comp}{\mathrm{comp}}

\newcommand{\wotimes}{\widetilde{\otimes}}
\newcommand{\starotimes}{\widehat\otimes^*}
\newcommand{\shriekotimes}{\widehat\otimes^!}

\DeclareRobustCommand{\assgraph}{\vcenter{\hbox{\includegraphics[height=1.2\fontcharht\font`\B]{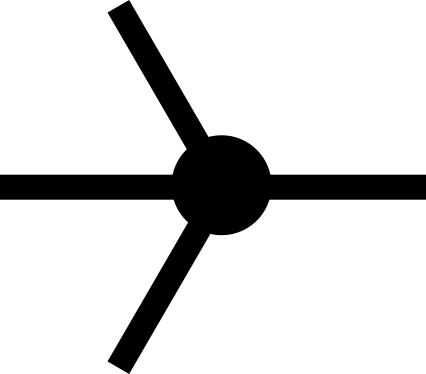}}}}
\DeclareRobustCommand{\coassgraph}{\vcenter{\hbox{\includegraphics[height=1.2\fontcharht\font`\B]{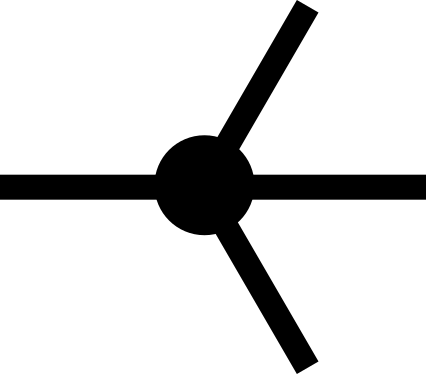}}}}
\DeclareRobustCommand{\idgraph}{\vcenter{\hbox{\includegraphics[height=0.4\fontcharht\font`\B]{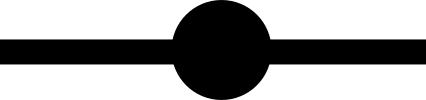}}}}
\DeclareRobustCommand{\frobgraph}{\vcenter{\hbox{\includegraphics[height=\fontcharht\font`\B]{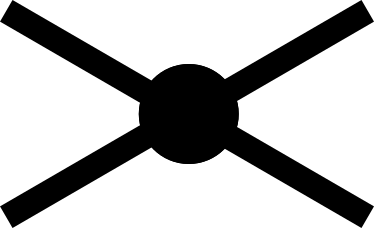}}}}
\DeclareRobustCommand{\symgraph}{\vcenter{\hbox{\includegraphics[height=1.2\fontcharht\font`\B]{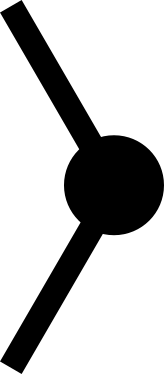}}}}
\DeclareMathAlphabet{\mathpzc}{OT1}{pzc}{m}{it}

\title{Graded Frobenius Algebras}
\author{Jonathan Clivio}
\date{\today}

\begin{document}
\begin{abstract}
    We construct a PROP which encodes 2D-TQFTs with a grading. This defines a graded Frobenius algebra as algebras over this PROP. We also give a description of graded Frobenius algebras in terms of maps and relations. This structure naturally arises as the cohomology of manifolds, loop homology and Hochschild homology of Frobenius algebras. In addition, we give a comprehensive description of the signs that arise in suspending algebras over PROPs.
\end{abstract}

\maketitle

\section{Introduction}
The goal of this paper is to clarify what is meant by a graded Frobenius algebra. To this end, we introduce a PROP that describes a Frobenius algebra with multiplication in a fixed degree $c$ and comultiplication of a fixed degree $d$.\\

There is a need for such a definition, because the straight-forward attempt of considering Frobenius algebras in graded abelian groups runs into problems. Defining a Frobenius algebra in graded abelian groups is unproblematic if all maps are degree-preserving, i.e., they are maps in degree 0. But this attempt fails to describe some non-trivial examples if maps of different degrees are considered: for example, we consider a graded abelian group $A$ with a unital multiplication $\mu\colon A\otimes A\to A$ and a counital comultiplication $\nu\colon A\to A\otimes A$ on it satisfying the usual Frobenius relation
\begin{align}\label{eq:frob intro}
    (\mu \otimes \id)\circ (\id\otimes \nu)=\nu\circ \mu=(\id\otimes \mu)\circ (\nu \otimes \id).
\end{align}
We now assume that $\mu$ and $\nu$ are of some degree $|\mu|=c$ and $|\nu|=d$. If $c+d$ is odd and we work over a field where $1\neq -1$, the only example that is described is the trivial algebra $A=0$. We prove this in Proposition \ref{prop:trivial}.\\

Another attempt to describe examples of Frobenius algebras with maps in non-zero degree is to suspend a Frobenius algebra. Suppose $A$ is a Frobenius algebra in graded abelian groups with multiplication $\mu$ and $\nu$ in degree 0. Suspending $A$ gives an algebra $\Sigma A$ with a multiplication $\Sigma \mu$ in degree $-1$ and a comultiplication $\Sigma \nu$ in degree $1$. Suspension and desuspension of algebras therefore leaves the sum of the degree of the multiplication and the degree of the comultiplication invariant. Thus an example of $A$ with a multiplication $\mu$ in degree $0$ and a comultiplication in degree $d\neq 0$ cannot be described as the suspension of an algebra with only operations in degree 0.\\
Moreover, the canonical map 
\begin{align}\label{eq:Sigma maps}
    \Sigma\colon \Hom(A^{\otimes k},A^{\otimes l})\to \Hom((\Sigma A)^{\otimes k},(\Sigma A)^{\otimes l})
\end{align}
is only functorial and symmetric monoidal up to sign. We compute these signs in Lemma \ref{lem:Sigma signs}. These signs also illustrate why Frobenius algebras in graded abelian groups are not the natural concept to consider. For example if $\mu$ and $\nu$ satisfy the Frobenius relation as in \eqref{eq:frob intro}, the signs in the functoriality and monoidality of $\Sigma$ in \eqref{eq:Sigma maps} induce the following sign in the Frobenius relation:
\begin{align*}
    (\Sigma\mu \otimes \id)\circ (\id\otimes \Sigma\nu)=(-1)\Sigma\nu\circ \Sigma\mu=(\id\otimes \Sigma\mu)\circ (\Sigma\nu \otimes \id).
\end{align*}
Therefore, the natural concept of graded Frobenius algebras will have signs appearing in the relations.

\subsection{The Graded Frobenius PROP}
In this paper, we give a definition of graded Frobenius algebras with a multiplication $\mu$ in degree $c$ and a comultiplication $\nu$ in degree $d$ for any choice $c,d\in \Z$ (see Definition \ref{def:graded Frob alg}). This definition solves the two problems described above: we can describe non-trivial examples for any choice of $c,d\in \Z$ and the signs that appear in this definition are stable under suspending algebras. We make this second claim precise when we discuss the relations in subsection \ref{subsec:maps rel}.\\
We describe a category of graph cobordisms in Definition \ref{def:ass 1-cat}. This lets us model every cobordism by an equivalence class of graphs.\\
For a graph $G$ representing a cobordism, we consider a twisting by the determinant $\bolddet(G,\partial_{in})$ of the homology relative to the incoming boundary. The determinant $\bolddet(G,\partial_{in})$ is a graded rank one free abelian group concentrated in degree $-\chi(G,\partial_{in})$, the negative Euler characteristic. Analogously, we consider $\bolddet(G,\partial_{out})$.\\
This lets us define the categories of graded graph cobordism $\GrCob_{c,d}$, fat graph cobordisms $\fatGrCob_{c,d}$ and planar graph cobordisms $\pGrCob_{c,d}$ in Definition \ref{def:grcob}. These categories have morphisms given by elements in the following graded abelian groups
\begin{align*}
    \bolddet_{c,d}(G):=\bolddet(G,\partial_{in})^{\otimes c}\otimes \bolddet(G,\partial_{out})^{\otimes d}.
\end{align*}
In these categories, there is a morphism $*\sqcup*\to *$. It models a multiplication of the graded Frobenius algebra $*$. It is given by a generator $\omega(\multi)\in \bolddet_{c,d}(\multi)$. Here $\multi$ is the graph with two inputs to the left, one output to the right and one three-valent vertex in the middle. It has relative Euler characteristic $-\chi(\multi,\partial_{in})=1$ and $-\chi(\multi,\partial_{out})=0$. Therefore $\bolddet_{c,d}(\multi)$ is concentrated in degree $c$ and the multiplication $\omega(\multi)$ is in degree $c$. Analogously, the comultiplication is given by a generator $\omega(\comulti)\in \bolddet_{c,d}(\comulti)$ in degree $d$.\\
The categories $\GrCob_{c,d}$, $\fatGrCob_{c,d}$ and $\pGrCob_{c,d}$ let us thus attain any combination in degrees $c,d$. Moreover, $\GrCob_{c,d}$ is modelled by graphs and describes $(c,d)$-graded \textit{commutative} Frobenius algebras, $\fatGrCob_{c,d}$ has morphisms given by fat graphs and describes $(c,d)$-graded \textit{symmetric} Frobenius algebras and $\pGrCob_{c,d}$ has morphisms given by planar graphs and describes $(c,d)$-graded Frobenius algebras\\

This twisting was previously described by Godin in \cite{godin2007higher} and by Wahl and Westerland in \cite{wahlwesterland}. However, in contrast to Godin who models commutative Frobenius algebras by closed surfaces, we work with graphs. In particular, twisting a surface $S$ by $\bolddet(S,\partial_{in})$ corresponds to a twisting by $\bolddet(G,\partial_{in})\otimes \bolddet(G,\partial_{out})$ where $G$ is a graph modelling the cobordism $S$. We prove this relation between the twistings in Proposition \ref{prop:1closed=2open twists}. In particular, working with graphs rather than closed surfaces for commutative graded Frobenius algebras allows us to describe examples of algebras where $c+d$ is odd.

\subsection{Suspending PROPs}

Above we observed that suspending an algebra with multiplication in degree $0$ and comultiplication in degree $0$ gives an algebra with multiplication in degree $-1$ and a comultiplication in degree $1$. These shifts in degree hold for any choice in degrees $c,d$.\\
This can be seen on the level of PROPs: we define the usual suspension PROP $\Sigma$ in Definition \ref{def:susp prop}. This PROP has the following universal property: a graded abelian group $A$ has the structure of an algebra over a PROP $P$ if and only if $\Sigma A$ has the structure of an algebra over $\Sigma\otimes P$ (see Proposition \ref{prop:univ susp prop}). On the other hand, we give natural isomorphisms of PROPs
\begin{align*}
    \Sigma \otimes \GrCob_{c,d}&\cong \GrCob_{c-1,d+1},\\
    \Sigma\otimes  \fatGrCob_{c,d}&\cong \fatGrCob_{c-1,d+1},\\
    \Sigma\otimes  \pGrCob_{c,d}&\cong \pGrCob_{c-1,d+1}
\end{align*}
in Proposition \ref{prop:suspending c,d}. This, in particular, shows that suspending a $(c,d)$-graded Frobenius algebra gives a $(c-1,d+1)$-graded Frobenius algebra.\\
This is in contrast to the relation between $(c,d)$-graded Frobenius algebras and $(c+1,d)$-graded Frobenius algebras. In Remark \ref{rem:one sided suspension}, we give a construction that produces a $(c+1,d)$-graded Frobenius algebra from a $(c,d)$-graded Frobenius algebra. This construction could be called a one-sided suspension of graded Frobenius algebras as it only affects the degree of the multiplication and not the coproduct. However in contrast to the usual suspension, this one-sided suspension is not an equivalence of categories.

\subsection{Maps, Relations and Signs}\label{subsec:maps rel}

In the ungraded case, there are several equivalent descriptions of a Frobenius algebra. We contrast two different types of definitions:
\begin{enumerate}
    \item Descriptions in terms of maps with relations: an object $A$ together with a selection of the following maps: a multiplication $\mu$, a unit $\eta$, a comultiplication $\nu$ and a counit $\varepsilon$ that satisfy certain relations.
    \item Descriptions in terms of PROPs: a 2-dimensional topological quantum field theory (TQFT). This is a (symmetric) monoidal functor out of a 2-dimensional cobordism category.
\end{enumerate}
There is a myriad of different variants of the two flavours of definitions with optional commutativity, unitality, counitality. The equivalence between these two types of definitions is described in \cite{abrams1996two,lauda2008open, kock2004frobenius}.\\
We explained above that we define a $(c,d)$-graded Frobenius algebra in terms of a PROP. In our main theorem, we relate our definition to a description in terms of maps and a choice of signs appearing in the relations:
\begin{MainThm}[Theorem \ref{thm:gen and rel}]
    Let $c,d\in \Z$ and let $(\mathcal{C},\otimes, \mathbbm{1},\tau)$ be a (strict) monoidal category enriched over graded abelian groups and denote $\mathcal{C}_n(X,Y)$ the abelian group of morphisms from $X$ to $Y$ of degree $n$.\\
    The data $(A,\mu,\eta,\nu,\varepsilon)$ where
    \begin{align*}
        A&\in \mathcal{C},\\
        \mu&\in \mathcal{C}_{c}(A\otimes A, A), & \eta&\in \mathcal{C}_{-c}(\mathbbm{1},A),\\
        \nu&\in \mathcal{C}_d(A,A\otimes A), & \varepsilon&\in \mathcal{C}_{-d}(A,\mathbbm{1})
    \end{align*}
    a monoidal functor out of $\pGrCob_{c,d}$ into $\mathcal{C}$ uniquely up to isomorphism defines if and only if the following relations are satisfied:
    \begin{enumerate}[(i)]
        \item graded associativity $\mu\circ (\mu \otimes \id)=(-1)^c\mu\circ (\id \otimes \mu)$;
        \item graded unitality $(-1)^c\mu\circ(\eta \otimes \id)=(-1)^{\frac{c(c-1)}{2}}\id=\mu\circ (\id \otimes \eta)$;
        \item graded coassociativity $(\id \otimes \nu)\circ \nu=(-1)^d(\nu \otimes \id)\circ \nu$;
        \item graded counitality $(-1)^d(\id \otimes \varepsilon)\circ \nu=(-1)^{\frac{d(d-1)}{2}}\id =(\varepsilon\otimes \id)\circ \nu$;
        \item graded Frobenius relation $(\mu \otimes \id)\circ (\id\otimes \nu)=(-1)^{cd}\nu\circ \mu=(\id \otimes \mu)\circ (\nu \otimes \id)$.
    \end{enumerate}
    If $\mathcal{C}$ is also symmetric, the data $(A,\mu,\eta,\nu,\varepsilon)$ defines a symmetric monoidal functor out of $\GrCob_{c,d}$ into $\mathcal{C}$ uniquely up to isomorphism if and only if the maps satisfy (i)-(v) and
    \begin{itemize}
        \item[\textit{(vi)}] graded commutativity $\mu \circ \tau=(-1)^c \mu$
    \end{itemize}
    and it defines a symmetric monoidal functor out of $\fatGrCob_{c,d}$ into $\mathcal{C}$ uniquely up to isomorphism if and only if the maps satisfy (i)-(v) and
    \begin{itemize}
        \item[\textit{(vi')}] graded symmetry $\varepsilon\circ \mu \circ \tau=(-1)^c\varepsilon\circ \mu$.
    \end{itemize}
\end{MainThm}
These signs are \textit{not} canonical. They depend on choices of generators $\omega(\multi)\in \bolddet_{c,d}(\multi)$ and $\omega(\comulti)\in \bolddet_{c,d}(\comulti)$. We explain in Remark \ref{rem:choice of orien} our choices. However, no choice can make the signs disappear.\\
We remarked above that our signs are stable under suspending algebras. This means that our choices of signs are such that if $(A,\mu,\eta,\nu,\varepsilon)$ satisfies the relations of $(c,d)$-graded Frobenius algebra, then $(\Sigma A,\Sigma \mu,\Sigma \eta,\Sigma \nu,\Sigma\varepsilon)$ satisfies the relations of a $(c-1,d+1)$-graded Frobenius algebra. This is proved in Proposition \ref{prop:shift maps}.\\
Our signs almost coincide with the signs found in the definition of a biunital coFrobenius bialgebra as defined in \cite{cieliebak2022cofrobenius}. We compare this in Remark \ref{rem:alex kai}.

\subsection{Examples}
In the literature, graded Frobenius algebras appear with different characterisations. Some have a more PROP-theoretic description and some have a description in terms of maps and relations. Our main theorem lets us relate the two types of description. In section \ref{sec:exmp}, we discuss a selection of the examples found in the literature.\\
The first example we discuss is the cohomology ring of an oriented manifold $M$ of dimension $d$ (see Example \ref{exmp:coh}). The intersection coproduct defines a degree $d$ comultiplication on cohomology. In particular, we discuss different sign conventions for the intersection coproduct found in the literature.\\
A second example, we discuss, is the Hochschild homology of a graded Frobenius algebra $A$ (see Example \ref{exmp:hochschild}). Wahl and Westerland give a description of a Frobenius algebra on $\overline{HH}_*(A)$ twisted by $\bolddet(S,\partial_{in})^{\otimes d}$. As discussed above, the twisting by the relative homology of the surface gives a $(d,d)$-graded Frobenius algebra.\\
In a similar vein, we discuss the homology of the loop space with the Chas-Sullivan product and the trivial coproduct. Godin gives operations which are twisted by $\bolddet(S,\partial_{in})^{\otimes d}$ in \cite{godin2007higher}. By our Proposition \ref{prop:1closed=2open twists}, this corresponds to a twisting by $\bolddet(G,\partial_{in})^{\otimes d}\otimes \bolddet(G,\partial_{out})^{\otimes d}$ (see Example \ref{exmp:LM w Bad coprod}).\\
Finally, we discuss the relation between the Chas-Sullivan product and the Goresky-Hingston coproduct in Rabinowitz loop homology as in the work by Cieliebak, Hingston and Oancea in \cite{cieliebak2020poincar,cieliebak2022cofrobenius,cieliebak2022reduced,cieliebak2023loop,cieliebak2024rabinowitz}.\\
We discuss the language of Tate vector spaces as presented in \cite{cieliebak2024rabinowitz}. Concretely for $A$ a graded Tate vector space, we show that the endomorphisms $\End_{A}(X,Y)=\Hom(A^{\starotimes X}, A^{\shriekotimes Y})$ have a canonical dioperad structure in Proposition \ref{prop:tate dioperad}.\\
On the other hand, we define the canonical subdioperad of $\GrCob_{c,d}$ in Definition \ref{def:subdiop}. This dioperad has the same operations as the full PROP in the case of a field of characteristic $\neq 2$. This is a consequence of our description of the subdioperad in Proposition \ref{prop:subdioperad}. We can thus describe the Rabinowitz loop homology as an algebra over the subdioperad in Example \ref{exmp:LM w good coprod}.

\subsection*{Organisation of the Paper.}
In Section \ref{sec:cob and graphs}, we describe a category of graph cobordisms modelling the category of 2-dimensional cobordisms. In Section \ref{sec:FrAlg}, we recall the classical relations of Frobenius algebras and show that they force algebras to be trivial in the odd dimensional graded case. In Section \ref{sec:gfa PROPs}, we construct the categories $\GrCob_{c,d}$, $\fatGrCob_{c,d}$ and $\pGrCob_{c,d}$ and study graded Frobenius algebras in terms of PROPs. In section \ref{sec:gfa maps and rel}, we state the main theorem and use it to study graded Frobenius algebras in terms of maps and relations. In Section \ref{sec:exmp}, we review and relate examples from the literature of graded Frobenius algebras. In Section \ref{sec:proof}, we set up machinery that helps us to prove the main theorem.

\subsection*{Conventions.}

Throughout this paper, we use the convention that maps act from the left. This in particular implies that for maps of graded abelian $f\colon A\to A'$ and $g\colon B\to B'$ evaluated on homogeneous elements $a\in A$ and $b\in B$, we have the following signs
\begin{align*}
    (f\otimes g)(a\otimes b)=(-1)^{|g||a|}f(a)\otimes g(b).
\end{align*}
Moreover, that the canonical composition map on morphism spaces is order-reversing
\begin{align*}
    \Hom(Y,Z)\otimes \Hom(X,Y)&\to \Hom(X,Z),\\
    f\otimes g&\mapsto f\circ g.
\end{align*}
We define the suspension of a graded abelian $A$ by $\Sigma A$. It is defined by $(\Sigma A)_{i+1}=A_{i}$.\\
We denote by $\mathbbm{1}$ the unit of the category of graded abelian groups: $\Z$ concentrated in degree $0$. We denote the dual of a graded abelian group by $A^\lor$ defined as $A^\lor_k:=\Hom_\bullet(A_{-k},\mathbbm1)$.\\
The natural isomorphism from the dual of a tensor product to the product of the duals is order-reversing, we have
\begin{align*}
    (A\otimes B)^{\lor}&\cong B^{\lor}\otimes A^{\lor},\\
    (a\otimes b)^{\lor}&\mapsto b^{\lor}\otimes a^{\lor}.
\end{align*}
This can be seen as $b^{\lor} \otimes a^{\lor}$ evaluated on $a\otimes b$ gives 1. But $a^{\lor}\otimes b^{\lor}$ evaluated on $a\otimes b$ gives $(-1)^{|a||b|}$.

\subsection*{Acknowledgments.}

I thank my PhD advisor Nathalie Wahl for suggesting the topic of signs in Frobenius algebras, many interesting discussions and reading my various drafts. I also like to thank Jan Steinebrunner for helpful discussions about operads, PROPs and Frobenius algebras. I further like to thank Kai Cieliebak and Alexandru Oancea for sharing their draft of \cite{cieliebak2022cofrobenius} and their insights on graded Frobenius algebras and its relation to string topology. Furthermore, I thank Alex Takeda for explaining the relation of our conventions to pre-Calabi-Yau structures. I would also like the referee for their insightful comments which helped clarify and extend this article. Finally, I thank Isaac Moselle and Janine Roshardt for proof-reading and encouraging discussions. I was supported by the Danish National Research Foundation through the Copenhagen Centre for Geometry and Topology (DNRF151).

\tableofcontents

\section{Cobordisms and Graphs}\label{sec:cob and graphs}

In this section, we recall the definition of closed, open and planar 2D-cobordism categories. Moreover, we describe an equivalent category coming from a 2-category of graph cobordism.
\subsection{Cobordisms}
\begin{figure}[H]
    \centering
    \includegraphics[width=0.24\linewidth]{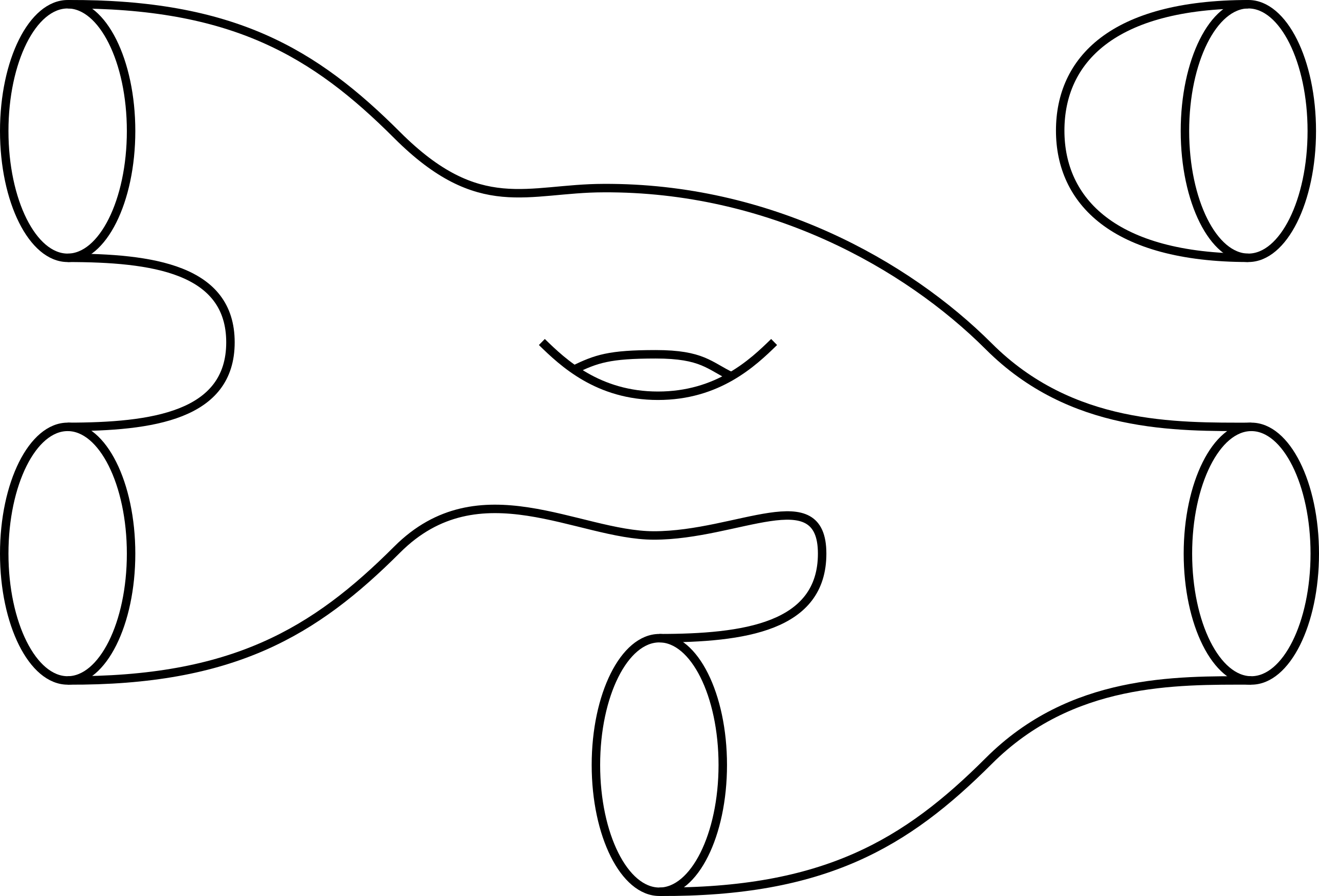}
    \caption{An example of a closed cobordism from $S^1\sqcup S^1 \sqcup S^1$ to $S^1\sqcup S^1$. We use the convention of drawing the incoming boundary to the left and the outgoing boundary to the right.}
    \label{fig:closed cob}
\end{figure}

\begin{figure}[H]
    \centering
    \includegraphics[width=0.6\linewidth]{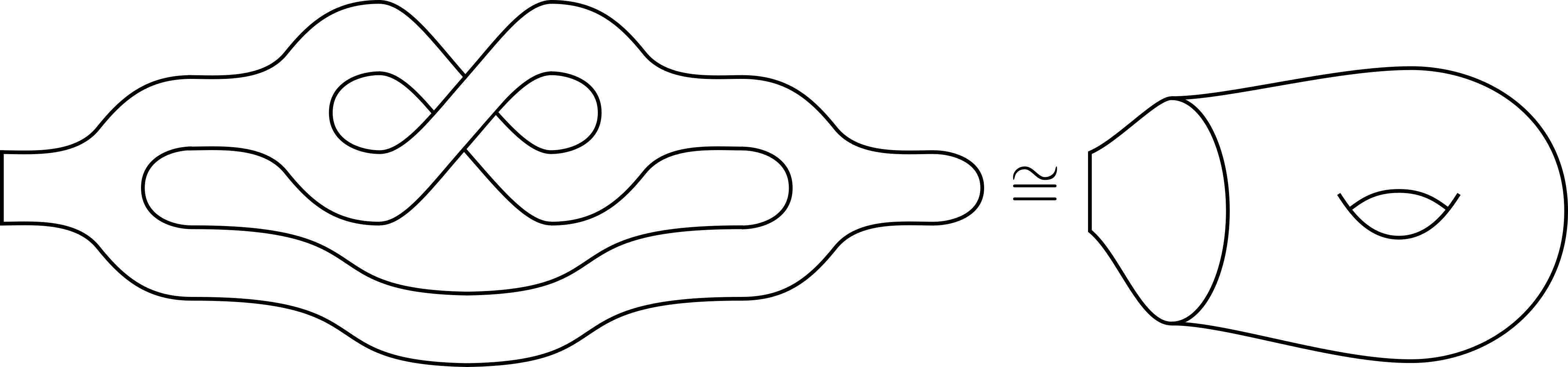}
    \caption{An example of an open cobordism from $[0,1]$ to $\varnothing$. The two surfaces are homeomorphic as they both have one boundary and Euler characteristic $-1$.}
    \label{fig:open cob}
\end{figure}

\begin{Def}
    The category of \textit{closed 2D-cobordisms} $\mathrm{Cob}_2^{\mathrm{closed}}$ has objects
    \begin{align*}
        \mathrm{Ob}(\mathrm{Cob}_2^{\mathrm{closed}})&:=\left\{\left. \coprod_{X} S^1\:\right| X \text{ is a finite set}\right\}.
    \end{align*}
    A morphism between $\coprod_X S^1$ and $\coprod_Y S^1$ is given by an equivalence class of pairs $(S,\varphi)$ where 
    \begin{itemize}
        \item $S$ is an oriented, compact topological 2-manifold possibly with boundary;
        \item $\varphi$ is a homeomorphism
        \begin{align*}
            \varphi\colon \coprod_X   S^1\sqcup\coprod_Y S^1 \to \partial S
        \end{align*}
    \end{itemize} 
    and $(S,\varphi)$ and $(S',\varphi')$ are equivalent if there exists a homeomorphisms $\psi:S\to S'$ such that $\psi\circ \varphi=\varphi'$. See Figure \ref{fig:closed cob}.
\end{Def}
\begin{Def}  
    The category of \textit{open 2D-cobordisms} $\mathrm{Cob}_2^{\mathrm{open}}$ has objects
    \begin{align*}
        \mathrm{Ob}(\mathrm{Cob}_2^{\mathrm{open}})&=\left\{\left.\coprod_{X} [0,1]\:\right|  X \text{ is a finite set}\right\}.
    \end{align*}
    A morphism between $\coprod_X [0,1]$ and $\coprod_Y [0,1]$ is given by an equivalence class of pairs $(S,\varphi)$ where 
    \begin{itemize}
        \item $S$ is an oriented, compact topological 2-manifold such that every connected component has non-empty boundary;
        \item $\varphi$ is an embedding
        \begin{align*}
            \varphi\colon \coprod_X [0,1]\sqcup \coprod_Y [0,1]\to \partial S.
        \end{align*}
    \end{itemize}
    The equivalence relation is given by $(S,\varphi)\sim (S',\varphi)$ if there exists a homeomorphisms $\psi:S\to S'$ such that $\psi\circ \varphi=\varphi'$. We call $\partial_{free} S=\partial S\setminus \im(\varphi)$ the free part of the boundary. See Figure \ref{fig:open cob}.
\end{Def}

\begin{figure}[H]
    \centering
    \includegraphics[width=0.4\linewidth]{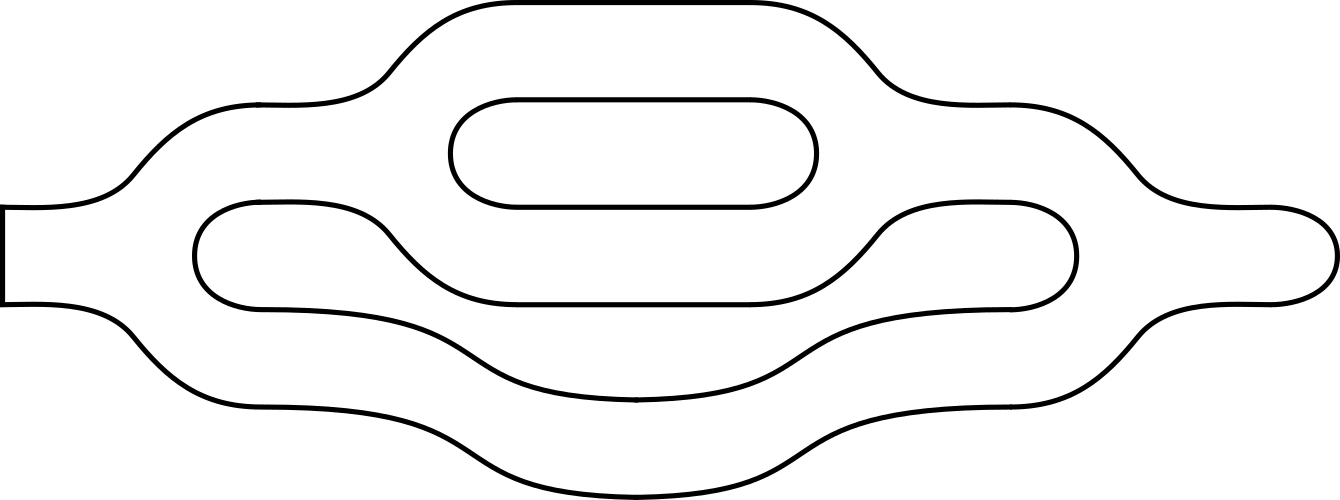}
    \caption{A planar cobordism between $[0,1]$ and $\varnothing$.}
    \label{fig:planar cob}
\end{figure}

We fix some orientation on $[0,1]\times \R$.
\begin{Def}  
    The category of \textit{planar 2D-cobordisms} $\mathrm{Cob}_2^{\mathrm{planar}}$ has objects
    \begin{align*}
        \mathrm{Ob}(\mathrm{Cob}_2^{\mathrm{planar}})&=\left\{\left.\coprod_{X} [0,1]\:\right|  X \text{ is a finite set with linear ordering}\right\}.
    \end{align*}
    A morphism between $\coprod_X [0,1]$ and $\coprod_Y [0,1]$ is given by an open cobordisms 
    \begin{align*}
        [(S,\varphi)]\in \Cob_2^{\mathrm{open}}\left(\coprod_X [0,1],\coprod_Y [0,1]\right)
    \end{align*}
    such that there exists an oriented embedding $\psi\colon S\hookrightarrow [0,1]\times \R$ satisfying that
    \begin{itemize}
        \item $\psi\circ \varphi|_{\coprod_X[0,1]}$ lands in $\{0\}\times \R$ and the image of the connected components is ordered (with respect to the canonical ordering on $\{0\}\times \R$) according to the linear ordering of $X$;  
        \item $\psi\circ \varphi|_{\coprod_Y[0,1]}$ lands in $\{1\}\times \R$ and the image of the connected components is ordered in $\{1\}\times \R$ according to the linear ordering of $Y$.  
    \end{itemize}
     See Figure \ref{fig:planar cob}.
\end{Def}

\subsection{Graphs}
We later give a description of those three categories using graphs. For this, we introduce what we mean by a graph.

\begin{Def}
    A \textit{graph} $G=(V,H,\sigma,s,L_{in},L_{out})$ consists of the following data:
    \begin{itemize}
        \item a finite set $V(G):=V$ called the \textit{vertices};
        \item a finite set $H(G):=H$ called the \textit{half-edges};
        \item an involution $\sigma\colon H\to H$ with no fixed points, the orbits $\{h,\sigma h\}$ are called the \textit{edges} and we denote the set of orbits $H_\sigma=:E=:E(G)$;
        \item a surjective map $s\colon H\to V$, where for a vertex $v\in V$ the cardinality $|s^{-1}(v)|$ is called its \textit{arity};
        \item subsets $L_{in}(G):=L_{in},L_{out}(G):=L_{out}\subseteq V(G)$ of the set of vertices of arity $\leq 1$ such that $L_{in}\cap L_{out}$ only contains vertices of arity $0$.
    \end{itemize} 
    A \textit{fat graph} is a graph together with a cyclic ordering of the half-edges in $s^{-1}(v)$ for all $v$.
\end{Def}
We call the vertices in $L_{in}$ and $L_{out}$ the \textit{external vertices}. The other vertices are called the \textit{internal vertices}. In our conventions, we allow for internal vertices of arity $\leq 1$.\\

A graph defines a 1-dimensional CW-complex, called the geometric realization, denoted $|G|$: the zero-cells are the vertices and the one-cells are the edges. Choosing an orientation on each one-cell corresponds to an ordering of the two half-edges making up the edge. Then the map $s$ describes how to glue the edges to the vertices.\\
We give some additional terminology which is useful, when talking about graphs.
\begin{Def}
    Let $G$ be a graph. 
    \begin{itemize}
        \item The graph $G$ is \textit{connected} if $|G|$ is connected.
        \item The graph $G$ is a \textit{forest} if $\pi_1(|G|)=1$ at every basepoint.
        %\item The graph $G$ is a \textit{tree} if $G$ is a connected forest.
        %\item The \textit{rank} $r(G)$ of the graph $G$ is the rank of the free group $H_1(|G|)$.
        \item An edge $e=\{h,\sigma h\}$ is called a \textit{tadpole} if $s(h)=s(\sigma h)$.
    \end{itemize}
\end{Def}

\begin{Def}
    A \textit{morphism of graphs} $f\colon G\to G'$ is a cellular map of the underlying CW-complexes $|G|\to|G'|$ which is on edges either the collapse or the identity and is a bijection on $L_{in}$ and $L_{out}$.
\end{Def}

\begin{Rem}
    Our definition of morphism of graphs is such that every morphism of graph is given by a sequence of edge collapses and a graph isomorphism. Importantly, the edges that are collapsed are not tadpoles and are not adjacent to an external vertex. 
\end{Rem}

\begin{figure}[H]
    \centering
    \includegraphics[width=0.8\linewidth]{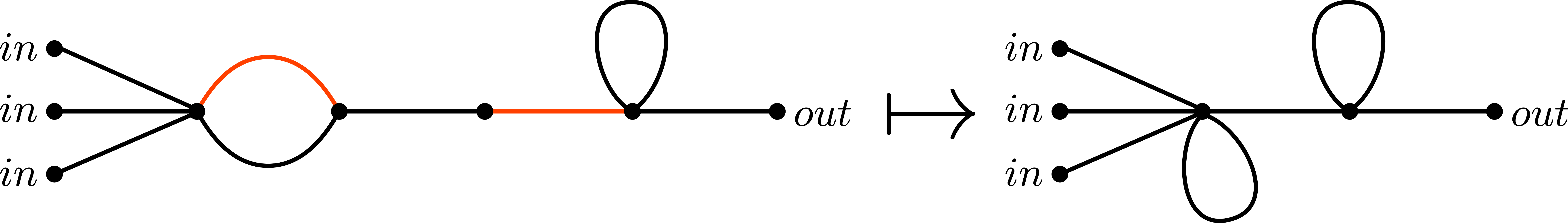}
    \caption{A morphism of graphs sends $L_{in}$ to $L_{in}$ and $L_{out}$ to $L_{out}$ and possibly collapses edges. Here the graph morphism collapses two edges.}
    \label{fig:graph morph}
\end{figure}

\begin{Def}
    A \textit{morphism of fat graphs} is a morphism of graphs $f\colon G\to G'$ that can be written as a sequence of edge collapses of edges
    \begin{align*}
        G=:G_0\to G/e_1=:G_1\to G_1/e_2=:G_2\to  \dots \to G_{k-1}/e_k=:G_k
    \end{align*}
    and a graph isomorphism
    \begin{align*}
        \varphi\colon G_k\to G'
    \end{align*}
    such that 
    \begin{itemize}
        \item if we have for $e_i=\{h_i,\sigma h_i\}$ the cyclic orderings $(h_i \; a_1  \dots a_m)=(a_1\dots a_m\; h_i)$ on $s^{-1}(s(h_i))$ and  $(\sigma h_i\; b_1\dots b_n)$ on $s^{-1}(s(\sigma h_i))$, then the vertex in $G_{i-1}/e_i$ corresponding to $e_i$ has cyclic ordering 
        \begin{align*}
            (a_1\dots a_m\; b_1\dots b_n);
        \end{align*}
        \item $\varphi$ commutes with the cyclic ordering.
    \end{itemize}
\end{Def}

\subsection{Graph Cobordisms}\label{subsec:graph cob}

We can now construct the 2-categories $\underline{\GrCob}$, $\underline{\fatGrCob}$ and $\underline{\pGrCob}$. This gives models for the cobordism categories.

\begin{Def}
    The \textit{graph cobordism category} $\GrCob$ is a 2-category, i.e.\ a category enriched over 1-categories with
    \begin{align*}
        \mathrm{Ob}(\underline{\GrCob})&=\{\text{finite sets}\},
    \end{align*}
    and the morphism categories $\underline{\GrCob}(X,Y)$ are given by the category of graphs $G$ with labellings $L_{in}(G)\cong X$ and $L_{out}(G)\cong Y$ and 2-morphisms graph morphisms which respect the labellings.\\
    The \textit{fat graph cobordism category} is a 2-category $\underline{\fatGrCob}$ with
    \begin{align*}
        \mathrm{Ob}(\underline{\fatGrCob})&=\{\text{finite sets}\},
    \end{align*}
    and the morphism categories $\underline{\fatGrCob}(X,Y)$ are given by the analogous category of fat graphs with labellings and morphisms are fat graph morphisms which respect the labellings.
\end{Def}

\begin{Exmp}\label{exmp:G_f}
    Let $f\colon X\cong Y$ be an isomorphism of finite sets. We define the graph $G_f$ 
    \begin{align*}
        V(G_f)&:=X, & H(G_f)&:=\varnothing, & L_{in}(G_f)&:=V(G_f)=:L_{out}(G_f).
    \end{align*}
    This has labellings $ L_{in}(G)\cong X$ via the identity and $L_{out}(G)\cong Y$ via $f$. Moreover, $G_f$ consists only of vertices of valence zero and thus has a canonical fat structure.
\end{Exmp}

\begin{figure}[H]
    \centering
    \includegraphics[width=1\linewidth]{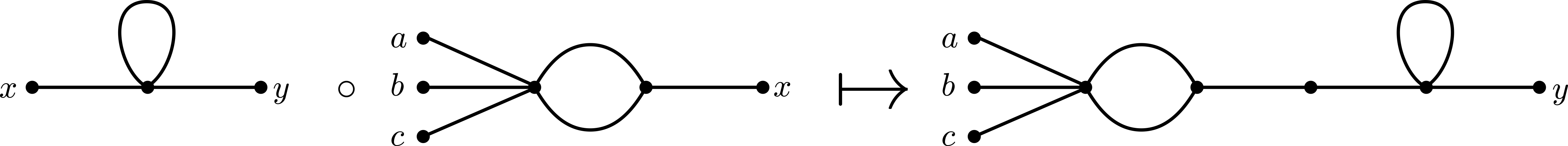}
    \caption{Composition of a graph cobordism between $\{a,b,c\}$ and $\{x\}$ and a graph cobordism between $\{x\}$ and $\{y\}$.}
    \label{fig:graph comp}
\end{figure}

These categories are well-defined as the composition along objects is given as follows: let $X,Y,Z$ be finite. We define the composition as
\begin{align*}
     \underline{\GrCob}(Y,Z)\times \underline{\GrCob}(X,Y) &\to \underline{\GrCob}(X, Z)
\end{align*}
which sends $(G',G)$ to the graph $G'\circ G:=G'\cup_Y G$. This is defined as the graph obtained from the disjoint union of $G$ and $G'$ where we associate the elements in $L_{out}(G)$ with $L_{in}(G')$ via the labelling by $Y$. In more visual terms, we glue the outgoing legs of $G$ to the incoming legs of $G'$ (as in Figure \ref{fig:graph comp}).\\
We note that this defines a functor as morphisms in $\underline{\GrCob}(X,Y)$ and $ \underline{\GrCob}(Y,Z)$ can be glued together in a similar manner.\\
The unit $\id_X\in \underline{\GrCob}(X,X)$ is given by the graph $G_\id$ as defined in Example \ref{exmp:G_f}.\\
Analogously, we can define composition maps
\begin{align*}
    \underline{\fatGrCob}(Y,Z)\times \underline{\fatGrCob}(X,Y) &\to \underline{\fatGrCob}(X,Z)
\end{align*}
for fat graphs.\\
This works as the newly combined vertices coming from $L_{out}(G)\cong L_{in}(G')$ have arity $\leq2$ and thus there is only one cyclic ordering of the half-edges at these vertices. The fat graph structure of $G'\circ G$ is thus well-defined.\\
Moreover, the unit $G_\id$ is also a fat graph.

\begin{Lem}\label{lem:GrCob is symmon}
    Disjoint unions define a symmetric monoidal structure on $\underline{\GrCob}$ and $\underline{\fatGrCob}$.
\end{Lem}
\begin{proof}
    We note that labellings 
    \begin{align*}
        L_{in}(G)&\cong X,& L_{out}(G)&\cong Y, & L_{in}(G')&\cong X', & L_{out}(G')&\cong Y'
    \end{align*}
    on (fat) graphs $G,G'$ induce labellings
    \begin{align*}
        L_{in}(G\sqcup G')&\cong X\sqcup X', & L_{out}(G\sqcup G')&\cong Y\sqcup Y'.
    \end{align*}
    This defines a coherent symmetric monoidal structure because of the following observation: the category $\mathrm{Fin}^{\cong}$ of finite sets with bijections is a canonical subcategory of $\underline{\GrCob}$ and $\underline{\fatGrCob}$ via $f\mapsto G_f$ as in Example \ref{exmp:G_f}.\\
    It is then a direct calculation that the unit, associator and twist map of the symmetric monoidal structure on $\mathrm{Fin}^{\cong}$ thus also induce unit, associator and twist map on $\underline{\GrCob}$ and $\underline{\fatGrCob}$.
\end{proof}
In \cite{chiche2014th}[Section 1.3], it is described that the inclusion of $1$-categories into $2$-categories has a left adjoint $\tau_\iota$ given on a 2-category $\underline{\mathcal{C}}$ by
\begin{align*}
    \mathrm{Ob}(\tau_{\iota}(\underline{\mathcal{C}}))&=\mathrm{Ob}(\underline{\mathcal{C}}),\\
    \tau_\iota(\underline{\mathcal{C}})(X,Y)&=\pi_0 (\underline{\mathcal{C}}(X,Y)).
\end{align*}

\begin{Def}\label{def:ass 1-cat}
    We define $\GrCob:=\tau_{\iota}(\underline{\GrCob})$, $\fatGrCob:=\tau_{\iota}(\underline{\fatGrCob})$ the associated 1-categories.
\end{Def}

An analogous argument to Lemma  \ref{lem:GrCob is symmon} shows that disjoint union induces a symmetric monoidal structure on $\GrCob$ and $\fatGrCob$.\\

\begin{figure}[H]
    \centering
    \includegraphics[width=0.5\linewidth]{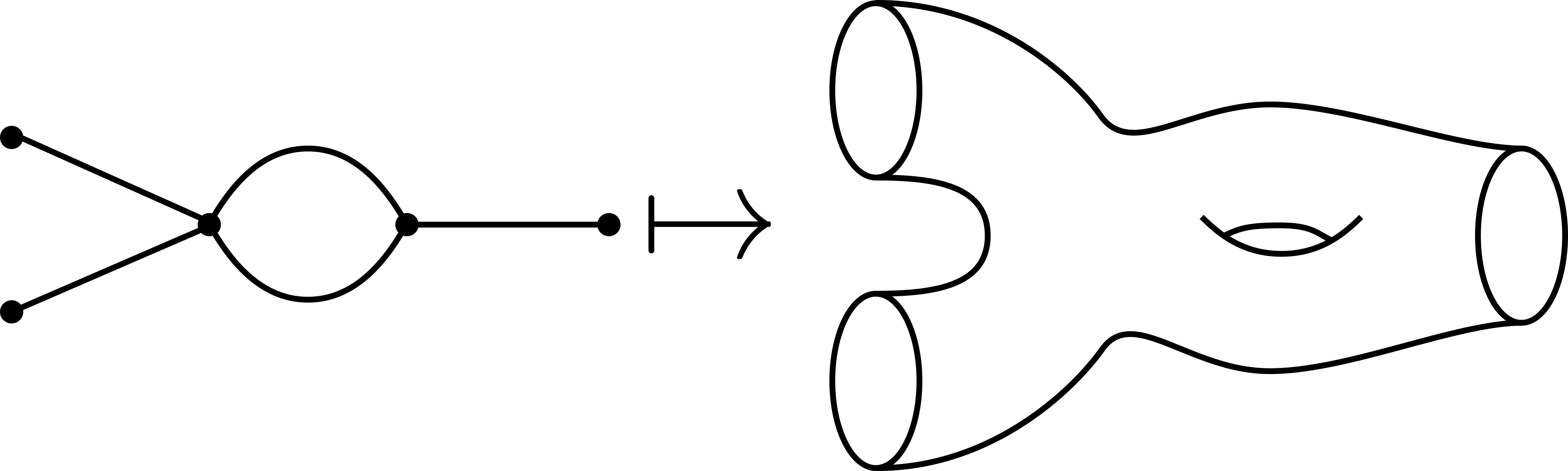}
    \caption{The closed cobordism $\Cob(G)$ corresponding to a graph $G$.}
    \label{fig:cob(G)}
\end{figure}

\begin{figure}[H]
    \centering
    \includegraphics[width=0.6\linewidth]{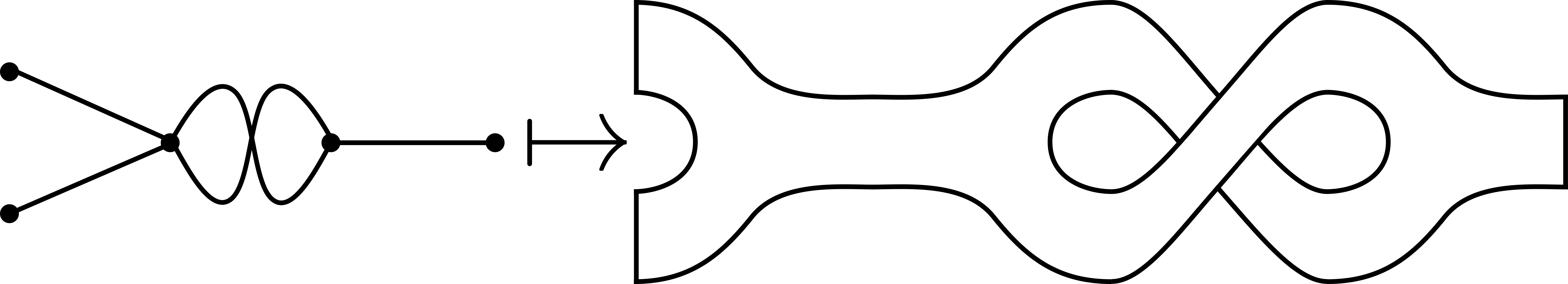}
    \caption{The open cobordism $\fatCob(G)$ corresponding to a fat graph $G$.}
    \label{fig:fatcob(G)}
\end{figure}

We can now construct equivalences of categories $\GrCob\to \Cob^{\mathrm{closed}}_2$ and $\fatGrCob\to \Cob^{\mathrm{open}}_2$.\\
A morphism in $\GrCob$ (or $\fatGrCob$) is given by an equivalence class $[G]$ of (fat) graphs up to zig-zags of (fat) graph isomorphisms.\\
Let $G\in\underline{\GrCob}(X,Y)$ be a graph cobordism. We construct the corresponding closed cobordism $\Cob(G)$ as follows.\\
An internal vertex of arity $0$ contributes a copy of $S^2$ in $\Cob(G)$. A vertex of arity $0$ which is only incoming or outgoing contributes a copy of $D^2$ in $\Cob(G)$ with its boundary labelled by the corresponding element in $X$ or $Y$. A vertex of arity $0$ in $L_{in}\cap L_{out}$ contributes a cylinder $S^1\times [0,1]$ with its two boundary components labelled by the corresponding elements in $X$ and $Y$. For a connected component of $G$ consisting of more than just one vertex, each edge gives a copy of the cylinder $S^1\times [0,1]$, each external vertex of arity $1$ is labelled by an element in $X$ or $Y$. This gives a labelled copy of $S^1$ which we glue to the corresponding cylinder. Finally, we can glue the edges at internal vertices to copies of $S^2\setminus(\sqcup_{s^{-1}(v)} D^2)$ (see Figure \ref{fig:cob(G)}).\\
The homeomorphism type of this construction is invariant under edge collapses and graph isomorphisms. Therefore $\mathrm{Cob}(G)$ only depends on the class that $G$ represents in $\pi_0(\underline{\GrCob}(X,Y))$. We thus get a well defined map 
\begin{align*}
    \GrCob(X,Y)=\pi_0(\underline{\GrCob}(X, Y))&\to \mathrm{Cob}_2^{\mathrm{closed}}(\coprod_{X} S^1,\coprod_Y S^1),\\
    [G]& \mapsto \Cob(G).
\end{align*}
Similarly we can construct for a fat graph cobordism between $X$ and $Y$ an open cobordism $\fatCob(G)$ between $\coprod_X [0,1]$ and $\coprod_Y [0,1]$. It is given by the following construction.\\ 
Vertices $v$ of arity 0 contribute a copy of $D^2$ with an embedding of copies of the $[0,1]$ labelled by the elements in $X,Y$ labelling $v$. For a connected component of $G$ consisting of more than just one vertex, each edge can be thickened to $[0,1]\times [0,1]$ and each external vertex corresponds to a closed interval. This gives us copies of the interval $[0,1]$ labelled by $X$ and $Y$. Finally, the cyclic ordering at each internal vertex gives a way to glue the thick edges to a closed disc $D^2$ (see Figure \ref{fig:fatcob(G)}).\\
The homeomorphism type $\fatCob(G)$ only depends on the equivalence class $[G]\in \pi_0(\underline{\fatGrCob}(X, Y))$, because the construction is invariant under fat edge collapses and graph isomorphisms.\\
We thus have defined a map:
\begin{align*}
    \fatGrCob(X,Y)=\pi_0(\underline{\fatGrCob}(X, Y))&\to {\mathrm{Cob}_2^{\mathrm{open}}}\left(\coprod_{X} [0,1],\coprod_Y [0,1]\right),\\
    [G]& \mapsto \fatCob(G).
\end{align*}

\begin{Prop}\label{prop:graph cob is a model}
    There are equivalences of categories
    \begin{align*}
        \GrCob & \to\mathrm{Cob}_2^{\mathrm{closed}},\\
        X &\mapsto \coprod_{X} S^1,\\
        [G]& \mapsto \mathrm{Cob}(G)
    \end{align*}
    and
    \begin{align*}
        \fatGrCob & \to\mathrm{Cob}_2^{\mathrm{open}},\\
        X &\mapsto \coprod_{X} [0,1],\\
        [G]& \mapsto \fatCob(G).
    \end{align*}
\end{Prop}
\begin{proof}
    We first prove that this map commutes with composition: composition of graphs $G\in \GrCob(X,Y)$ and $G'\in\GrCob(Y,Z)$ is given by gluing the outgoing legs of $G$ to the ingoing legs of $G'$ according to the $Y$ labelling. We do a case distinction for the vertices we glue together.
    \begin{itemize}
        \item If one of them is an arity zero vertex which is both incoming and outgoing, gluing the vertex does not change (the isomorphism type of) $G$. Similarly after applying $\Cob$ this corresponds to gluing a cylinder $S^1\times [0,1]$ to $S^1$. This does not change the homeomorphism type of $\Cob(G)$.
        \item If one of the vertices is an arity zero vertex $v$ which is either incoming or outgoing, gluing $v$ to $G$ corresponds to forgetting the label at the glued vertex. After applying $\Cob$, this gives $S^2\setminus D^2$ glued along $S^1$. This is homeomorphic to gluing $\Cob(v)=D^2\cong S^2\setminus D^2$ along $S^1$.
        \item If both glued vertices are of arity $1$, gluing them gives a valence 2 vertex. Under our construction a valence 2 gives $S^2\setminus (D^2\sqcup D^2)$ in $\Cob(G' \circ G)$ to which we glue two cylinders. This is homeomorphic to just gluing the two cylinders together along $S^1$.
    \end{itemize}
    A similar argument to the closed case shows that $\fatCob$ commutes with composition. For this we note that $[0,1]\times [0,1]$ to a disk $D^2$ along the boundary is homeomorphic to $[0,1]\times [0,1]$.\\
    The unit in the graph categories are given by $G_{\id_X}$ as in Example \ref{exmp:G_f}. We have
    \begin{align*}
        \Cob\left( G_{\id_X}\right)&=\coprod_X S^1\times [0,1]=\id_X,\\
        \fatCob\left(G_{\id_X}\right)&=\coprod_X [0,1]\times [0,1]=\id_X
    \end{align*}
    and thus $\Cob$ and $\fatCob$ are functors.\\
    Moreover, the functors are surjective on objects.\\
    It thus remains to show that $\Cob$ and $\fatCob$ are bijective on the set of morphisms. For this, we note that both functors are monoidal because the monoidal structure on the source and the target is given by disjoint union which commutes with our construction.\\
    Kock shows in \cite{kock2004frobenius} that $\Cob_2^{\mathrm{closed}}$ is generated by 
    \begin{align*}
        \closedmulti,\closedunit,\closedcomulti, \closedcounit
    \end{align*}
    as a symmetric monoidal category. These generators are in the image of $\Cob$ which shows fullness of $\Cob$.\\
    Similarly, Lauda and Pfeiffer show in \cite{lauda2008open} that $\Cob_2^{\mathrm{open}}$ are generated by
    \begin{align*}
        \openmulti,\openunit,\opencomulti,\opencounit
    \end{align*}
    which are in the image. This shows fullness of $\fatCob$.\\
    For faithfulness of $\Cob$, we assume that $\Cob(G)=\Cob(G')$. Without loss of generality, we may assume that both graphs are connected and get sent to a surface of genus $g$. We know that $G$ and $G'$ have the same number of incoming and outgoing legs. We choose trees in $T\subseteq G$ and $T\subseteq G'$ that contain exactly all internal vertices. Collapsing all edges in $T$ and $T'$ gives graphs $G/T$ and $G'/T'$ with one internal vertex. Because $[G]=[G/T]$, we have $\Cob(G/T)=\Cob(G)$ is a surface of genus $g$. In particular, $G/T$ has $g$ tadpoles. Similarly, we note that $G'/T'$ has $g$ tadpoles. This shows $G/T\cong G'/T'$ and thus $[G]=[G/T]=[G'/T']=[G']$ shows faithfulness of $\Cob$.\\
    For faithfulness of $\fatCob$, we use the category of open-closed fat graphs $\mathpzc{Fat}^{\mathpzc{oc}}$ as defined in \cite{santander2015comparing}. Egas shows in \cite[Theorem A]{santander2015comparing} 
    \begin{align*}
        \mathpzc{Fat}^{\mathpzc{oc}}\simeq \coprod_S B\mathrm{Mod}(S)
    \end{align*}
    where the disjoint union runs over all topological types of open-closed cobordisms in which each connected component has at least one boundary component which is not free. Restricting to the open part and taking $\pi_0$ shows that if $\fatCob(G)=\fatCob(G')$ implies $[G]= [G']$ for $G$ and $G'$ where all connected components have at least one external vertex.\\
    For general $G$ and $G'$ with $\fatCob(G)=\fatCob(G')$, we can add an edge and an external vertex to each connected component and obtain graphs $G_+$ and $G'_+$. We have $\fatCob(G_+)=\fatCob(G_+')$ by the above argument. We then have a zig-zag of fat graph morphisms connecting $G_+$ and $G_+'$. This zig-zag sends the added external vertices to each other. Thus we obtain by restriction a zig-zag connecting $G$ and $G'$ and thus $[G]=[G']$. This shows faithfulness of $\fatCob$ and concludes the proof.
\end{proof}

We can now define the planar graph cobordism category.

\begin{Def}
    The \textit{planar graph cobordism category} is a 2-category $\underline{\pGrCob}$ with
    \begin{align*}
        \mathrm{Ob}(\underline{\pGrCob})&=\{\text{finite sets with linear ordering}\},
    \end{align*}
    and the morphism categories $\underline{\pGrCob}(X,Y)$ are given by the full subcategory of  $\underline{\fatGrCob}(X,Y)$ spanned by all fat graphs $G$ such that $\fatCob(G)$ is a planar cobordism in $\Cob_2^{\mathrm{planar}}(X,Y)$.
\end{Def}

We denote the associated $1$-category $\pGrCob:=\tau_{\iota}(\underline{\pGrCob})$ with 
\begin{align*}
    \pGrCob(X,Y)=\pi_0(\underline{\pGrCob}(X,Y)).
\end{align*}
The monoidal structure of $\mathrm{Fin}_\leq$, the category of finite ordered sets with order preserving maps, induces a monoidal structure on $\underline{\pGrCob}$ and $\pGrCob$. The argument is analogous to Lemma \ref{lem:GrCob is symmon} as for any isomorphism of finite ordered sets $f\colon X\cong Y$, the graph $G_f$ as in Example \ref{exmp:G_f} is a planar graph.\\
However, the monoidal structure is not symmetric. This is because the monoidal structure on $\mathrm{Fin}_\leq$ is also not symmetric.\\

We have a forgetful functor $\pGrCob\to \fatGrCob$. The category $\pGrCob$ is defined such that 
\begin{center}
    \begin{tikzcd}
        \pGrCob\ar[r] &\fatGrCob \ar[r, "\fatCob"]& \Cob_2^{\mathrm{open}}
    \end{tikzcd}
\end{center}
factors through $\Cob_2^{\mathrm{planar}}\to \Cob_2^{\mathrm{open}}$. Moreover, the following Corollary is immediate from Proposition \ref{prop:graph cob is a model} and the definitions of $\pGrCob$ and $\Cob_2^{\mathrm{planar}}$.
\begin{Cor}\label{cor:planar is model}
    The functor
    \begin{align*}
        \pGrCob &\to \Cob_2^{\mathrm{planar}},\\
        X&\mapsto \coprod_X [0,1],\\
        [G]&\mapsto \fatGrCob(X)
    \end{align*}
    is an equivalence of categories.
\end{Cor}

\section{Frobenius Algebras}\label{sec:FrAlg}

In this section, we recall the definition of Frobenius algebras. We give definitions both in terms of maps and generators and, on the other hand, as functors out of cobordism categories. Then in \ref{subsec:graded context}, we apply those definitions to the graded context and show in Proposition \ref{prop:trivial} that in the odd dimensional case, there exist only trivial Frobenius algebras.

\begin{Def}
    Let $(\mathcal{C},\otimes,\mathbbm{1}, \tau)$ be a (strict) monoidal category.\\
    A \textit{Frobenius algebra} $(A,\mu,\eta,\nu,\varepsilon)$ in $\mathcal{C}$ is an object $A$ in $\mathcal{C}$ together with maps
    \begin{align*}
        \mu\colon&A\otimes A\to A, & \eta\colon&\mathbbm{1}\to A,\\
        \nu\colon&A\to A\otimes A, & \varepsilon\colon&A\to \mathbbm{1}
    \end{align*}
    satisfying the following relations:
    \begin{enumerate}[(i)]
        \item associativity $\mu\circ (\mu \otimes \id)=\mu\circ (\id \otimes \mu)$;
        \item unitality $\mu\circ(\eta \otimes \id)=\id_A=\mu\circ (\id \otimes \eta)$;
        \item coassociativity $(\id \otimes \nu)\circ \nu=(\nu \otimes \id)\circ \nu$;
        \item counitality $(\id \otimes \varepsilon)\circ \nu=\id =(\varepsilon\otimes \id)\circ \nu$;
        \item Frobenius relation $(\mu \otimes \id)\circ (\id\otimes \nu)=\nu\circ \mu=(\id \otimes \mu)\circ (\nu \otimes \id)$.
    \end{enumerate}
    If $\mathcal{C}$ is also symmetric with twist map $\tau$, a \textit{commutative Frobenius algebra} $(A,\mu,\eta,\nu,\varepsilon)$ in $\mathcal{C}$ is an object $A$ in $\mathcal{C}$ together with maps as above satisfying the relations (i)-(v) and
    \begin{itemize}
        \item[(vi)] commutativity $\varepsilon\circ \mu \circ \tau=\mu$
    \end{itemize}
    and a \textit{symmetric Frobenius algebra} $(A,\mu,\eta,\nu,\varepsilon)$ in $\mathcal{C}$ is an object $A$ in $\mathcal{C}$ together with maps as above satisfying the relations (i)-(v) and
    \begin{itemize}
        \item[(vi')] symmetry $\varepsilon\circ \mu \circ \tau=\varepsilon \circ \mu$.
    \end{itemize}
    A \textit{morphism of Frobenius algebras} $(A,\mu,\eta,\nu,\varepsilon)\to (A',\mu',\eta',\nu',\varepsilon')$ is a morphism $\varphi\colon A\to A'$ in $\mathcal{C}$ such that 
    \begin{itemize}
        \item $\mu\circ (\varphi\otimes \varphi)=\varphi\circ \mu'$;
        \item $\varphi\circ \eta=\eta'$;
        \item $\nu \circ \varphi= (\varphi\otimes \varphi)\circ \nu'$;
        \item $\varepsilon=\varepsilon'\circ \varphi$.
    \end{itemize}
    Denote by $\mathrm{Frob}(\mathcal{C})$ the category of Frobenius algebras in $\mathcal{C}$ with morphisms of Frobenius algebras, by $\mathrm{ComFrob}(\mathcal{C})$ the full subcategory spanned by commutative Frobenius algebras and by $\mathrm{SymmFrob}(\mathcal{C})$ the full subcategory spanned by symmetric Frobenius algebras.
\end{Def}

It is a well-known fact that an equivalent description of commutative Frobenius algebras is as symmetric monoidal functor out of a closed 2D-cobordism category (see \cite{abrams1996two,kock2004frobenius}). Similarly, a symmetric Frobenius algebra can be characterised as a symmetric monoidal functor out of a open 2D-cobordism category and a Frobenius algebra can be characterised as a monoidal functor out of a planar 2D-cobordism category  (see \cite{lauda2008open}).

\begin{Not}
    For (strict) symmetric monoidal categories $\mathcal{C},\mathcal{D}$, we denote $\Fun^\otimes(\mathcal{C},\mathcal{D})$ for the category of (strict) symmetric monoidal functors from $\mathcal{C}$ to $\mathcal{D}$ and symmetric monoidal natural transformations between them.\\
    For (strict) monoidal categories $\mathcal{C},\mathcal{D}$, we denote $\Fun^{\mathrm{mon}}(\mathcal{C},\mathcal{D})$ for the category of (strict) monoidal functors from $\mathcal{C}$ to $\mathcal{D}$ and monoidal natural transformations between them.\\
\end{Not}

\begin{Prop}[{\cite[Corollary 4.5, Corollary 4.6, Corollary 4.7.]{lauda2008open}}]
    There are natural equivalences of categories
    \begin{align*}
        \Fun^\otimes(\mathrm{Cob}_2^{\mathrm{closed}},\mathcal{C})&\to \mathrm{ComFrob}(\mathcal{C}),\\
        F &\mapsto (F\left(\closedobject\right), F\left(\closedmulti\right),F\left(\closedunit\right),F\left(\closedcomulti\right), F\left(\closedcounit\right)),
    \end{align*} 
    \begin{align*}
        \Fun^\otimes(\mathrm{Cob}_2^{\mathrm{open}},\mathcal{C})&\to \mathrm{SymFrob}(\mathcal{C}),\\
        F &\mapsto (F\left(\openobject\right), F\left(\openmulti\right),F\left(\openunit\right),F\left(\opencomulti\right), F\left(\opencounit\right))
    \end{align*}
    and
    \begin{align*}
        \Fun^\mathrm{mon}(\mathrm{Cob}_2^{\mathrm{planar}},\mathcal{C})&\to \mathrm{Frob}(\mathcal{C}),\\
        F &\mapsto (F\left(\openobject\right), F\left(\openmulti\right),F\left(\openunit\right),F\left(\opencomulti\right), F\left(\opencounit\right)).
    \end{align*}
    In particular, a natural transformation gives a morphism of Frobenius algebras.
\end{Prop}

We can also apply Proposition \ref{prop:graph cob is a model} and Corollary \ref{cor:planar is model} to get a characterisation in terms of graph cobordisms.

\begin{Cor}
    Let $\mathcal{C}$ be a monoidal category. There exists an equivalence of categories
    \begin{align*}
        \Fun^{\mathrm{mon}}(\pGrCob,\mathcal{C})&\to \mathrm{Frob}(\mathcal{C}),\\
        F &\mapsto (F\left(*\right), F\left(\multi\right),F\left(\unit\right),F\left(\comulti\right), F\left(\counit\right)).
    \end{align*}
    If $\mathcal{C}$ is also symmetric, there are natural equivalences of categories
    \begin{align*}
        \Fun^\otimes(\GrCob,\mathcal{C})&\to \mathrm{ComFrob}(\mathcal{C}),\\
        F &\mapsto (F\left(*\right), F\left(\multi\right),F\left(\unit\right),F\left(\comulti\right), F\left(\counit\right))
    \end{align*}
    and 
    \begin{align*}
        \Fun^\otimes(\fatGrCob,\mathcal{C})&\to \mathrm{SymFrob}(\mathcal{C}),\\
        F &\mapsto (F\left(*\right), F\left(\multi\right),F\left(\unit\right),F\left(\comulti\right), F\left(\counit\right)).
    \end{align*}
\end{Cor}

Here $\multi$ denotes the following graph
\begin{figure}[H]
    \centering
    \includegraphics[width=0.2\linewidth]{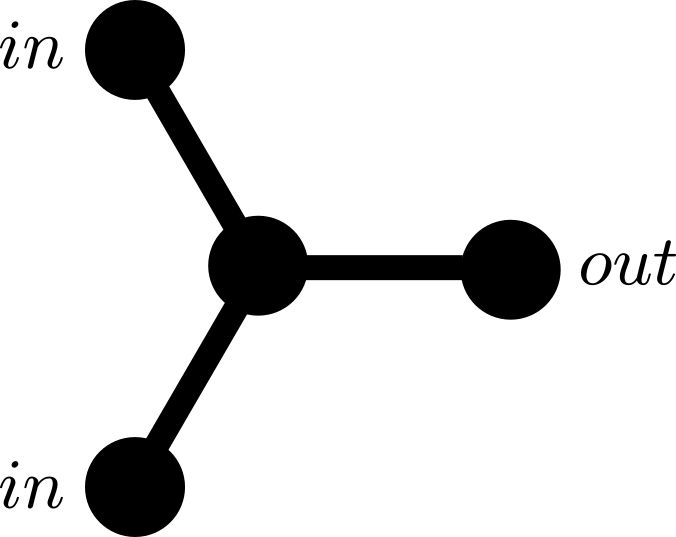}
\end{figure}
where we do not denote the external vertices. Similarly, $\unit$ denotes the following graph
\begin{figure}[H]
    \centering
    \includegraphics[width=0.15\linewidth]{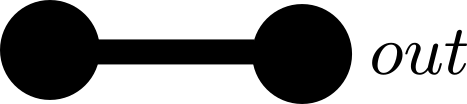}
\end{figure}
without the external vertex. The graphs $\comulti$ and $\counit$ are defined as their respective duals.\\
For a functor $F$ out of $\GrCob$ (or $\fatGrCob$ or $\pGrCob$), we define the following notation:
\begin{align*}
    A&:=F(*),\\
    \mu&:=F(\multi), & \eta&:=F(\unit),\\
    \nu&:=F(\comulti), & \varepsilon&:=F(\counit),\\
    p&:= \varepsilon\circ \mu, & c&:= \nu \circ \eta.
\end{align*} 

%The reader may be rightfully confused why we did all that work only to end up with a third characterisation of Frobenius algebras. A description that seems more complicated than the previous two. The reason why we went through this ordeal is that in the graded context the characterisations of Frobenius algebras we have developed so far fail to describe certain interesting examples. 
One advantage of the graph models of cobordism categories $\GrCob$, $\fatGrCob$ and $\pGrCob$ compared to the geometric cobordism categories $\mathrm{Cob}_2^{\mathrm{closed}}$, $\mathrm{Cob}_2^{\mathrm{open}}$ and $\Cob_2^{\mathrm{planar}}$ is that it becomes convenient to state that all commutative Frobenius algebras are symmetric and all symmetric Frobenius algebras are Frobenius algebras.

\begin{Prop}\label{prop:comm to symm}
    There is a canonical symmetric monoidal functor $\fatGrCob\to \GrCob$ given by forgetting the fat graph structure. Precomposition with this functor encodes the forgetful functor $\mathrm{ComFrob}(\mathcal{C})\to \mathrm{SymFrob}(\mathcal{C})$ for all symmetric monoidal categories $\mathcal{C}$.\\
    Moreover, there is a canonical monoidal functor $\pGrCob\to \fatGrCob$ given by forgetting the linear ordering on objects. Precomposition with this functor encodes the forgetful functor $\mathrm{SymFrob}(\mathcal{C})\to \mathrm{ComFrob}(\mathcal{C})$ for all monoidal categories $\mathcal{C}$.
\end{Prop}

\subsection{Graded Context}\label{subsec:graded context}

A first attempt to define a Frobenius algebra in a graded context may be to give a functor out of the cobordism categories we described until now. This fails to describe certain interesting examples in the graded context. We first specify what we mean by the graded context.\\

Let $R$ be a commutative ring and let $A=\{A_i\}_{i\in \Z}$, $B=\{B_i\}_{i\in \Z}$ be two $\Z$-graded $R$-modules. We  define the $\Z$-graded $R$-module
\begin{align*}
    \Rgrmod_\bullet(A,B)&=\{\Rgrmod_{d}(A,B)\}_{d\in \Z},\\
    \Rgrmod_d(A,B)&=\{\phi=(\phi_i)_{i\in \Z}\mid \phi_i\colon A_i\to B_{i+d} \ R\text{-linear}\}.
\end{align*}
We consider two notions of categories of graded $R$-module. Firstly, $\Rgrmod_0$ the category of $\Z$-graded $R$-modules with morphisms $\Rgrmod_0(A,B)$. Secondly, $\Rgrmod_\bullet$ the category with morphisms $\Rgrmod_\bullet(A,B)$.\\
We note that both categories $\Rgrmod_\bullet$ and $\Rgrmod_0$ can be enriched over $\Rgrmod_0$ (or $\Rgrmod_\bullet$). Furthermore, if $R$ is an $S$-algebra, both are enriched over $S$$\mathrm{Mod}_0$ (or $S\mathrm{Mod}_\bullet$).\\
Both are symmetric monoidal categories with unit
\begin{align*}
    \mathbbm{1}_d:=\begin{cases}
        \Z &\text{if } d=0,\\
        0&\text{else,}
    \end{cases}
\end{align*}
and the tensor product:
\begin{align*}
    (A\otimes B)_k=\bigoplus_{i+j=k} A_i\otimes_R B_j 
\end{align*}
with the twist induced by the Koszul sign
\begin{align*}
    \tau\colon A_i\otimes B_j &\to B_j\otimes A_i,\\
    a\otimes b &\mapsto (-1)^{ij} b\otimes a.
\end{align*}
Recall that we denote the shift $\Sigma A$ as the object with $(\Sigma A)_{i+1}=A_i$.\\

We fix a (symmetric) monoidal functor $F\colon \GrCob\to \Rgrmod_\bullet$ (or out of $\fatGrCob$ or $\pGrCob$) enriched over sets. In particular, we allow all degrees for $\mu=F(\multi)$ and $\nu=F(\comulti)$. We note that unitality implies that $|\mu|+|\eta|=|\id|=0$ and counitality implies $|\nu|+|\varepsilon|=|\id|=0$.
\begin{Prop}\label{prop:finite rank}
    Let $R$ be a principal ideal domain. Let $F\colon \pGrCob\to \Rgrmod_\bullet$ be a monoidal functor. Then $A$ is free of finite rank.
\end{Prop}
\begin{proof}
    A consequence of unitality, counitality and the Frobenius relation is that
    \begin{align}\label{eq:dual}
        (p\otimes \id)\circ (\id \otimes c)=\id=(\id \otimes p)\circ (c\otimes \id).
    \end{align}
    This means that $A$ is dualizable which is equivalent to $A$ being projective and finitely generated. As $R$ is a principal ideal domain, this is equivalent to $A$ being free of finite rank.
\end{proof}

\begin{Prop}\label{prop:trivial}
    Let $R$ be a principal ideal domain of characteristic not 2.\\
    Let $F\colon \GrCob\to \Rgrmod_\bullet$, $F\colon \fatGrCob\to \Rgrmod_\bullet$ or $F\colon \pGrCob\to \Rgrmod_\bullet$ be a monoidal functor. If $|\mu|-|\nu|$ is odd, then $A$ is the zero module.
\end{Prop}
\begin{proof}
    By Proposition \ref{prop:comm to symm}, it suffices to show the statement for $\pGrCob$ and, by Proposition \ref{prop:finite rank}, we can fix a basis $\alpha_1,\dots, \alpha_k$ of $A$.\\
    We note that $n:=|p|=|\mu|+|\varepsilon|=|\mu|-|\nu|$ and $|c|=|\nu|+|\eta|=|\nu|-|\mu|=-n$ are both of odd degree. We write 
    \begin{align*}
        c(1)=\sum_j\alpha_j\otimes \beta_j
    \end{align*}
    for some $\beta_j\in A$. By duality, $\beta_1,\dots,\beta_k$ also defines a basis. The fact that $|c|=-n$ shows that $|\beta_i|=-|\alpha_i|-n$ for all $i$.\\
    If $A\neq 0$, we may assume that $k>0$ and thus $\alpha_1\neq 0$. We want to exploit \eqref{eq:dual} for $\alpha_1$ and $\beta_1$. We first apply the right hand side to $\alpha_1$:
    \begin{align*}
        \alpha_1=&\ (\id \otimes p)\circ (c\otimes \id)(\alpha_1)=(\id \otimes p)\left(\sum_{j}\alpha_j\otimes \beta_j\otimes \alpha_1\right)=\sum_j (-1)^{n|\alpha_j|}\alpha_j p(\beta_j\otimes \alpha_1).
    \end{align*}
    As $\alpha_1,\dots,\alpha_k$ is a basis, we must have
    \begin{align}\label{eq:triv1}
         p(\beta_1\otimes \alpha_1)=(-1)^{n|\alpha_1|}.
    \end{align}   
    On the other hand, when we apply the left hand side of \eqref{eq:dual} to $\beta_1$, we find
    \begin{align*}
        \beta_1= (p\otimes \id)\circ (\id \otimes c)(\beta_1)=(p\otimes \id)\left((-1)^{n|\beta_1|}\beta_1\otimes \sum_{j} \alpha_j\otimes \beta_j\right)=(-1)^{n|\beta_1|}\sum_{j}p(\beta_1\otimes \alpha_j)\beta_j.
    \end{align*}
    As $\beta_1,\dots,\beta_k$ is a basis, we must have
    \begin{align}\label{eq:triv2}
        p(\beta_1\otimes \alpha_1)=(-1)^{n(-n-|\alpha_1|)}=(-1)^{n+n|\alpha_1|}.
    \end{align}
    As $n$ is odd and $(-1)^n=-1\neq 1$ in $R$, we have the desired contradiction between \eqref{eq:triv1} and \eqref{eq:triv2}. Therefore there is no element $\alpha_1\neq 0$ and thus $A=0$.
\end{proof}

In section \ref{sec:exmp}, we give examples of structures that look like Frobenius algebras such that $|\mu|-|\nu|$ is odd and are non-zero. In other words, this characterisation of a Frobenius algebra misses key examples.

\section{Graded Frobenius Algebras: PROPs}\label{sec:gfa PROPs}

The goal of this section is to define graph cobordism categories which already encode the signs and degrees. Concretely, we define categories $\GrCob_{c,d}$, $\fatGrCob_{c,d}$ and $\pGrCob$ enriched over graded abelian groups where the product has degree $c$ and the coproduct degree $d$. This gives us a new definition of graded Frobenius algebras as functors out of $\GrCob_{c,d}$, $\fatGrCob_{c,d}$ or $\pGrCob_{c,d}$ to a category enriched over graded abelian groups $\Zgrmod_0$.\\
In Subsection \ref{subsec:prop}, we prove some properties immediately from the definition: suspending an algebra over $\GrCob_{c,d}$ gives an algebra in $\GrCob_{c-1,d+1}$. Tensoring algebras over $\GrCob_{c,d}$ and $\GrCob_{c',d'}$ gives an algebra over $\GrCob_{c+c',d+d'}$. Analogous statements for $\fatGrCob_{c,d}$ and $\pGrCob_{c,d}$ can be proved in the same way.\\
In Subsection \ref{subsec:dioperad}, we study the canonical subdioperad of $\GrCob_{c,d}$ which is generated by forests. If $c+d$ is odd, then this subdioperad has $2$-torsion cokernel in the full PROP. This dioperad naturally shows up to describe an infinite-dimensional structure carrying a Frobenius relation as in \ref{subsec:Rabinowitz}.

\subsection{Definition}
In the following computations and definitions, we heavily use the concept of a top exterior power of an abelian group. Unfortunately, we require multiple slightly different versions of a top exterior power. We hope to ameliorate this situation by giving a clear overview of the versions we use.
\begin{itemize}
    \item Let $M$ be a free abelian group of rank $n$, we denote
\begin{itemize}
    \item $\Or(M):=\Lambda^n M$ the top exterior power as a graded abelian group concentrated in degree 0;
    \item $\det(M):=\Sigma^{-n}\Or(M)$ the top exterior power as a graded abelian group concentrated in degree $-n$.
\end{itemize}
\item Let $A$ be a free graded abelian group of rank 1. It is invertible with respect to $\otimes$, the tensor product of graded abelian groups, with its inverse $A^{-1}:=\Zgrmod_\bullet(A,\mathbbm{1})$.
\item Let $A=\{A_i\}_{i\in \Z}$ be a free graded abelian group of finite rank, we denote 
\begin{align*}
    \bolddet(A):=\bigotimes\limits_{i\in \Z} \det(A_{i})^{(-1)^i}:=\bigotimes_{i=-n}^n \det(A_i)^{(-1)^{i}}
\end{align*}
where $n\geq 0$ is such that $A_i=0$ for $|i|>n$.
\item Let $G$ be a graph and $c,d\in \Z$, we denote
\begin{itemize}
    \item $\bolddet(G,\partial_{in}):=\bolddet(H_*(|G|,|L_{in}(G)|))$;
    \item $\bolddet(G,\partial_{out}):=\bolddet(H_*(|G|,|L_{out}(G)|))$;
    \item ${\bolddet}_{c,d}(G):=\bolddet(G,\partial_{in})^{\otimes c}\otimes \bolddet(G,\partial_{out})^{\otimes d}$.
\end{itemize}
\end{itemize}

\begin{Lem}\label{lem:det}
    Let $A,B,C\in \Zgrmod_\bullet$ free of finite total rank. Then the following holds:
    \begin{enumerate}[(a)]
        \item\label{det:1} the determinant $\bolddet(A)$ is of rank 1 and concentrated in degree $-\chi(A)$, where $\chi(A)$ is the Euler characteristic of $A$;
        \item\label{det:2} there exists a natural isomorphism 
        \begin{align*}
            \bolddet(A\oplus B)\cong \bolddet(A)\otimes \bolddet(B);
        \end{align*}
        \item\label{det:3} there exists a natural isomorphism
        \begin{align*}
            \bolddet(A\oplus \Sigma A)\cong \mathbbm{1}.
        \end{align*}
        This implies a natural isomorphism $\bolddet(A)^{-1}\cong \bolddet(\Sigma A)$;
        \item\label{det:4} if $0\to A\to B\to C\to 0$ is a short exact sequence, there exists a canonical isomorphism $\bolddet(A)\otimes \bolddet(C)\cong \bolddet(B)$; 
        \item\label{det:5} if $\partial$ a differential on $A$ such that $H_*(A):=H_*(A,\partial)$ and $A$ are finitely generated free $\Z$-modules, then there exists a natural isomorphism: 
        \begin{align*}
            \bolddet(A)\cong \bolddet(H_*(A)).
        \end{align*}
    \end{enumerate}
\end{Lem}
\begin{proof}
    Statements \eqref{det:1}, \eqref{det:2} and \eqref{det:3} follow directly from the definition and the following facts for finitely generated free $\Z$-modules $M_1,\dots,M_k$:
    \begin{align*}
        \det\left(\bigoplus_{i=1}^k M_i\right)\cong \bigotimes_{i=1}^k \det(M_i) \text{ and }  \bolddet(\Sigma M_i)\cong\det(M_i)^{-1}.
    \end{align*}
    For statement \eqref{det:4}, we note that we can choose a splitting $f\colon C\to B$ because $C$ is free. This gives an isomorphism 
    \begin{align*}
        \bolddet(A)\otimes\bolddet(C)\cong \bolddet(A\oplus C)\cong \bolddet(B).
    \end{align*}
    This is independent of the choice of splitting because two splittings of $B\to C$ differ by elements in $A$ and thus the difference between the two induced isomorphisms is zero after tensoring with $\bolddet(A)$.\\ 
    For statement \eqref{det:5}, the short exact sequence
    \begin{align*}
        0\to \im(\partial\colon A_{i+1}\to A_i)\to \ker(\partial\colon A_i\to A_{i-1})\to H_i(A)\to 0
    \end{align*}
    gives
    \begin{align*}
        \bolddet(\ker(\partial\colon A_i\to A_{i-1}))&\overset{\eqref{det:4}}\cong \bolddet(\im(\partial\colon  A_{i+1}\to A_i))\otimes \bolddet(H_i(A))\\
        \implies \bolddet(H_i(A))&\!\;\cong \bolddet(\ker(\partial\colon A_i\to A_{i-1}))\otimes \bolddet(\im(\partial\colon  A_{i+1}\to A_i))^{-1}\\
        &\overset{\eqref{det:3}}{\cong} \bolddet(\ker(\partial\colon A_i\to A_{i-1}))\otimes \bolddet(\Sigma\:\im(\partial\colon  A_{i+1}\to A_i))\\
        &\!\;\cong \bolddet(\ker(\partial\colon A_i\to A_{i-1}))\otimes \bolddet(\im(\partial\colon  A_{i}\to A_{i-1}))
    \end{align*}
    for all $i\in \Z$.
    On the other hand, the short exact sequence
    \begin{align*}
        0\to \ker(\partial\colon A_i\to A_{i-1})\to A_i\overset{\partial}\to  \im(\partial\colon  A_{i}\to A_{i-1})\to 0
    \end{align*}
    gives 
    \begin{align*}
        \bolddet(A)\overset{\eqref{det:3}}{\cong} \bolddet(\ker(\partial\colon A_i\to A_{i-1}))\otimes \bolddet(\im(\partial\colon  A_{i}\to A_{i-1}))
    \end{align*}
    which concludes the proof.
\end{proof}

We note that a morphism of graph cobordisms $G\to G'$ specifies homotopy equivalences
\begin{align*}
    (|G|,|L_{in}(G)|) &\to (|G'|,|L_{in}(G')|),\\
    (|G|,|L_{out}(G)|) &\to (|G'|,|L_{out}(G')|)
\end{align*}
and thus isomorphisms 
\begin{align*}
    \bolddet(G,\partial_{in})&\to \bolddet(G',\partial_{in}),\\
    \bolddet(G,\partial_{out})&\to \bolddet(G',\partial_{out}),\\
    {\bolddet}_{c,d}(G)&\to {\bolddet}_{c,d}(G').
\end{align*}
Therefore $\det_{c,d}\colon \underline{\GrCob}(X,Y)\to\Zgrmod_0$ defines a functor for all $c,d\in \Z$ and finite sets $X,Y$.\\
Moreover, we can glue orientations together:
\begin{Lem}\label{lem:def of compos}
    Let $G,G'$ be composable graphs in $\underline{\GrCob}$, there is a natural isomorphism
    \begin{align*}
        \mathrm{comp}\colon {\bolddet}_{c,d}(G')\otimes {\bolddet}_{c,d}(G)\cong {\bolddet}_{c,d}(G'\circ G).
    \end{align*}
\end{Lem}
\begin{proof}
    We recall:
    \begin{align*}
        {\bolddet}_{c,d}(G') \otimes {\bolddet}_{c,d}(G)&\cong\bolddet(G',\partial_{in})^{\otimes c}\otimes \bolddet(G',\partial_{out})^{\otimes d}\otimes \bolddet(G,\partial_{in})^{\otimes c}\otimes \bolddet(G,\partial_{out})^{\otimes d}\\
        &\cong \bolddet(G',\partial_{in})^{\otimes c}\otimes \bolddet(G,\partial_{in})^{\otimes c}\otimes \bolddet(G',\partial_{out})^{\otimes d}\otimes \bolddet(G,\partial_{out})^{\otimes d}
    \end{align*}
    and
    \begin{align*}
        {\bolddet}_{c,d}(G'\circ G)\cong \bolddet(G'\circ G,\partial_{in})^{\otimes c}\otimes \bolddet(G'\circ G,\partial_{out})^{\otimes d}.
    \end{align*}
    It thus suffices to give isomorphisms 
    \begin{align*}
        \bolddet(G',\partial_{in})\otimes \bolddet(G,\partial_{in})&\cong \bolddet(G'\circ G,\partial_{in})\\
        \bolddet(G',\partial_{out})\otimes \bolddet(G,\partial_{out})&\cong \bolddet(G'\circ G,\partial_{out})
    \end{align*}
    We give the first isomorphism, the second one is constructed analogously.\\
    We can canonically embed $|G|\subseteq |G'\circ G|$. Let $C_*(|G|,|\partial_{in}|)$, $C_*(|G'|,|\partial_{in}|)$, $C_*(|G'\circ G|,|\partial_{in}|)$ and $C_*(|G'\circ G|,|G|)$ be the cellular chain complexes of the pairs. We note that the inclusion
    \begin{align*}
         C_*(|G'\circ G|,|G|)\to C_*(|G'|,|\partial_{in}|)
    \end{align*}
    is an isomorphism.\\
    The triple of spaces $(|G'\circ G|,|G|,|\partial_{in}|)$ gives a short exact sequence of chain complexes:
    \begin{center}
        \begin{tikzcd}
            0 \ar[r] & C_*(|G|,|\partial_{in}|) \ar[r] & C_*(|G'\circ G|,|\partial_{in}|)  \ar[r] & \overset{ C_*(|G'|,|\partial_{in}|)  }{\overbrace{C_*(|G'\circ G|, |G|)}} \ar[r] & 0.
        \end{tikzcd}
    \end{center}
    Lemma \ref{lem:det} \eqref{det:4} and \eqref{det:5} now give the desired isomorphism.
\end{proof}

\begin{Def}\label{def:grcob}
    Let $c,d\in \Z$.\\
    The category $\GrCob_{c,d}$ is defined with objects finite sets and morphism spaces
    \begin{align*}
        \GrCob_{c,d}(X,Y)=\colim_{\underline{\GrCob}(X, Y)} {\bolddet}_{c,d}(G). %This is a 1-colimit!
    \end{align*}
    The category $\fatGrCob_{c,d}$ is defined with objects finite sets and morphism spaces
    \begin{align*}
        \fatGrCob_{c,d}(X,Y)=\colim_{\underline{\fatGrCob}(X, Y)} {\bolddet}_{c,d}(G). %This is a 1-colimit!
    \end{align*}
    The category $\fatGrCob_{c,d}$ is defined with objects finite sets and morphism spaces
    \begin{align*}
        \pGrCob_{c,d}(X,Y)=\colim_{\underline{\pGrCob}(X, Y)} {\bolddet}_{c,d}(G). %This is a 1-colimit!
    \end{align*}
    Lemma \ref{lem:def of compos} lets us define the composition in these categories.
\end{Def}

We can now define a graded Frobenius algebra:

\begin{Def}\label{def:graded Frob alg}
    Let $\mathcal{C}$ be a monoidal category enriched over $\Zgrmod_0$ and let $c,d\in \Z$.\\
    A \textit{$(c,d)$-graded Frobenius algebra} is a monoidal functor $F\colon\GrCob_{c,d}\to \mathcal{C}$ enriched over $\Zgrmod_0$, i.e.\ degree-preserving on morphism spaces. The category is denoted by 
    \begin{align*}
        \mathrm{grFrob}_{c,d}(\mathcal{C})=\Fun_{\Zgrmod_0}^{\mathrm{mon}}(\pGrCob_{c,d},\mathcal{C}).
    \end{align*}
    If $\mathcal{C}$ is also symmetric, a \textit{$(c,d)$-graded commutative Frobenius algebra} is a symmetric monoidal functor $F\colon\GrCob_{c,d}\to \mathcal{C}$ and the category is denoted by 
    \begin{align*}
        \mathrm{grComFrob}_{c,d}(\mathcal{C})=\Fun_{\Zgrmod_0}^\otimes(\GrCob_{c,d},\mathcal{C}).
    \end{align*}
    and a \textit{$(c,d)$-graded symmetric Frobenius algebra} is defined as such a functor $F\colon \GrCob_{c,d}\to \mathcal{C}$ and the category is denoted by
    \begin{align*}
        \mathrm{grSymFrob}_{c,d}(\mathcal{C})=\Fun_{\Zgrmod_0}^\otimes(\fatGrCob_{c,d},\mathcal{C}).
    \end{align*}
\end{Def}

\subsection{Properties}\label{subsec:prop}

In this Subsection, we prove some first properties of these category. The definitions of the categories $\GrCob_{c,d}$, $\fatGrCob_{c,d}$ and $\pGrCob_{c,d}$ we gave are admittedly hard to work with. However, the morphism spaces can be described in a more convenient way.

\begin{Prop}\label{prop:hom sets direct sum}
    There exist isomorphisms
    \begin{align*}
        {\GrCob_{c,d}}(X,Y)\cong\bigoplus_{[G]\in \GrCob(X, Y)} {\bolddet}_{c,d}(G)_{\pi_1(|\underline{\GrCob}(X,Y)|,G)},
    \end{align*}
    \begin{align*}
        {\fatGrCob_{c,d}}(X,Y)\cong \bigoplus_{[G]\in\fatGrCob(X,Y)} {\bolddet}_{c,d}(G)_{\pi_1(|\underline{\fatGrCob}(X,Y)|,G)}.
    \end{align*}
    and 
    \begin{align*}
        {\pGrCob_{c,d}}(X,Y)\cong \bigoplus_{[G]\in\pGrCob(X,Y)} {\bolddet}_{c,d}(G)_{\pi_1(|\underline{\pGrCob}(X,Y)|,G)}.
    \end{align*}
\end{Prop}
\begin{proof}
    We note that the functor $\bolddet_{c,d}$ sends all maps in $\underline{\GrCob}(X,Y)$ to isomorphisms. Therefore, $\bolddet_{c,d}$ factors through the localisations $\underline{\GrCob}(X,Y)\to |\underline{\GrCob}(X, Y)|$ and $\underline{\fatGrCob}(X,Y)\to |\underline{\fatGrCob}(X, Y)|$. Here $|\mathcal{C}|$ is interpreted as the $\infty$-groupoid given by the category $\mathcal{C}$ with all morphisms formally inverted.\\
    We can thus compute the colimit giving the morphism space between finite sets $X$ and $Y$ as follows:
    \begin{align*}
        \colim_{G\in \underline{\GrCob}(X, Y)} {\bolddet}_{c,d}(G)=\colim_{G\in |\underline{\GrCob}(X,Y)|} {\bolddet}_{c,d}(G)
    \end{align*}
    We denote by $\pi_1(|\underline{\GrCob}(X,Y)|,G)$ the fundamental group of $|\underline{\GrCob}(X,Y)|$ based at the point corresponding to $G$. An element in $\pi_1(|\underline{\GrCob}(X,Y)|,G)$ is represented by a zig-zag of graph morphisms starting and ending at $G$:
    \begin{align*}
        G=G_0\to G_1\leftarrow G_2\to \dots \leftarrow G_k=G.
    \end{align*}
    Such a zig-zag gives a zig-zag of maps
    \begin{align*}
        \bolddet_{c,d}(G_0)\to \bolddet_{c,d}(G_1)\leftarrow \bolddet_{c,d}(G_2)\to \dots \leftarrow \bolddet_{c,d}(G_k).
    \end{align*}
    But all these maps are isomorphism and we thus obtain an action of $\pi_1(|\underline{\GrCob}(X,Y)|,G)$ on $\bolddet_{c,d}(G)$.\\
    The functor $\bolddet_{c,d}$ lands in a $1$-category. Therefore, we can compute the colimit on $\tau_{\leq1}\underline{\GrCob}(X,Y)$, the $1$-truncation of the $\infty$-groupoid $|\underline{\GrCob}(X,Y)|$.\\
    Moreover, two graphs $G,G'\in \underline{\GrCob}(X, Y)$ are equivalent in $|\underline{\GrCob}(X,Y)|$ if and only if they represent the same element in ${\GrCob}(X,Y)=\pi_0(\underline{\GrCob}(X,Y))$. We thus can compute the colimit defining ${\GrCob_{c,d}}(X,Y)$ by:
    \begin{align*}
        {\GrCob_{c,d}}(X,Y)\cong\colim_{G\in \tau_{\leq1}\underline{\GrCob}(X,Y) }\bolddet_{c,d}(G)\cong\bigoplus_{[G]\in \GrCob(X, Y)} {\bolddet}_{c,d}(G)_{\pi_1(|\underline{\GrCob}(X,Y)|,G)}
    \end{align*}
    The proof for $\fatGrCob$ and $\pGrCob$ is analogous.
\end{proof}

\begin{Rem}
    The group $\pi_1(|\underline{\GrCob}(X,Y)|,G)$ is isomorphic to $\pi_0(\mathrm{hAut}(|G|))$: a zig-zag of graphs morphisms is the same as a simple homotopy equivalence from the CW-complex $|G|$ to itself. For the $1$-dimensional CW complex $|G|$, any homotopy equivalence is homotopic to a simple homotopy equivalence.\\
    On the other hand, the group $\pi_1(|\underline{\fatGrCob}(X,Y)|,G)$ is isomorphic to the mapping class group $\pi_0(\mathrm{Homeo}^+(\fatCob(G),\partial_{in}( \fatCob(G))\cup \partial_{out}(\fatCob(G)))$ that fixes the boundary point-wise. This can be seen for example by the open part of Theorem A in \cite{santander2015comparing}.\\
    We note that $\underline{\pGrCob}(X,Y)$ is a full subcategory of $\underline{\fatGrCob}(X,Y)$ where if $G\to G'\in \underline{\fatGrCob}(X,Y)$ and either $G$ or $G'$ is a planar graph, then so is the other. Therefore, we find also
    \begin{align*}
        \pi_1(|\underline{\pGrCob}(X,Y)|,G)\cong  \pi_0(\mathrm{Homeo}^+(\fatCob(G),\partial_{in}( \fatCob(G))\cup \partial_{out}(\fatCob(G))).
    \end{align*}
\end{Rem}

% The composition maps in $\GrCob_{c,d}$ and $\fatGrCob_{c,d}$ are given as follows: let $X,Y,Z$ be finite sets. The maps
% \begin{align*}
%     {\bolddet}_{c,d}(G)\otimes {\bolddet}_{c,d}(G')\to \colim_{(G,G')\in \Gr_{X\to Y}\times \Gr_{Y\to Z}} {\bolddet}_{c,d}(G)\otimes {\bolddet}_{c,d}(G')
% \end{align*}
% assemble to a map
% \begin{align*}
%     \colim_{G\in \Gr_{X\to Y}} {\bolddet}_{c,d}(G)\otimes \colim_{G'\Gr_{Y\to Z}} {\bolddet}_{c,d}(G')&\to \colim_{(G,G')\in \Gr_{X\to Y}\times \Gr_{Y\to Z}} {\bolddet}_{c,d}(G)\otimes {\bolddet}_{c,d}(G').
% \end{align*}
% Lemma \ref{lem:def of compos} gives us a map
% \begin{align*}
%     \colim_{(G,G')\in \Gr_{X\to Y}\times \Gr_{Y\to Z}} {\bolddet}_{c,d}(G')\otimes {\bolddet}_{c,d} (G)\cong& \colim_{(G,G')\in \Gr_{X\to Y}\times \Gr_{Y\to Z}} {\bolddet}_{c,d}(G'\circ G).
% \end{align*}
% Finally by the composition functor
% \begin{align*}
%     \Gr_{X\to Y}\times \Gr_{Y\to Z}&\to \Gr_{X\to Z},\\
%     (G,G')&\mapsto G'\circ G,
% \end{align*}
% we obtain a map
% \begin{align*}
%     \colim_{(G,G')\in \Gr_{X\to Y}\times \Gr_{Y\to Z}} {\bolddet}_{c,d}(G'\circ G)\to \colim_{G''\in \Gr_{X\to Z}} {\bolddet}_{c,d}(G'').
% \end{align*}
% This combines to the desired composition map
% \begin{align}\label{eq:def compo}
%     \mathrm{comp}\colon\colim_{G\in \Gr_{X\to Y}} {\bolddet}_{c,d}(G)\otimes \colim_{G'\Gr_{Y\to Z}} {\bolddet}_{c,d}(G')\to \colim_{G''\in \Gr_{X\to Z}} {\bolddet}_{c,d}(G'').
% \end{align}

\begin{Exmp}\label{exmp:GrCob 00}
    In the case $c=d=0$, we have ${\bolddet}_{0,0}(G)=\mathbbm{1}$ and the automorphisms of $G$ act trivially. We therefore compute:
    \begin{align*}
        {\GrCob_{0,0}}(X,Y)&\cong\bigoplus_{[G]\in \GrCob(X, Y)} \mathbbm{1},\\
        {\fatGrCob_{0,0}}(X,Y)&\cong\bigoplus_{[G]\in \fatGrCob(X, Y|} \mathbbm{1}
    \end{align*}
    and 
    \begin{align*}
        {\pGrCob_{0,0}}(X,Y)&\cong\bigoplus_{[G]\in \pGrCob(X, Y|} \mathbbm{1}
    \end{align*}
    We find
    \begin{align*}
        \mathrm{grComFrob}_{0,0}(\Rgrmod_\bullet)&= \Fun_{\Zgrmod_0}^\otimes(\GrCob_{0,0},\Rgrmod_\bullet)\\
        &\cong\Fun^\otimes(\GrCob,\Rgrmod_0)\\
        &\cong\mathrm{ComFrob}(\Rgrmod_0).
    \end{align*}
    Therefore a $(0,0)$-graded commutative Frobenius algebra coincides with the classical definition of a commutative Frobenius algebras with maps in degree $0$. An analogous statement holds for (symmetric) Frobenius algebras.
\end{Exmp}

In the following, we study the relation between $\GrCob_{c,d}$ for different choices of $c,d\in \Z$. For this, we use the symmetric monoidal structure on PROPs. 
\begin{Rem}
    Analogous statements hold for $\fatGrCob$ and $\pGrCob$ with the caveat that $\pGrCob$ is not symmetric monoidal and thus not a PROP but a PRO. However, the construction for tensor product of PROs and the suspension PRO is analogous to the PROP case.
\end{Rem}

\begin{Not}
    We denote the \textit{tensor product} of two PROPs in $\Zgrmod_0$ for the coordinate-wise tensor product: for PROPs $P,Q$
    \begin{align*}
        (P\otimes Q)(X,Y):=P(X,Y)\otimes Q(X,Y)
    \end{align*}
    for finite sets $X,Y$ with the component-wise composition maps.
\end{Not}

\begin{Exmp}\label{exmp:End AotimesB}
    Let $(\mathcal{C},\otimes, \mathbbm{1},\tau)$ be a symmetric monoidal category enriched over $\Zgrmod_0$ and let $A,B\in \mathcal{C}$ be objects. There is an isomorphism of PROPs 
    \begin{align*}
        \End_A\otimes \End_B\cong \End_{A\otimes B}
    \end{align*}
    given on maps $f\colon A^{\otimes k}\to A^{\otimes l}$ and $g\colon B^{\otimes k}\to B^{\otimes l}$ by
    \begin{align*}
        (A\otimes B)^{\otimes k}\to A^{\otimes k}\otimes B^{\otimes k}\overset{f\otimes g}{\to}A^{\otimes l}\otimes B^{\otimes l}\to (A\otimes B)^{\otimes l}.
    \end{align*}
    Indeed, this map is equivariant and commutes with composition and the product $\otimes$.
\end{Exmp}

\begin{Prop}\label{prop:c+c' d+d'}
    There exist canonical maps of PROPs or PROs
    \begin{align*}
        \GrCob_{c+c',d+d'}&\to \GrCob_{c,d}\otimes \GrCob_{c',d'},\\
        \fatGrCob_{c+c',d+d'}&\to \fatGrCob_{c,d}\otimes \fatGrCob_{c',d'}
    \end{align*}
    and 
    \begin{align*}
        \pGrCob_{c+c',d+d'}\to \pGrCob_{c,d}\otimes \pGrCob_{c',d'}
    \end{align*}
\end{Prop}
\begin{proof}
    There are canonical maps
    \begin{align}
    \begin{split}\label{eq:c+c' d+d'}
        \bolddet(G,\partial_{in})^{\otimes c+c'}&\otimes \bolddet(G,\partial_{out})^{\otimes d+d'}\\
        \cong&\bolddet(G,\partial_{in})^{\otimes c}\otimes \bolddet(G,\partial_{out})^{\otimes d}\otimes \bolddet(G,\partial_{in})^{\otimes c'}\otimes \bolddet(G,\partial_{out})^{\otimes d'}.
    \end{split}
    \end{align}
    This then gives maps 
    \begin{align*}
        \bolddet_{c+c',d+d'}(G)\to \colim_{G'\in \underline{\GrCob}(X,Y)}\bolddet_{c,d}(G')\otimes \colim_{G'\in \underline{\GrCob}(X,Y)}\bolddet_{c',d'}(G')
    \end{align*}
    for all $G\in \underline{\GrCob}(X,Y)$. Taking the colimit over $\underline{\GrCob}(X,Y)$ thus gives the desired map for the $\GrCob$. The proof for $\fatGrCob$ and $\pGrCob$ is analogous.
\end{proof}

\begin{Rem}
    Even though the map $\eqref{eq:c+c' d+d'}$ used in the construction is an isomorphism, the resulting map of PROPs is not a isomorphism in general. For example if $c+d$ and $c'+d'$ are odd, Proposition \ref{prop:subdioperad} shows that $\GrCob_{c,d}$ and $\GrCob_{c',d'}$ is $2$-torsion on summands coming from graphs which are not forests. This is not true for $\GrCob_{c+c',d+d'}$.
\end{Rem}

\begin{Cor}\label{cor:c+c' d+d'}
    Let $A$ be an algebra over $\GrCob_{c,d}$ and $B$ be an algebra over $\GrCob_{c',d'}$. Then $A\otimes B$ inherits the structure of an algebra over $\GrCob_{c+c',d+d'}$. Analogous statements hold for $\fatGrCob$ and $\pGrCob$.
\end{Cor}
\begin{proof}
    We show the statement for $\GrCob_{c,d}$. The proof for $\fatGrCob$ and $\pGrCob$ is analogous.\\
    Giving a $\GrCob_{c,d}$-algebra structure on $A$ is the same thing as giving a map of PROPs $\GrCob_{c,d}\to \mathrm{End}_A$ (see \cite{markl2008operads}). Similarly, we get a map of PROPs $\GrCob_{c',d'}\to \End_B$. Precomposing with the map from Proposition \ref{prop:c+c' d+d'} and postcomposing with the map from Example \ref{exmp:End AotimesB}, we thus get a map of PROPs
    \begin{align*}
        \GrCob_{c+c',d+d'}&\to\GrCob_{c,d}\otimes \GrCob_{c',d'}\to \End_A\otimes \End_{B}\to \End_{A\otimes B} 
    \end{align*}
    and thus a $\GrCob_{c+c',d+d'}$-algebra structure on $A\otimes B$.
\end{proof}

Another relation between $\GrCob_{c,d}$-algebras for different $c,d$ comes from suspending algebras. This operation has a corresponding operation on the level of PROPs:

\begin{Def}\label{def:susp prop}
    The \textit{suspension PROP} in $\Zgrmod_0$ is defined as the endomorphism PROP $\End_{\Sigma \mathbbm{1}}$ of $\Sigma \mathbbm{1}$ in $\Zgrmod_\bullet$. Explicitly, it is given by
    \begin{align*}
        \Sigma(X,Y):={\Zgrmod_\bullet}((\Sigma \mathbbm{1})^{\otimes X},(\Sigma \mathbbm{1})^{\otimes Y})
    \end{align*}
    for finite sets $X,Y$.\\
    For a PROP $P$ in $\Zgrmod_0$, the \textit{suspended PROP} is defined as $\Sigma P:=\Sigma  \otimes P$.
\end{Def}

\begin{Prop}\label{prop:univ susp prop}
    Let $P$ be a PROP in $\Zgrmod_0$ and $A\in \Zgrmod_\bullet$. Then giving a $P$-algebra structure on $A$ is equivalent to giving a $\Sigma P$-algebra structure on $\Sigma A$. 
\end{Prop}
\begin{proof}
    We again use the characterisation of $A$ being a $P$-algebra as a PROP maps into $\End_A$. By definition, this is equivalent to giving for all finite sets $X,Y$ equivariant maps 
    \begin{align*}
        P(X,Y)\to \mathcal{C}(A^{\otimes X},A^{\otimes Y})
    \end{align*} 
    which commute with composition $\circ$ and the symmetric monoidal product $\otimes$. This is equivalent to giving equivariant maps
    \begin{align*}
        \Zgrmod_\bullet((\Sigma \mathbbm{1})^{\otimes X},(\Sigma \mathbbm{1})^{\otimes Y})\otimes P(X,Y)\to \Zgrmod_\bullet((\Sigma \mathbbm{1})^{\otimes X},(\Sigma \mathbbm{1})^{\otimes Y}) \otimes \Zgrmod_\bullet(A^{\otimes X},A^{\otimes Y})
    \end{align*}
    which commute with $\circ$ and $\otimes$. After postcomposing with the isomorphism from Example \ref{exmp:End AotimesB}, we get equivariant maps
    \begin{align*}
        \Sigma P(X,Y)\to \Zgrmod_\bullet((\Sigma A)^{\otimes X},(\Sigma A)^{\otimes Y})
    \end{align*}
    which commute with $\circ$ and $\otimes$. This is by definition equivalent to giving a PROP morphism $\Sigma P\to \mathrm{End}_{\Sigma A}$, which is equivalent to giving a $\Sigma P$ structure on $\Sigma A$.
\end{proof}

\begin{Rem}
    This construction is completely analogous as for PROs, operads, cyclic operads, dioperads and properads. However in the category of modular operads, there exists no suspension modular operad in $\Zgrmod_0$. Indeed there is a canonical pairing $\Sigma\mathbbm{1}\otimes \Sigma\mathbbm{1}\to \Sigma \mathbbm{1}$ in degree $-1$. This lets us identify 
    \begin{align*}
        \End_{\Sigma\mathbbm{1}} (X)=\Zgrmod_\bullet((\Sigma \mathbbm{1})^{\otimes X},\Sigma \mathbbm{1})\cong \Sigma^{1-|X|}\mathbbm{1}.
    \end{align*}
    For any elements $x\in X$ and $y\in Y$, this gives a map
    \begin{align*}
        \End_{\Sigma \mathbbm{1}}(X)\otimes\End_{\Sigma\mathbbm{1}}(Y)\to \End_{\Sigma \mathbbm{1}}((X\setminus x)\sqcup (Y\setminus y))
    \end{align*}
    in $\Zgrmod_0$. However, the map for elements $x\neq x'\in X$
    \begin{align*}
        \End_{\Sigma \mathbbm{1}}(X)\to \End_{\Sigma\mathbbm{1}}(X\setminus \{x,x'\})
    \end{align*}
    is of degree $-1$ and thus not in $\Zgrmod_0$. One way to amend this is to work with the language of twisted modular operads as introduced in \cite{Getzler}.
\end{Rem}

\begin{Prop}\label{prop:suspending c,d}
    For integers $c,d\in \Z$, there exist canonical isomorphisms
    \begin{align*}
        \Sigma \GrCob_{c,d}&\cong \GrCob_{c-1,d+1},\\
        \Sigma \fatGrCob_{c,d}&\cong \fatGrCob_{c-1,d+1}
    \end{align*}
    and
    \begin{align*}
        \Sigma \pGrCob_{c,d}&\cong \pGrCob_{c-1,d+1}
    \end{align*}
\end{Prop}
\begin{proof}
    We fix integers $c,d\in \Z$ and finite sets $X,Y$. The key observation is that for every graph cobordism $G$ between $X$ and $Y$ there is an isomorphism
    \begin{align}\label{eq:suspending c,d}
        {\bolddet}_{c-1,d+1}(G)\cong \Zgrmod_\bullet((\Sigma \mathbbm{1})^{\otimes X},(\Sigma \mathbbm{1})^{\otimes Y})\otimes {\bolddet}_{c,d}(G) 
    \end{align}
    which is natural in $X$ and $Y$.\\
    To see the isomorphism, we observe
    \begin{align*}
        {\bolddet}_{c-1,d+1}(G)=&\  \bolddet(G,\partial_{in})^{\otimes c-1}\otimes\bolddet(G,\partial_{out})^{\otimes d+1}\\
        \cong &\   \bolddet(G,\partial_{in})^{-1}\otimes \bolddet(G,\partial_{out})\otimes {\bolddet}_{c,d}(G).
    \end{align*}
    We consider the short exact sequences 
    \begin{center}
        \begin{tikzcd}[row sep=0]
            0 \ar[r ] & C_*(|\partial_{in}|)\ar[r] & C_*(|G|) \ar[r] & C_*(|G|,|\partial_{in})\ar[r ] &0,\\
            0 \ar[r ] & C_*(|\partial_{out}|)\ar[r] & C_*(|G|) \ar[r] & C_*(|G|,|\partial_{out})\ar[r ] &0.
        \end{tikzcd}
    \end{center}
    Then Lemma \ref{lem:det} \eqref{det:4} gives the natural isomorphisms
    \begin{align*}
        \bolddet(G,\partial_{in})^{-1}\otimes \bolddet(G,\partial_{out})\cong&\    \bolddet(|\partial_{in}|)\otimes\bolddet(|G|)^{-1}\otimes\bolddet(|G|)\otimes \bolddet(|\partial_{out}|)^{-1}\\
        \cong&\   \bolddet(|\partial_{in}|)\otimes \bolddet(|\partial_{out}|)^{-1} \\
        \cong & \det(\Z^{X})\otimes \det(\Z^{Y})^{-1} \\
        \cong &\   ((\Sigma \mathbbm{1})^{\otimes X})^{-1}\otimes (\Sigma \mathbbm{1})^{\otimes Y}\\
        \cong & \  \Zgrmod_\bullet((\Sigma \mathbbm{1})^{\otimes X},(\Sigma \mathbbm{1})^{\otimes Y})=\Sigma(X,Y).
    \end{align*}
    We recall that, in our convention, the determinant lives in the degree equal to the negative rank. Thus, we have for a finite set $X$ a natural isomorphism $\det(\Z^{X})\cong ((\Sigma \mathbbm{1})^{\otimes X})^{-1}$.\\
    This gives \eqref{eq:suspending c,d} and thus the proof.
\end{proof}
Combining Proposition \ref{prop:univ susp prop} and Proposition \ref{prop:suspending c,d}, we obtain the following:
\begin{Cor}\label{cor:suspending GCob alg}
    Giving a $(c,d)$-graded (commutative or symmetric) Frobenius algebra structure on $A:=F(*)\in \Zgrmod_\bullet$ is equivalent to giving a $(c-1,d+1)$-graded (commutative or symmetric) Frobenius algebra structure on $\Sigma A$.
\end{Cor}

\begin{Rem}\label{rem:one sided suspension}
    We note that $|\mu|+|\nu|=|\Sigma(\mu)|+|\Sigma(\nu)|$. Therefore this defines an invariant under suspension. In particular if $|\mu|+|\nu|\neq 0$, the algebra is not a suspension of an algebra with operations in degree 0. Adding a factor $\bolddet(G,\partial_{in})$ (or $\bolddet(G,\partial_{out})$) in the definition of the PROP $\GrCob_{c,d}$ helps us describe other algebras because this only changes the degree of the multiplication and not the comultiplication (or vice-versa). However this is an operation on the PROP and contrary to suspension does not have an obvious analogue on the level of algebras. See Remark \ref{rem:one sided susp on algebras}.
\end{Rem}

\subsection{Canonical Dioperad}\label{subsec:dioperad}

The PROP $\GrCob_{c,d}$ has a canonical subdioperad generated by forests. In Subsection \ref{subsec:Rabinowitz}, we study an example where this dioperad shows up naturally rather than the whole PROP. We study the structure of this canonical subdioperad.\\

For this, we recall the definition of a dioperad:
\begin{Def}
    A \textit{dioperad} consists of the following data a collection $\mathcal{D}(X,Y)$ for all finite sets $X,Y$ with
    \begin{itemize}
        \item a $(S_X,S_Y)$-bimodule structure on $\mathcal{D}(X,Y)$ where $S_X$ and $S_Y$ denote the symmetric group of $X$ and $Y$, respectively;
        \item a unit
        \begin{align*}
            1\colon \mathbbm{1}\to \mathcal{D}(*,*);
        \end{align*}
        \item composition maps 
        \begin{align*}
            {}_y\circ_x\colon \mathcal{D}(Y,Z)\otimes \mathcal{D}(W,X)\to \mathcal{D}( (Y\setminus \{y\})\sqcup W,  Z\sqcup (X\setminus \{x\}))
        \end{align*}
        for all finite sets $W,X,Y,Z$ and elements $x\in X$, $y\in Y$.
    \end{itemize}
    This data has to satisfy some unitality, associativity and equivariance relations. See \cite{gan2002koszul} for the precise formulation of those relations.
\end{Def}

\begin{Def}\label{def:subdiop}
    The \textit{forest dioperad} $\GrCob_{c,d}^{\mathrm{forest}}$ of $\GrCob_{c,d}$ is given by the collections 
    \begin{align}\label{eq:closed dioperad in PROP}
        \GrCob_{c,d}^{\textit{forest}}(X,Y)\subseteq {\GrCob_{c,d}}(X,Y)
    \end{align}
    generated by the graph cobordisms between $X$ and $Y$ which are forests with the restricted composition.
    % The \textit{forest dioperad} $\fatGrCob_{c,d}^{forest}$ is given by the collections 
    % \begin{align*}
    %     \fatGrCob_{c,d}^{forest}(X,Y)\subseteq \Hom_{\fatGrCob_{c,d}}(X,Y)
    % \end{align*}
    % generated by the fat graphs cobordisms between $X$ and $Y$ which are forests with the restricted composition.
\end{Def}

In general this dioperad carries less information than the PROP. But in some cases, we have the following proposition:
\begin{Prop}\label{prop:subdioperad}
    If $c,d\in \Z$ are of different parity, then the PROP $\GrCob_{c,d}$ is an extension by $\Z/2$ factors of dioperad $\GrCob_{c,d}^{\mathrm{forest}}$. In explicit terms, the inclusion \eqref{eq:closed dioperad in PROP} finite sets $X,Y$
    \begin{align*}
        \GrCob_{c,d}^{{\mathrm{forest}}}(X,Y)\to {\GrCob_{c,d}}(X,Y)
    \end{align*}
    has a cokernel which is $2$-torsion, for all finite sets $X,Y,Z$.
\end{Prop}
\begin{proof}
    We fix $c,d\in \Z$ of different parity. We use the description
    \begin{align*}
        {\GrCob_{c,d}}(X,Y)\cong\bigoplus_{[G]\in \GrCob(X, Y)} {\bolddet}_{c,d}(G)_{\pi_1(|\underline{\GrCob}(X,Y)|,G)}
    \end{align*}
    as in Proposition \ref{prop:hom sets direct sum}. We have an analogous description:
    \begin{align*}
        \GrCob_{c,d}^{{\mathrm{forest}}}(X,Y)\cong \bigoplus_{\substack{[G]\in \GrCob(X,Y)\\ G \text{ a forest}}} {\bolddet}_{c,d}(G)_{\pi_1(|\underline{\GrCob}(X,Y)|,G)}.
    \end{align*}
    Let $[G]\in \GrCob(X,Y)$ be an equivalence class of graphs cobordisms between $X$ and $Y$ with $G$ not a forest. We want to show that the summand
    \begin{align*}
        {\bolddet}_{c,d}(G)_{\pi_1(|\underline{\GrCob}(X,Y)|,G)}
    \end{align*}
    is $\Z/2$.\\
    After collapsing some edges, we can assume that $G$ has a tadpole. Then there is a graph automorphism of $G$ flipping the direction of the tadpole. This represents an element in $\pi_1(|\underline{\GrCob}(X,Y)|,G)$. This automorphism acts by $(-1)^c$ on $\bolddet(G,\partial_{in})^{\otimes c}$ and by $(-1)^d$ on $\bolddet(G,\partial_{out})^{\otimes d}$. See Subsection \ref{subsec:computational tools} for a detailed explanation how to derive these signs.\\
    It thus acts by $(-1)^{c+d}=-1$ on ${\bolddet}_{c,d}(G)$ which is a free group of rank one. This shows ${\bolddet}_{c,d}(G)_{\pi_1(|\underline{\GrCob}(X,Y)|,G)}=\Z/2$ and thus concludes the proof.
\end{proof}

\section{Graded Frobenius Algebra: Maps and Relations}\label{sec:gfa maps and rel}

In this section, we state the main theorem which gives a description of graded Frobenius algebras in terms of maps and relations. We will use this theorem to study some of the simplest examples of graded Frobenius algebras. Moreover, we discuss the suspension of algebras from the perspective of maps and relations.

\begin{Thm}\label{thm:gen and rel}
    Let $c,d\in \Z$ and let $(\mathcal{C},\otimes, \mathbbm{1},\tau)$ be a (strict) monoidal category enriched over graded abelian groups.\\
    The data $(A,\mu,\eta,\nu,\varepsilon)$ where
    \begin{align*}
        A&\in \mathcal{C},\\
        \mu&\in \mathcal{C}_{c}(A\otimes A, A), & \eta&\in \mathcal{C}_{-c}(\mathbbm{1},A),\\
        \nu&\in \mathcal{C}_d(A,A\otimes A), & \varepsilon&\in \mathcal{C}_{-d}(A,\mathbbm{1})
    \end{align*}
    defines a monoidal functor out of $\pGrCob_{c,d}$ into $\mathcal{C}$ uniquely up to isomorphism if and only if the following relations are satisfied:
    \begin{enumerate}[(i)]
        \item\label{item:thm ass} graded associativity $\mu\circ (\mu \otimes \id)=(-1)^c\mu\circ (\id \otimes \mu)$;
        \item\label{item:thm unit} graded unitality $(-1)^c\mu\circ(\eta \otimes \id)=(-1)^{\frac{c(c-1)}{2}}\id=\mu\circ (\id \otimes \eta)$;
        \item\label{item:thm coass} graded coassociativity $(\id \otimes \nu)\circ \nu=(-1)^d(\nu \otimes \id)\circ \nu$;
        \item\label{item:thm counit} graded counitality $(-1)^d(\id \otimes \varepsilon)\circ \nu=(-1)^{\frac{d(d-1)}{2}}\id =(\varepsilon\otimes \id)\circ \nu$;
        \item graded Frobenius relation $(\mu \otimes \id)\circ (\id\otimes \nu)=(-1)^{cd}\nu\circ \mu=(\id \otimes \mu)\circ (\nu \otimes \id)$.
    \end{enumerate}
    If $\mathcal{C}$ is also symmetric, the data $(A,\mu,\eta,\nu,\varepsilon)$ defines a symmetric monoidal functor out of $\GrCob_{c,d}$ into $\mathcal{C}$ uniquely up to isomorphism if and only if the maps satisfy (i)-(v) and
    \begin{itemize}
        \item[\textit{(vi)}] graded commutativity $\mu \circ \tau=(-1)^c \mu$
    \end{itemize}
    and it defines a symmetric monoidal functor out of $\fatGrCob_{c,d}$ into $\mathcal{C}$ uniquely up to isomorphism if and only if the maps satisfy (i)-(v) and
    \begin{itemize}
        \item[\textit{(vi')}] graded symmetry $\varepsilon\circ \mu \circ \tau=(-1)^c\varepsilon\circ \mu$.
    \end{itemize}
\end{Thm}

We postpone the proof of this theorem to section \ref{sec:proof} as it involves a lot of technical details. Instead we focus on what are the choices involved in those signs and what are the consequences.\\
The strategy for the theorem is to choose generators of the categories $\GrCob_{c,d}$, $\fatGrCob_{c,d}$ and $\pGrCob_{c,d}$ and study the generating relations between them. For this, we choose generators in $\bolddet_{c,d}(G)$ for a list of graphs $G$ such that all other graphs can be decomposed into them.

\begin{Def}
    A generator $\omega(G)\in\bolddet_{c,d}(G)$ is called an \textit{orientation} of the graph $G$.
\end{Def}

\begin{Rem}[Choices of Orientations]\label{rem:choice of orien}
    The relations are \textit{not} canonical. Namely, they depend on a choice of orientations
    \begin{align*}
        \omega(\multi) \in {\bolddet}_{c,d}(\multi),\:  \omega(\unit) \in {\bolddet}_{c,d}(\unit),\:\omega(\comulti) \in {\bolddet}_{c,d}(\comulti),\:  \omega(\counit)\in {\bolddet}_{c,d}(\counit).
    \end{align*}
    By definition, this depends on generators of
    \begin{align*}
        &\bolddet(\multi,\partial_{in})^{\otimes c}\otimes \bolddet(\multi,\partial_{out})^{\otimes d},\: & &\bolddet(\unit,\partial_{in})^{\otimes c}\otimes \bolddet(\unit,\partial_{out})^{\otimes d},\\
        &\bolddet(\comulti,\partial_{in})^{\otimes c}\otimes \bolddet(\comulti,\partial_{out})^{\otimes d},\; & &\bolddet(\counit,\partial_{in})^{\otimes c}\otimes \bolddet(\counit,\partial_{out})^{\otimes d}.
    \end{align*}
    We thus choose generators $\omega_{in}(G)\in \bolddet(G,\partial_{in})$ and $\omega_{out}(G)\in\bolddet(G,\partial_{out})$ for those four graphs and then define $\omega(G):=\omega_{in}(G)^{\otimes c}\otimes \omega_{out}(G)^{\otimes d}$.\\
    We note that $\bolddet(\multi,\partial_{out})\cong \bolddet(\unit,\partial_{out})\cong \bolddet(\comulti,\partial_{in})\cong \bolddet(\counit,\partial_{in})\cong \mathbbm{1}$. We thus choose the canonical generator 1. Moreover for $\bolddet(\unit,\partial_{in})$ and $\bolddet(\counit,\partial_{out})$, we choose the canonical generator coming from the canonical generator of $H_0(\unit,\partial_{in})$ and $H_0(\counit,\partial_{out})$. However, there is no canonical choice for generators
    \begin{align*}
        \omega_{in}(\multi)&\in \bolddet(\multi,\partial_{in}),& \omega_{out}(\comulti)&\in \bolddet(\comulti,\partial_{out}).
    \end{align*}
    Therefore, there exists no canonical choice of orientation
    \begin{align*}
        \omega(\multi):=\omega_{in}(\multi)^{\otimes c}\otimes 1^{\otimes d}
    \end{align*}
    for $c$ odd and no canonical choice of orientation
    \begin{align*}
        \omega(\comulti):=1^{\otimes c}\otimes \omega_{out}(\comulti)^{\otimes d}
    \end{align*}
    for $d$ odd. We make a choice in the proof in \ref{subsec:proof}. A different choice would lead to different signs in the unitality and counitality relation.
\end{Rem}

\begin{Rem}\label{rem:alex kai}
    In \cite{cieliebak2022cofrobenius}, Cieliebak and Oancea introduce a type of algebra called biunital coFrobenius. In the finite dimensional case, this corresponds to our notion of graded Frobenius algebra (see \cite[Proposition 5.14]{cieliebak2022cofrobenius}). The infinite dimensional case is discussed in \ref{subsec:Rabinowitz}.\\
    They define the graded unitality relation
    \begin{align*}
        \mu\circ(\eta \otimes \id)=\id=(-1)^c\mu\circ (\id \otimes \eta)
    \end{align*}
    compared to our graded unitality relation
    \begin{align*}
        (-1)^c\mu\circ(\eta \otimes \id)=(-1)^{\frac{c(c-1)}{2}}\id=\mu\circ (\id \otimes \eta).
    \end{align*}
    The two relations agree for $c\equiv 3,4$ (mod 4) and are different for $c\equiv 1,2$ (mod 4).\\
    Similarly, the signs in their graded counitality relation agree with our signs for $d\equiv 3,4$ (mod 4) and are different for $d\equiv 1,2$ (mod 4).\\
    For $(A,\mu,\eta,\nu,\varepsilon)$ satisfying the signs in our relations, the product
    \begin{align*}
        \mu':=(-1)^{\frac{c(c+1)}{2}}\mu
    \end{align*}
    and the coproduct
    \begin{align*}
        \nu':= (-1)^{\frac{d(d+1)}{2}}\nu
    \end{align*}
    together with the same unit and counit satisfy the signs in their relations.
\end{Rem}

\begin{Cor}\label{cor:no comulti}
    Let $(\mathcal{C},\otimes, \mathbbm{1},\tau)$ be a monoidal category enriched over $\Zgrmod_0$ and let $c,d\in \Z$.\\
    Let the data $(A,\mu,\eta,\varepsilon)$ be as in Theorem \ref{thm:gen and rel}. This data uniquely defines a $(c,d)$-graded Frobenius algebra if and only if it satisfies graded associativity as in \eqref{item:thm ass}, graded unitality as in \eqref{item:thm unit} and the map
    \begin{align*}
        A\otimes A\overset{\mu}{\to} A\overset{\varepsilon}{\to} \mathbbm{1}
    \end{align*}
    is a non-degenerate pairing.
\end{Cor}
\begin{proof}
    The proof strategy is analogous to the ungraded situation as in \cite[Section 2.4]{kock2004frobenius}:\\
    We first use Theorem \ref{thm:gen and rel}. To show that a Frobenius algebra $(A,\mu,\eta,\nu,\varepsilon)$ as in the theorem induces a non-degenerate pairing $p:=\varepsilon\circ \mu$, one can show that $q:=(-1)^{cd+\frac{c(c+1)}{2}+\frac{d(d+1)}{2}}\nu\circ \eta$ is the corresponding copairing that satisfies
    \begin{align*}
        (\id\otimes p)\circ(q\otimes \id)=\id.
    \end{align*}
    On the other hand, we assume that the data $(A,\mu,\eta,\varepsilon)$ is as in the statement. Let $q\colon \mathbbm{1}\to A\otimes A$ be the associated copairing. Then one can show that the comultiplication $\nu:=(-1)^{\frac{d(d+1)}{2}}(\id\otimes \mu)\circ(q\otimes \id)$ satisfies coassociativity as in \eqref{item:thm coass} and counitality as in \eqref{item:thm counit}.
\end{proof}

\begin{Rem}
    Another perspective on graded Frobenius algebras is the following. Assume, we are given a degree $0$ multiplication $\mu$ and a non-degenerate pairing $p$ with associated copairing $q$ in degree $-d$ and $d$. Then one can use the language of oriented marked ribbon quivers as in \cite[Section 6]{kontsevich2021pre}. In this language, one can define a coproduct as the following quiver with a $0$-orientation:
    \begin{figure}[H]
        \centering
        \includegraphics[width=0.5\linewidth]{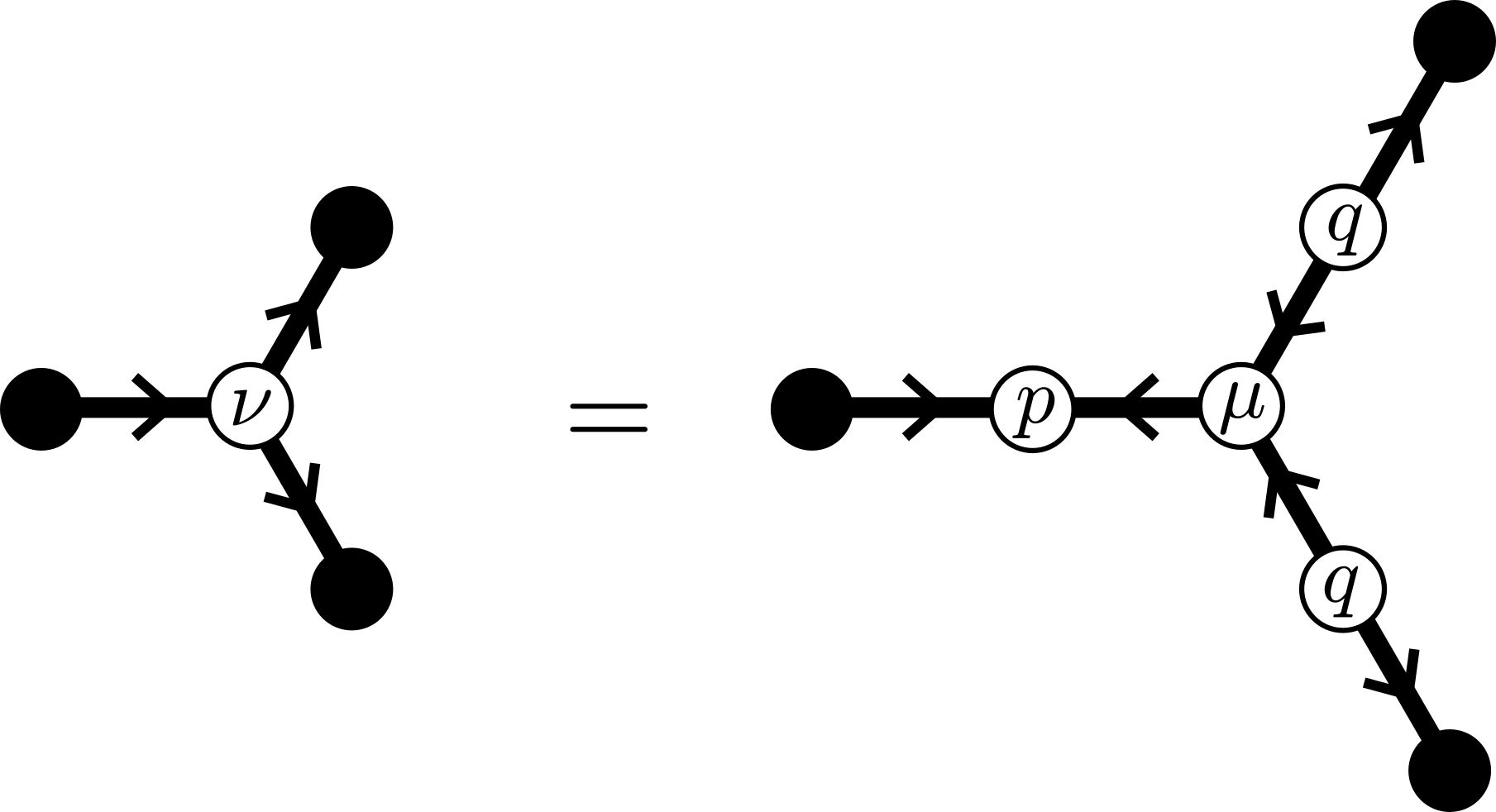}
    \end{figure}
    The signs in coassociativity and the Frobenius relation then come from permutations in the orientations.
\end{Rem}

We can now study the simplest non-trivial $\GrCob_{c,d}$-algebras. The simplest example is given for $c=d=0$:

\begin{Exmp}\label{exmp:unit is algabra}
    In $\mathcal{C}$, the unit $\mathbbm{1}$ has a canonical structure of a $\GrCob_{0,0}$-algebra given by the canonical maps
    \begin{align*}
        \mu\colon \mathbbm{1}\otimes \mathbbm{1}&\to \mathbbm{1}, & \eta\colon \mathbbm{1}&\to \mathbbm{1},\\
        \nu \colon \mathbbm{1}&\to \mathbbm{1}\otimes \mathbbm{1}, & \varepsilon\colon \mathbbm{1}&\to \mathbbm{1}.
    \end{align*}
\end{Exmp}

\begin{Exmp}\label{exmp:simplest algebra}
    Let $A$ be a $\GrCob_{c,d}$-algebra in $\Rgrmod_\bullet$ for some principal ideal domain $R$. We have unit maps $R\to A_{-c}$ and counit map $A_d\to R$. If $A$ is non-zero, non-degeneracy implies thus that we have a factor $R$ in degree $-c$ generated by $x:=\eta(1)$ and a factor $R$ in degree $d$ generated by $y$ with $\varepsilon(y)=1_R$.\\
    We distinguish between two cases: $c+d=0$ and $c+d\neq 0$.\\
    For $c+d=0$, we apply Proposition \ref{prop:suspending c,d} to Example \ref{exmp:unit is algabra} to find that $\Sigma^{d}\mathbbm{1}$ inherits the structure of a $\GrCob_{c,d}$-algebra. This is thus the $\GrCob_{c,d}$-algebra in $\Rgrmod_\bullet$ of smallest rank which is not zero. Moreover, it is the $\fatGrCob_{c,d}$-algebra of smallest rank bigger which is non-trivial.\\
    For $c+d\neq 0$, the algebra $A$ is not the suspension of a $\GrCob_{0,0}$-algebra. Moreover, every non-trivial $\GrCob_{c,d}$-algebra in $\Rgrmod_\bullet$ thus contains at least $\Sigma^{-c} R\oplus \Sigma^d R$. It turns out that we can define a $\GrCob_{c,d}$-algebra structure on $R_{c,d}:=\Sigma^{-c} R\oplus \Sigma^d R$.\\
    We denote the generator of $\Sigma^{-c} R$ by $x$ and the generator of $\Sigma^d R$ by $y$. By assumption, we have $\eta(1)=x$ and $\varepsilon(y)=1$. The multiplication $\mu$ is then uniquely defined by unitality
    \begin{align*}
        \mu(x\otimes y)&=(-1)^{\frac{c(c+1)}{2}}y, & \mu(y\otimes x)&= (-1)^{cd+\frac{c(c-1)}{2}}y, &\mu(x\otimes x)&=(-1)^{\frac{c(c+1)}{2}}x
    \end{align*}
    and degree reasons
    \begin{align*}
        \mu(y\otimes y)=0.
    \end{align*}
    On the other hand, the comultiplication $\nu$ is defined by counitality as
    \begin{align*}
        \nu(x)&=(-1)^{cd+\frac{d(d+1)}{2}}x\otimes y+ (-1)^{\frac{d(d-1)}{2}}y\otimes x, & \nu(y)&=(-1)^{\frac{d(d-1)}{2}}y\otimes y.
    \end{align*}
    One can check that this satisfies graded associativity, graded coassociativity, the graded Frobenius relation and graded commutativity.\\
    This is also the $\fatGrCob_{c,d}$- and $\pGrCob_{c,d}$-algebra of smallest rank which is non-trivial.
\end{Exmp}

\begin{Rem}\label{rem:one sided susp on algebras}
    As discussed in Remark \ref{rem:one sided suspension}, twisting by $\bolddet(G,\partial_{in})$ and $\bolddet(G,\partial_{out})$ changes only the degree of one operation, either multiplication or comultiplication while suspension changes degrees of both. On the level of algebras, suspension is simpler than twisting by $\bolddet(G,\partial_{in})$.\\
    We have an equivalence
    \begin{align*}
        \mathrm{grComFrob}_{c,d}(\Zgrmod_\bullet)\cong \mathrm{grComFrob}_{c-1,d+1}(\Zgrmod_\bullet)
    \end{align*}
    given by suspending the algebra. One can think of this in the following way. We tensor an algebra $A$ in $\mathrm{grComFrob}_{c,d}(\Zgrmod_\bullet)$ with $\Sigma \mathbbm{1}$ the $\GrCob_{-1,1}$-algebra of smallest rank which is not zero. We then apply Corollary \ref{cor:c+c' d+d'} to get a $\GrCob_{c-1,d+1}$-algebra structure on $\Sigma \otimes A\cong \Sigma A$.\\
    We can do an analogous construction to get a functor 
    \begin{align*}
        \mathrm{grComFrob}_{c,d}(\Zgrmod_\bullet)\to \mathrm{grComFrob}_{c,d+1}(\Zgrmod_\bullet).
    \end{align*}
    For this we tensor a $\GrCob_{c,d}$-algebra $A$ with $\Z_{0,1}$ as defined in Example \ref{exmp:simplest algebra}. This is the $\GrCob_{0,1}$-algebra of smallest rank which is not zero. By Corollary \ref{cor:c+c' d+d'}, this gives a $\GrCob_{c,d+1}$-algebra structure on $\Z_{0,1}\otimes A$ which has underlying module $A\oplus \Sigma A$.\\
    In contrast to suspension, this functor is not an equivalence. This can be seen as $\Z_{0,1}$ is not $\otimes$-invertible in contrast to $\GrCob_{1,-1}$-algebra $\Sigma \mathbbm{1}$.
\end{Rem}

The suspension of Frobenius algebras in $\Zgrmod_\bullet$ can also be studied on the level of maps and relations. Corollary \ref{cor:suspending GCob alg} shows that a shifted $\GrCob_{c,d}$-algebra is a $\GrCob_{c-1,d+1}$-algebra. This can also be stated in terms of maps and relations. However, it is a priori not clear that the suspended maps $\Sigma \mu,\Sigma\eta,\Sigma \nu$ and $\Sigma\eta$ satisfy the signs as in our main theorem because those signs depend on choices of orientations. It turns out that our choices are stable under suspending algebras in this sense.\\

To make this precise, we describe how a map  $\Zgrmod_\bullet(A^{\otimes k},A^{\otimes l})$ induces a suspended map in $\Zgrmod_\bullet((\Sigma A)^{\otimes k},(\Sigma A)^{\otimes l})$.\\
For this, we note that there is a canonical generator $\sigma_{k,l}$ in $\Zgrmod_\bullet((\Sigma \mathbbm{1})^{\otimes k},(\Sigma\mathbbm{1})^{\otimes l})$ which sends the canonical generator $s(1)\otimes \dots \otimes s(1)\in (\Sigma \mathbbm{1})^{\otimes k}$ to the canonical generator $s(1)\otimes \dots \otimes s(1)\in (\Sigma \mathbbm{1})^{\otimes l}$. We thus get a map
\begin{align*}
    \Sigma\colon \Zgrmod_\bullet(A^{\otimes k},A^{\otimes l})&\to \Zgrmod_\bullet((\Sigma \mathbbm{1})^{\otimes k},(\Sigma \mathbbm{1})^{\otimes l})\otimes \Zgrmod_\bullet(A^{\otimes k},A^{\otimes l})\cong \Zgrmod_\bullet((\Sigma A)^{\otimes k},(\Sigma A)^{\otimes l}),\\
    f&\mapsto \sigma_{k,l}\otimes f.
\end{align*}
We note that this is coherent with our choice that maps (here tensoring with $\sigma_{k,l}$) act from the left. However, this map is functorial and monoidal \textit{only up to sign}:
\begin{Lem}\label{lem:Sigma signs}
    Let $f\colon A^{\otimes k}\to A^{\otimes l}$ and $g\colon A^{\otimes m}\to A^{\otimes n}$ be maps in $\Zgrmod_\bullet$. Then we have
    \begin{enumerate}[(a)]
        \item\label{eq:tensor sign} $\Sigma f\otimes \Sigma g=(-1)^{(n-m)(|f|+k)}\Sigma (f\otimes g)$;
        \item\label{eq:comp sign} $\Sigma f\circ \Sigma g=(-1)^{(n-m)|f|}\Sigma(f\circ g)$ if $n=k$;
        \item\label{eq:id sign} $\Sigma \id_{A^{\otimes k}}=\id_{(\Sigma A)^{\otimes k}}$;
        \item\label{eq:twist sign} $\Sigma\tau_A=-\tau_{\Sigma A}$ for $\tau_A\colon A\otimes A\to A\otimes A$ the twist induced by the Koszul sign.
    \end{enumerate}
\end{Lem}
\begin{proof}
    We first note that the maps $\sigma_{k,l}$ are in degree $l-k$ and satisfy
    \begin{align*}
        \sigma_{k,l}\otimes \sigma_{m,n}(s(1)^{\otimes k}\otimes s(1)^{\otimes m})=(-1)^{(n-m)k}s(1)^{\otimes l}\otimes s(1)^{\otimes n}=(-1)^{(n-m)k}\sigma_{k+m,l+n}(s(1)^{\otimes k+m})
    \end{align*}
    and 
    \begin{align*}
        \sigma_{k,l}\circ \sigma_{m,k}=\sigma_{m,l}.
    \end{align*}
    We can thus compute 
    \begin{align*}
        \Sigma f\otimes \Sigma g&=(\sigma_{k,l}\otimes f)\otimes (\sigma_{m,n}\otimes g)=(-1)^{(n-m)|f|}(\sigma_{k,l}\otimes \sigma_{m,n})\otimes (f\otimes g)\\
        &=(-1)^{(n-m)(|f|+k)}(\sigma_{k+m,l+n})\otimes (f\otimes g)=(-1)^{(n-m)(|f|+k)}\Sigma(f\otimes g)
    \end{align*}
    and if $n=k$
    \begin{align*}
        \Sigma f\circ \Sigma g&=(\sigma_{k,l}\otimes f)\circ (\sigma_{m,n}\otimes g)=(-1)^{|\sigma_{m,n}||f|}(\sigma_{k,l}\circ \sigma_{m,n})\otimes (f\circ g)\\
        &=(-1)^{(n-m)|f|}\sigma_{m,l}\otimes (f\circ g)=(-1)^{(n-m)|f|}\Sigma(f\circ g).
    \end{align*}
    The equation \eqref{eq:id sign} follows directly from \eqref{eq:comp sign}. For \eqref{eq:twist sign}, we denote $s(a)\in (\Sigma A)_{i+1}$ for the element coming from $a\in A_{i}$ and similarly $s(b)\in (\Sigma A)_{j+1}$. The map $\Sigma \tau_A$ is defined as the composition
    \begingroup
    \setlength{\tabcolsep}{0pt}
    \begin{center}
        \begin{tabular}{rcl}
            $s(a)\otimes s(b)$ & $\mapsto$ & $(-1)^{|a|}s(1)\otimes s(1)\otimes a\otimes b$\\
            & $\overset{\sigma_{2,2}\otimes \tau_A}{\mapsto}$&$(-1)^{|a||b|+|a|}(s(1)\otimes s(1))\otimes (b\otimes a)$\\
            & $\mapsto$&$(-1)^{|a||b|+|a|+|b|}s(b)\otimes s(a)$\\
            &$=$&$(-1)^{|s(a)||s(b)|-1}s(b)\otimes s(a)=-\tau_{\Sigma A}(s(a)\otimes s(b))$.
        \end{tabular}
    \end{center}
    \endgroup
\end{proof}

\begin{Prop}\label{prop:shift maps}
    Let $(A,\mu,\eta,\nu,\varepsilon)$ in $\Zgrmod_\bullet$ be as in Theorem \ref{thm:gen and rel} and $c,d\in \Z$.\\
    This data satisfies the relations of a $(c,d)$-graded (commutative or symmetric) Frobenius algebra if and only if $(\Sigma A, \Sigma\mu,\Sigma\eta,\Sigma\nu,\Sigma\varepsilon)$ satisfies the relations of a $(c-1,d+1)$-graded (commutative or symmetric) Frobenius algebra.
\end{Prop}
\begin{proof}
    This is a direct application of Lemma \ref{lem:Sigma signs}. For example, if $(A,\mu,\eta,\nu,\varepsilon)$ satisfies $(c,d)$-graded left unitality, then we check compute $(c-1,d+1)$-graded left unitality for $(\Sigma A,\Sigma\mu,\Sigma\eta,\Sigma\nu,\Sigma\varepsilon)$:
    \begin{align*}
        (-1)^{c-1}\Sigma\mu\circ (\Sigma\eta\otimes \id_{\Sigma A})\overset{\eqref{eq:id sign}}=&(-1)^{c-1}\Sigma\mu\circ (\Sigma\eta\otimes \Sigma\id_A)\overset{\eqref{eq:tensor sign}}{=}(-1)^{c-1}\Sigma\mu\circ \Sigma(\eta\otimes \id_A)\\
        \overset{\eqref{eq:comp sign}}=&(-1)^{c-1+c}\Sigma(\mu\circ (\eta\otimes \id_A))=(-1)^{c-1+\frac{c(c-1)}{2}}\Sigma\id_{A}\\
        \overset{\eqref{eq:id sign}}=&(-1)^{\frac{(c-1)(c-1-1)}{2}}\id_{\Sigma A}.
    \end{align*}
\end{proof}

\section{Examples}\label{sec:exmp}
In this section, we discuss different examples of graded Frobenius algebras found in the literature. Some of the examples have a description more akin to our Definition \ref{def:graded Frob alg} in terms of PROPs, while some of the examples have a description in terms of maps and relations similar to Theorem \ref{thm:gen and rel}.\\

Throughout this section, we use Sweedler notation of a coproduct $\nu$ on an element $a$:
\begin{align*}
    \nu(a)=\sum_{(a)} a'\otimes a''.
\end{align*}

\subsection{Cohomology of an Oriented Manifold}\label{subsec:coh(M)}

Let $M$ be a connected, closed manifold of dimension $d$ and let $K$ be a field. Denote by $H^*(M)$ the cohomology of $M$ with coefficients in $K$. Assume that $M$ is $K$-orientable and fix an orientation $[M]\in H_d(M)$.\\
We denote 
\begin{align*}
    \mu\colon H^*(M)\otimes H^*(M)&\to H^*(M),\\
    \alpha \otimes \beta &\mapsto \alpha \cup \beta
\end{align*} 
and
\begin{align*}
    \varepsilon\colon H^*(M) &\to K,\\
    \alpha &\mapsto \begin{cases}
        \alpha ([M]) &\text{if } *=d;\\
        0 &\text{else.}
    \end{cases}
\end{align*}
The multiplication $\mu$ has degree $0$ and $\varepsilon$ has degree $-d$. Moreover, $\mu$ is associative and has unit map
\begin{align*}
    \eta \colon K&\to H^*(M),\\
    1&\mapsto1\in H^0(M).
\end{align*}
Moreover, it is graded commutative as it satisfies $\mu\circ\tau=\mu$. Poincaré duality implies that the pairing $\varepsilon\circ \mu$ is non-degenerate. We can thus apply Corollary \ref{cor:no comulti} and get the following example
\begin{Exmp}\label{exmp:coh}
    The data $(H^*(M),\mu,\eta,\varepsilon)$ defines a $(0,d)$-graded commutative Frobenius algebra
\end{Exmp}
The associated comultiplication $\nu$ is the intersection coproduct. There are two common constructions of the intersection coproduct leading to different signs: we refer to them as the Poincaré coproduct $\nu_P$ and the Thom coproduct $\nu_{Th}$ (see e.g.\ \cite[Appendix B]{hingston2017product}).\\
We first define the Poincaré coproduct. As we work with field coefficients, the Künneth maps
\begin{align*}
    \times \colon H_*(M)\otimes H_*(M)&\to H_*(M^2),\\
    \times \colon H^*(M)\otimes H^*(M) &\to H^*(M^2)
\end{align*}
are isomorphisms. Denote by $\theta$ their inverses. We thus have a coproduct on \textit{homology} given by
\begin{align*}
    \nu_{Hom}:=\theta\circ \Delta_*\colon H_*(M)\to H_*(M^2)\to H_*(M)\otimes H_*(M).
\end{align*}
We denote the Poincaré duality by 
\begin{align*}
    PD\colon H^*(M)&\overset{\cong}{\to} H_{d-*}(M),\\
    \alpha &\mapsto \alpha \cap [M].
\end{align*}
We then define the Poincaré coproduct as Poincaré dual to the homology coproduct $\nu_{Hom}$
\begin{align*}
    \nu_{P}(\alpha)=\sum_{(PD(\alpha))} PD^{-1}(PD(\alpha)')\otimes PD^{-1}(PD(\alpha)'')
\end{align*}
where the sum comes from $\nu_{Hom}$. This formula in terms of elements seems rather straight forward. But if we describe it as a composition of maps, we need the tensor product without Koszul sign: for graded maps $f\colon V\to V'$ and $g\colon W\to W'$ denote $f\wotimes g$ as the as the map
\begin{align*}
    (f\wotimes g)(a\otimes b)=f(a)\otimes g(b).
\end{align*}
Then we write the Poincaré coproduct as
\begin{align*}
    \nu_P:= (PD^{-1}\wotimes PD^{-1})\circ\theta\circ \Delta_*\circ PD.
\end{align*}

The Thom coproduct is defined via an orientation $[M^2]$ of the manifold $M^2$. To get the signs as in our main theorem, we define $[M^2]:=(-1)^{\frac{n(n-1)}{2}}[M]\times [M]$.\\
This sign convention might seem surprising. It reflects the fact that we are shuffling the $i$-th dimension of the first copy of $M$ to be the $2i-1$-th dimension of $M^2$ and the $i$-th dimension of the second copy of $M$ to be the $2i$-th dimension of $M^2$. We call such a shuffle of two ordered sets of the same size a\textit{ riffle shuffle}.\\
In terms of coordinates, take a ordered basis $(v_1,\dots,v_n)$ of $\R^n$. We then order the basis
\begin{align*}
    \{(v_1,0),\dots, (v_n,0),(0,v_1),\dots,(0,v_n)\}
\end{align*}
of $\R^n\oplus \R^n$ as 
\begin{align*}
    ((v_1,0),(0,v_1),(v_2,0),\dots,(0,v_n)).
\end{align*}
This choice has the advantage that for two oriented manifolds $M,N$ of dimension $d,d'$, we have $[(M\times N)^2]=[M^2]\times [N^2]$. Corollary \ref{cor:c+c' d+d'} then gives that $H^*(M)\otimes H^*(N)$ inherits the structure of a $\GrCob_{0,d+d'}$-algebra. Our choice of signs for $\nu_{Th}$ then ensures that 
\begin{align*}
    H^*(M)\otimes H^*(N) \to H^*(M\times N)
\end{align*}
is a map of coalgebras.\\
We consider the diagonal $\Delta\colon M\to M\times M$ and fix a tubular neighbourhood $U\subseteq M\times M$ of $\Delta(M)$. We denote $j\colon M\to U$ and $i\colon U\to M\times M$ for the canonical inclusions. We thus have $\Delta=i\circ j$.\\
We define $[U]\in H_{2d}(U,\partial U)$ as the image of $[M^2]$ under the map
\begin{align*}
      H_{2d}(M^2)\to H_{2d}(M^2,\Delta(M)^c)\overset{(i_*)^{-1}}{\cong}H_{2d}(U,\partial U).
\end{align*}
Let $\tau\in H^d(U,\partial U)$ be the Thom class that satisfies
\begin{align*}
    \tau \cap [U]=j_*([M]).
\end{align*}
Now, we can define the Thom coproduct as the following composition
\begin{align*}
    \nu_{Th}\colon H^*(M)\cong H^*(U)\overset{\tau \cup}{\to} H^{*+d}(U,\partial U) \cong H^{*+d}(M^2,\Delta(M)^c)\to H^{*+d}(M^2)\cong(H^*(M)^{\otimes 2})^{*+d}.
\end{align*}
We point the readers attention to the fact that this is coherent with our convention that maps act from the left as we are cupping with $\tau$ from the left in the second map.\\

To compare the signs that arise in the two definitions, we consider this diagram
\begin{center}
   \begin{tikzcd}[row sep= 1cm]
       H^k(M)\arrow[dd,"\cap {[M]}"']  & H^k(U)\arrow[dd, "\cap (j_*{[M]})"', near start] \arrow[l, "j^*"', "\cong"] \arrow[r, "\cup \tau", bend right]\ar[r, phantom, "(-1)^{dk}"] \arrow[r, " \tau\cup ", bend left] & H^{k+d}(U,\partial U) \arrow[ldd, "\cap {[U]}"'] &H^{k+d}(M^2,\Delta(M)^c)\arrow[l, "\cong", "i^*"']\ar[d] \\
       & & H_{d-k}(M^2)\ar[d, equal] \ar[rd, phantom, "{(-1)^{\frac{d(d-1)}{2}}}"] &H^{k+d}(M^2)\arrow[d, equal] \ar[l, "{\cap [M^2]}"', "\cong "]\\
       H_{d-k}(M)\arrow[r,"j_*"]& H_{d-k}(U) \arrow[r,"i_*"] &H_{d-k}(M^2)\ar[rd, phantom, "{(-1)^{d+dk+di}}"] &H^{k+d}(M^2)\arrow[l, "{\cap ([M]\times [M])}"', "\cong"]\\
       & & (H_*(M)^{\otimes 2})_{d-k}\ar[u, "\times", "\cong"'] & (H^*(M)^{\otimes 2})^{k+d}. \arrow[l, "\cap {[M]\widetilde\otimes \cap [M]}"', "\cong"] \ar[u, "\times", "\cong"']
   \end{tikzcd}
\end{center}
\begin{itemize}
    \item The composition along the bottom left gives the Poincaré coproduct $\nu_P$.
    \item The composition along the top right gives the Thom coproduct $\nu_{Th}$.
    \item The rectangle on the left hand side and the hexagon on the top right commute due to the naturality of the cap product.
    \item The triangle in the centre commutes by the definition of the Thom class $\tau$. 
    \item The sign $(-1)^{dk}$ on top comes from graded commutativity of the cup product: $\alpha\cup \beta=(-1)^{|\alpha||\beta|}\beta\cup \alpha$.
    \item The sign $(-1)^{\frac{d(d-1)}{2}}$ in the rectangle on the right in the middle comes from our definition of $[M^2]$.
    \item The sign in the bottom right rectangle means that, on elements in $H^i(M)\otimes H^{d+k-i}(M)$, we have the sign $(-1)^{d+dk+di}$ because the cross product and the cap product commute:
    \begin{align*}
        (\alpha \cap [M])\times (\beta \cap [M])= (-1)^{d|\beta|}(\alpha \times \beta)\cap ([M]\times [M]).
    \end{align*}
\end{itemize}
Combining this, we get on the $(i,d+k-i)$-coordinate in the output the sign
\begin{align}\label{eq:Thom vs Poincare}
    \nu_{P}=(-1)^{\frac{d(d-1)}{2}+d+di}\nu_{Th}.
\end{align}
One can check that $(A,\mu,\eta,\nu_{Th},\varepsilon)$ satisfies the relations of a $(0,d)$-graded commutative Frobenius algebra as in Theorem \ref{thm:gen and rel}. Cieliebak and Oancea prove this with their slightly different sign conventions in \cite[Proposition 7.1]{cieliebak2022cofrobenius}.\\
The signs of the Poincaré coproduct can be explained in the following way: as discussed in Proposition \ref{prop:suspending c,d} if $A$ has the structure of a $(c,d)$-graded commutative Frobenius algebra then $\Sigma A$ inherits the structure of a $(c-1,d+1)$-graded commutative Frobenius algebra. In this example, this means that
\begin{align*}
    \Sigma^{-d} H^*(M)=H^{*+d}(M)
\end{align*}
inherits the structure of a $(d,0)$-graded Frobenius algebra. Then the coproduct one obtains is given on an element $s^{-d}(\alpha)\in H^{*+d}(M)$ by:
\begin{align*}
    \Sigma^d\nu_{Th}(s^{-d}(\alpha))&=(\sigma_{1,2}^{\otimes -d}\otimes\nu_{Th})(s(1)^{\otimes {-d}}\otimes \alpha)\\
    &=(-1)^{d}\sigma_{1,2}^{\otimes -d}(s(1)^{\otimes -d})\otimes \nu_{Th}(\alpha)\\
    &=(-1)^{d}\sigma_{1,2}(s(1))^{\otimes -d}\otimes \sum_{(\alpha)}\alpha'\otimes \alpha''\\ %No riffle shuffle here because we evaluate from the inside out
    &=(-1)^{d}(s(1)^{\otimes 2})^{\otimes -d}\otimes \sum_{(\alpha)}\alpha'\otimes \alpha''\\
    &= \sum_{(\alpha)}(-1)^{d+|\alpha'|+|s(\alpha')|+\dots +|s^{d-1}(\alpha')|} s^{-d}(\alpha')\otimes s^{-d}(\alpha'')\\
    &=\sum_{(\alpha)}(-1)^{d+d|\alpha'|+\frac{d(d-1)}{2}}s^{-d}(\alpha')\otimes s^{-d}(\alpha'').
\end{align*}
Here $\sigma_{1,2}^{\otimes -d}$ denotes the $d$-fold tensor product of the canonical map $\Sigma^{-1}\mathbbm{1}\to \Sigma^{-1}\mathbbm{1}\otimes \Sigma^{-1}\mathbbm{1}$.\\
These are precisely the signs, we observed in \eqref{eq:Thom vs Poincare}. We thus conclude $\Sigma^{-d}\nu_{Th}=\nu_{P}$. Therefore, the Poincaré coproduct is the suspended Thom coproduct on $\Sigma^{-d}H^*(M)\cong  H^{*+d}(M)$. In $H^{*+d}(M)$, the coproduct is a degree 0 map, which makes the signs appearing in the graded coassociativity and counitality considerably more comfortable. One pays for this comfort by getting less appealing signs for the associativity and unitality of the corresponding shifted multiplication.

\subsection{Hochschild Homology of a Frobenius Algebra}\label{subsec:Hochschild}

\begin{Def}
    Let $\GrCob_{c,d}^{\partial_{out}\neq \varnothing}$ be the subcategory of $\GrCob_{c,d}$ which only has morphisms coming from graphs $G$ satisfying the plumber's condition, i.e., all connected components of $G$ have at least one outgoing external vertex.\\
    Let $\GrCob_{c,d}^{\partial_{in}\neq \varnothing}$ be the subcategory of $\GrCob_{c,d}$ which only has morphisms coming from graphs $G$ satisfying the reverse plumber's condition, i.e., all connected components of $G$ have at least one incoming external vertex.
\end{Def}

In \cite{wahlwesterland}, Wahl and Westerland show that if $A$ is a $(0,d)$-graded symmetric Frobenius algebra, then the normalized Hochschild homology $\overline{HH}_*(A)$ and $A$ has the structure of an open-closed TQFT. Isolating the closed part, we get a TQFT twisted by $\bolddet(S,\partial_{out}S)^{\otimes d}$. We can restate this in our language.
\begin{Exmp}\label{exmp:hochschild}
    Let $A$ be a $(0,d)$-graded symmetric Frobenius algebra. Then $\overline{HH}_*(A)$ has the structure of a $\GrCob^{\partial_{in}\neq \varnothing}_{d,d}$-algebra.
\end{Exmp}
The key input for this is the following observation: the twisting $\bolddet(\Cob(G),\partial_{out}(\Cob(G)))^{\otimes d}$ corresponds to our $(d,d)$-twisting, $\bolddet(G,\partial_{in})^{\otimes d}\otimes \bolddet(G,\partial_{out})^{\otimes d}$. (Recall that $\Cob(G)$ is the corresponding closed surface cobordism represented by the graph $G$.)\\
To make this precise, we recall that $\bolddet_{1,1}$ is a functor from $\underline{\GrCob}(X,Y)$ which commutes with composition along $X$ and $Y$. This can be stated by saying that 
\begin{align*}
    \bolddet_{1,1}\colon \underline{\GrCob}&\to B\mathrm{Pic}(\Z),\\
    X&\mapsto *,\\
    G&\mapsto \bolddet_{1,1}(G)
\end{align*}
is a functor of $2$-categories. Here $B\mathrm{Pic}(\Z)$ denotes the $2$-category with one object $*$ and morphism category
\begin{align*}
    B\mathrm{Pic}(\Z)(*,*)=\mathrm{Pic}(\Z)
\end{align*}
the category of $\otimes$-invertible objects in $\Zgrmod_0$ with isomorphisms between them.\\
The twistings $\bolddet(\Cob(G),\partial_{in}(\Cob(G)))^{\otimes d}$ and $\bolddet(\Cob(G),\partial_{out}(\Cob(G)))^{\otimes d}$ satisfy analogous functoriality and gluing properties as $\bolddet_{1,1}$ and thus also define $2$-functors.\\
That these twistings correspond to our twisting by $\bolddet_{1,1}$ is thus described in the following proposition:
\begin{Prop}\label{prop:1closed=2open twists}
    The three functors of $2$-categories
    \begin{align*}
        \bolddet_{1,1},\bolddet(\mathrm{Cob}(-),\partial_{in}),\bolddet(\mathrm{Cob}(-),\partial_{out})\colon \underline{\GrCob}\to B\mathrm{Pic}(\Z)
    \end{align*}
    are isomorphic.\\
    Explicitly, there exist natural isomorphisms 
    \begin{align}
    \begin{split}\label{eq:in+out twist is closed in twist}
        \bolddet(G,\partial_{in})\otimes \bolddet(G,\partial_{out})&\cong \bolddet(\mathrm{Cob}(G),\partial_{in}(\mathrm{Cob}(G)))\\
        \bolddet(G,\partial_{in})\otimes \bolddet(G,\partial_{out})&\cong \bolddet(\mathrm{Cob}(G),\partial_{out}(\mathrm{Cob}(G)))
    \end{split}
    \end{align}
    for all $G\in \underline{\GrCob}(X,Y)$ which commute with gluing of graphs.
\end{Prop}
\begin{proof}
    We only prove the first isomorphism. The second one is proved analogously.\\
    Both the left and right hand side of \eqref{eq:in+out twist is closed in twist} represent functors out of $\underline{\GrCob}(X,Y)$:
    \begin{align*}
        &\bolddet(-,\partial_{in})\otimes \bolddet(-,\partial_{out}),& \bolddet(\Cob(-),\partial_{in}(\Cob(-))).
    \end{align*}
    After precomposition with the forgetful functor $\underline{\fatGrCob}(X,Y)\to \underline{\GrCob}(X,Y)$, they define functors out of $\underline{\fatGrCob}(X,Y)$. What we construct in the following is a natural isomorphism of the two functors out of $\underline{\fatGrCob}(X,Y)$.\\
    This lifts to a natural transformation of functors out of $\underline{\GrCob}(X,Y)$. This can be seen as all edge collapses can be lifted to an edge collapse in a fat graph. It remains to check that graph isomorphisms can be lifted. For this we note that all graph isomorphisms can be decomposed into a graph isomorphism that can be lifted to a fat graph isomorphism and a graph automorphism. Finally all graph automorphisms act trivially on the left and right hand side of \eqref{eq:in+out twist is closed in twist}. Therefore the following isomorphism can be lifted to a natural transformation of functors out of $\underline{\GrCob}(X,Y)$.
    
    Let $G\in \underline{\fatGrCob}(X, Y)$ be a fat graph. We construct isomorphisms which commute with gluing of graphs:
    \begin{align*}
        f_G\colon \bolddet(\fatCob(G),\partial_{out})&\overset{\cong}{\to} \bolddet(\fatCob(G),\partial_{in}\cup \partial_{free});\\
        g_G\colon \bolddet(\fatCob(G),\partial_{in})\otimes \bolddet(\fatCob(G),\partial_{in}\cup \partial_{free})&\overset{\cong}{\to} \bolddet(\Cob(G),\partial_{in}).
    \end{align*}
    We recall that $\partial_{free}$ is the part of the boundary which is neither in the incoming boundary nor in the outgoing boundary. Then the proposition follows by taking $g_G\circ (\id\otimes f_G)$.\\
    For $f_G$, we first note that $H_*(|G|,|L_{out}|)$ is a free graded $\Z$-module. Therefore any isomorphism
    \begin{align*}
        H_*(|G|,|L_{out}|)\cong (H_*(|G|,|L_{out}|))^{\lor} \cong H^*(|G|,|L_{out}|)
    \end{align*}
    gives an isomorphism 
    \begin{align*}
        \bolddet(G,\partial_{out})=\bolddet(H_*(|G|,|L_{out}|))\cong \bolddet(H^*(|G|,|L_{out}|)).
    \end{align*}
    Next, we note that the homotopy equivalence of pairs
    \begin{align*}
        (|G|,|L_{out}|)\simeq (\fatCob(G),\partial_{out})
    \end{align*}
    gives an isomorphism
    \begin{align*}
        \bolddet(H^*(|G|,|L_{out}|))\cong \bolddet(H^*(\fatCob(G),\partial_{out})).
    \end{align*}
    All those isomorphisms so far give an isomorphism
    \begin{align*}
        \bolddet(G,\partial_{out})\cong \bolddet(H^*(\fatCob(G),\partial_{out}))
    \end{align*}
    which commutes with gluing.\\
    For the next step, we use generalized Lefschetz duality (see for example \cite[Theorem 3.43]{hatcher2002algebraic}). This gives an isomorphism
    \begin{align*}
        H^*(\fatCob(G),\partial_{out})\cong H_{2-*}(\fatCob(G),\partial_{in}\cup \partial_{free})
    \end{align*}
    induced by capping with $ [\fatCob(G)]$, the orientation in $H_2(\fatCob(G),\partial\fatCob(G))$. As the fat graphs induce oriented open cobordisms, this isomorphism also commutes with gluing. We thus get an isomorphism
    \begin{align*}
        \bolddet(H^*(\fatCob(G),\partial_{out}))\cong \bolddet(H_{2-*}(\fatCob(G),\partial_{in}\cup \partial_{free}))\cong \bolddet(H_{*}(\fatCob(G),\partial_{in}\cup \partial_{free}))
    \end{align*}
    which commutes with gluing. This gives $f_G$.\\
    For $g_G$, we consider $\fatCob(G)$ as a subspace of $\Cob(G)$. Figuratively speaking, we consider it as the lower half. See Figure \ref{fig:Lefdual}.
    \begin{figure}[H]
        \centering
        \includegraphics[width=0.4\linewidth]{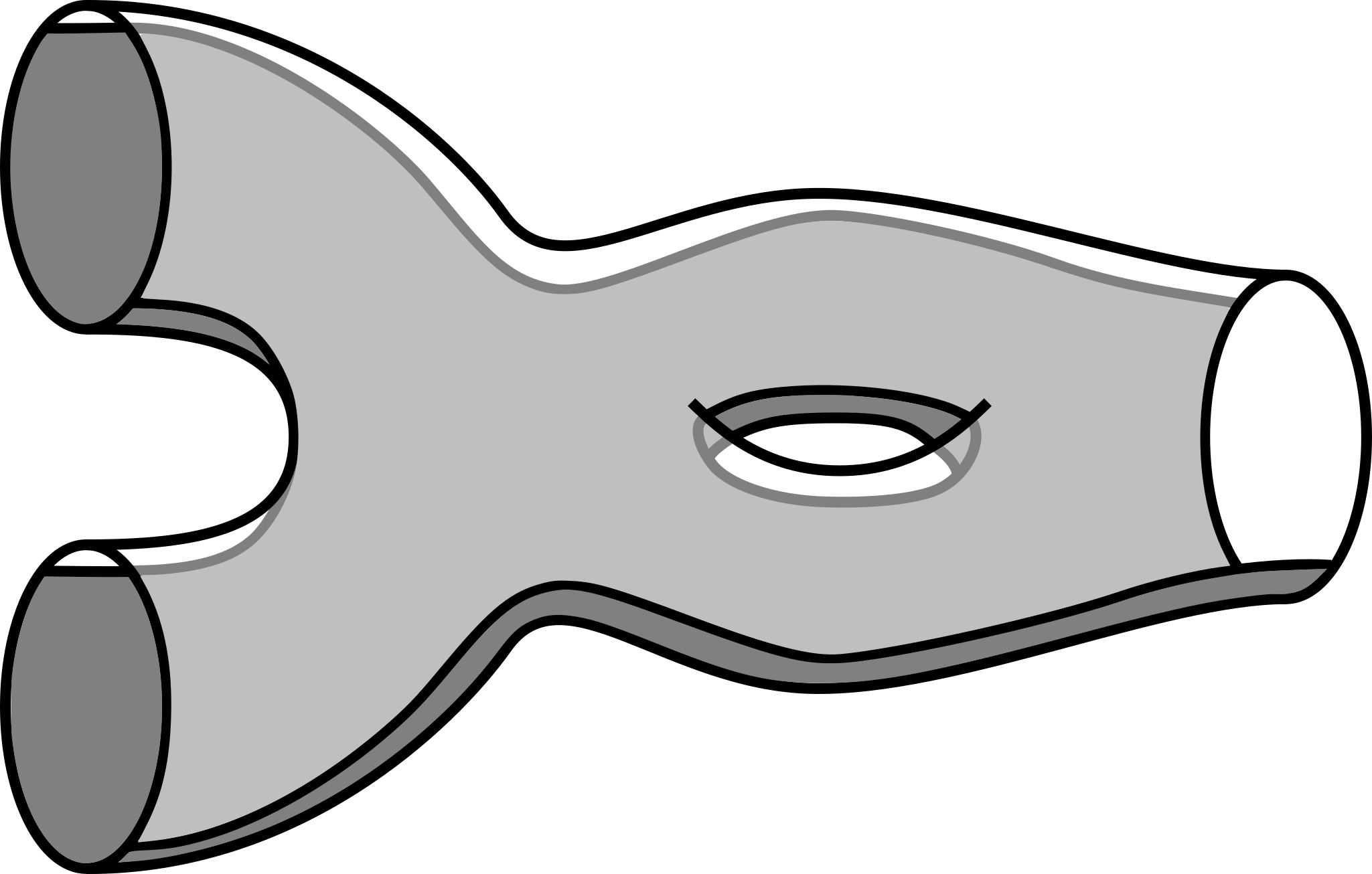}
        \caption{Embedding $\fatCob(G)$ in grey as the bottom half of $\Cob(G)$.}
        \label{fig:Lefdual}
    \end{figure}
    Consider the triple of spaces $(\Cob(G),\fatCob(G)\cup \partial_{in}(\Cob(G)),\partial_{in}(\Cob(G)))$. We note that by excision
    \begin{align*}
         H_*(\fatCob(G),\partial_{in}(\fatCob(G)))&\cong H_*(\fatCob(G)\cup \partial_{in}(\Cob(G)),\partial_{in}(\Cob(G)))\\
         H_*(\fatCob(G),\partial_{in}(\fatCob(G))\cup \partial_{free}(\fatCob(G)))&\cong H_*(\Cob(G),\fatCob(G)\cup\partial_{in}(\Cob(G))).
    \end{align*}
    We thus have
    \begin{align*}
        \bolddet(&\fatCob(G),\partial_{in}(\fatCob(G)))\otimes \bolddet(\fatCob(G),\partial_{in}(\fatCob(G))\cup \partial_{free}(\fatCob(G)))\\
        &\overset{\cong}{\to} \bolddet(\fatCob(G)\cup \partial_{in}(\Cob(G)),\partial_{in}(\Cob(G)))\otimes \bolddet(\Cob(G),\fatCob(G)\cup\partial_{in}(\Cob(G))).
    \end{align*}
    We note that collapsing the incoming boundary and the lower half in $\Cob(G)$ is the same space as collapsing the free and incoming boundary of the upper half. See Figure \ref{fig:Lefdual}.\\
    The short exact sequence of the triple of spaces and Lemma \ref{lem:det} \eqref{det:4} then gives 
    \begin{align*}
        \bolddet(&\fatCob(G)\cup \partial_{in}(\Cob(G)),\partial_{in}(\Cob(G)))\otimes \bolddet(\Cob(G),\fatCob(G)\cup\partial_{in}(\Cob(G)))\\
        &\overset{\cong}{\to}\bolddet(\Cob(G),\partial_{in}(\Cob(G))).
    \end{align*}
    This gives $g_G$ and thus concludes the proof.
\end{proof}

We note that the category $\GrCob_{d,d}^{\partial_{in}\neq \varnothing}$ is generated by $\omega(\multi)$, $\omega(\comulti)$ and $\omega(\counit)$. A symmetric monoidal functor $F\colon \GrCob_{d,d}^{\partial_{in}\neq \varnothing}\to \Z\mathrm{Mod}_\bullet$ is uniquely determined up to natural isomorphism by an object $F(*)$ and
\begin{align*}
    \mu_{HH}:=F(\omega(\multi))&\in \Zgrmod_{d}(F(*)\otimes F(*),F(*)),\\
    \nu_{HH}:=F(\omega(\comulti)) &\in \Zgrmod_d(F(*),F(*)\otimes F(*)),\\ \varepsilon_{HH}:=F(\omega(\counit))&\in \Zgrmod_{-d}(F(*),\mathbbm{1}).
\end{align*}
Those maps satisfy graded associativity, coassociativity, counitality and commutativity and the graded Frobenius relation as in Theorem \ref{thm:gen and rel}.\\
Following \cite[Section 6.5]{wahlwesterland} using our sign conventions that maps act from the left, we have the description
\begin{align*}
    \mu_{HH}((a_0\otimes \dots \otimes a_k)\otimes  (b_0\otimes \dots \otimes b_l))=&\ \begin{cases}
        0 &\text{if }k>0,\\
        \sum_{(a_0)}(-1)^{(|a_0'|+d)|a_0''|} a_0''a_0'b_0\otimes b_1\otimes \dots \otimes b_l &\text{if } k=0,
    \end{cases}\\
    \nu_{HH}(a_0\otimes \dots \otimes a_k)=&\ \sum_{\substack{i=0,\dots, k\\ (a_0)}}(-1)^{(|a_0'|+k-i)(|a_0''|+|a_1|+\dots +|a_i|)}\\
    &(a_0''\otimes a_1\otimes  \dots \otimes a_i)\otimes (a_0'\otimes a_{i+1}\otimes \dots \otimes a_k),\\
    \varepsilon_{HH}(a_0\otimes \dots \otimes a_k)=&\ \begin{cases}
        0 &\text{if }k>0,\\
        \varepsilon(a_0) &\text{if }k=0.
    \end{cases}
\end{align*}
We note that these are different signs compared to \cite{wahlwesterland}. One reason is that they have a different sign convention that maps act from the right. On the other hand, the signs in their coproducts are missing the Koszul sign coming from moving $a_0'$ past $a_0''\otimes a_1\otimes \dots\otimes a_i$.

\subsection{Homology of the Loop Space}\label{subsec:loop bad coprod}
Let $M$ be a connected, closed, oriented manifold of dimension $d$ and let $K$ be a field. Let $LM=\Map(S^1,M)$ be its loop space.\\
In \cite{cohen2004polarized}, Cohen and Godin define operations $\mu_S$ on $H_*(LM):=H_*(LM;K)$ for every closed cobordism $S$ such that each connected component has at least one outgoing boundary. In \cite{godin2007higher}, Godin moreover explains that this fits into a positive boundary open-closed TQFT which is twisted by $\bolddet(S,\partial_{in}S)^{\otimes -d}$. This minus sign comes from our convention that the determinant is in negative degree while Godin has the convention that the determinant is in positive degree.\\
We only focus on the closed underived part of this TQFT. Using the first isomorphism in Proposition \ref{prop:1closed=2open twists}, we can rephrase this in our language as the following example:
\begin{Exmp}\label{exmp:LM w Bad coprod}
    The homology of the loop space $H_*(LM)$ has the structure of a $\GrCob_{-d,-d}^{\partial_{out}\neq \varnothing}$-algebra with operations given by the Chas-Sullivan product and the trivial coproduct.
\end{Exmp}
We thus get a multiplication on $H_*(LM)$, that we denote by $\mu_{Go}$, and a comultiplication on $H_*(LM)$, that we denote by $\nu_{Go}$.\\ 
We can compare this to the results of Tamanoi in \cite{tamanoi2010loop}. Tamanoi defines a coproduct
\begin{align*}
    \Psi\colon H_{*+d}(LM)\to H_*(LM\times LM).
\end{align*}
As we work over a field, we can assume that the Künneth map $\times\colon H_*(LM)\otimes H_*(LM)\to H_*(LM\times LM)$ is invertible. We therefore can denote  $\nu_{Ta}:=\times^{-1}\circ \Psi$ a comultiplication in our sense. Tamanoi moreover defines multiplications
\begin{align*}
    \mu_{Ta}\colon H_*(LM)\otimes H_*(LM)&\to H_*(LM),\\
    \cdot \colon H_*(LM)\otimes H_*(LM\times LM)&\to H_*(LM\times LM),\\
    \cdot \colon H_*(LM\times LM)\otimes H_*(LM)&\to H_*(LM\times LM)
\end{align*}
which are defined such that for $a,b,c\in H_*(LM)$, we have
\begin{align*}
    a\cdot (b\times c)&=\mu_{Ta}(a\otimes b)\times c\\
    (a\times b)\cdot c&=a\times \mu_{Ta}(b\otimes c).
\end{align*}
In \cite[Theorem 2.2]{tamanoi2010loop}, Tamanoi observes the following Frobenius relation
\begin{align*}
    \Psi(\mu_{Ta}(a\otimes b))=(-1)^{d(|a|-d)}a\cdot \Psi(b)=\Psi(a)\cdot b\in H_*(LM\times LM).
\end{align*}
We interpret everything as living in $H_{*+d}(LM)$ and $H_{*+d}(LM\times LM)$. Then the Künneth map is a degree $d$ map and its inverse a degree $-d$ map. We can thus rewrite the Frobenius relation as
\begin{align*}
    \nu_{Ta}\circ \mu_{Ta}=(\mu_{Ta}\otimes \id)\circ(\id \otimes \nu_{Ta})=(\id \otimes \mu_{Ta})\circ (\nu_{Ta}\otimes \id).
\end{align*}
This is exactly the relation, we expect for a the $\GrCob_{0,-2d}^{\partial_{out}\neq \varnothing}$-algebra structure on $H_{*+d}(LM)=\Sigma^{-d}H_*(LM)$. This structure is predicted by shifting the $\GrCob_{-d,-d}^{\partial_{out}\neq \varnothing}$-structure in Example \ref{exmp:LM w Bad coprod} and applying Proposition \ref{prop:suspending c,d}.\\

It is worth noting that the coproducts, $\nu_{Go}$ and its shifted version $\nu_{Ta}$, are the so called trivial coproduct on $H_*(LM)$. This is in contrast to the Goresky-Hingston coproduct constructed in \cite{goresky2009loop,hingston2017product}.\\
Explicitly, Hingston and Wahl give a definition for the Chas-Sullivan product and the Goresky-Hingston product in \cite{hingston2017product}.\\
Their definition of the product is graded associative, unital and commutative. This lifts from the respective properties of the intersection product $\nu_{Th}$ discussed in \ref{subsec:coh(M)} (see the proof of \cite[Theorem 2.5]{hingston2017product}). This is in line with Example \ref{exmp:LM w Bad coprod}.\\
For the coproduct, we slightly amend their definition. Hingston and Wahl define in \cite[Definition 1.4]{hingston2017product} the loop coproduct as the degree $1-n$ map given for a chain $A\in C_*(LM,M)$ by
\begin{align*}
    \mathrm{AW}(\mathrm{cut}(R_{\mathrm{GH}}([\tau_{\mathrm{GH}}\cap](A\times I)))).
\end{align*}
For a detailed explanation of this notation, see \cite{hingston2017product}. We note that the first map in this composition $\times I$ acts from the right while the second map $[\tau_{\mathrm{GH}}\cap]$ acts from the left.\\
As we chose the convention that maps act from the left, we define the map as the composition
\begin{align*}
    \lor(A):=\mathrm{AW}(\mathrm{cut}(R_{\mathrm{GH}}([\tau_{\mathrm{GH}}\cap](I\times A)))).
\end{align*}
An analogous proof strategy as for \cite[Theorem 2.14]{hingston2017product} shows that $\lor$ is a degree $1-n$ coproduct that satisfies the graded coassociativity and cocommutativity relation as in Theorem \ref{thm:gen and rel}.\\

However, the Chas-Sullivan product and the Goresky-Hingston coproduct do not satisfy a Frobenius relation in this setting. In particular, the product is an operation on $H_*(LM)$ while the coproduct is an operation on $H_*(LM,M)$ the homology relative the constant loops.

\subsection{Rabinowitz Loop Homology}\label{subsec:Rabinowitz}

In the series of articles \cite{cieliebak2020poincar,cieliebak2022cofrobenius,cieliebak2022reduced,cieliebak2023loop,cieliebak2024rabinowitz}, Cieliebak, Hingston and Oancea give an interpretation where the Chas-Sullivan product and the Goresky-Hingston coproduct can be interpreted in a Frobenius type structure, what Cieliebak and Oancea call a biunital coFrobenius bialgebra. As discussed in Remark \ref{rem:alex kai}, the signs in Theorem \ref{thm:gen and rel} correspond to the signs appearing in the definition of biunital coFrobenius algebras. This gives a correspondence of their notion of a biunital coFrobenius bialgebras and our graded Frobenius algebras in the finite dimensional case. But for the key example of Rabinowitz loop homology, Cieliebak and Oancea need infinite dimensional algebras. For this they use the language of graded Tate vector spaces to get a Frobenius relation.\\
They give two tensor products on graded Tate vector spaces: $\widehat\otimes^*,\widehat{\otimes}^!$. All maps that appear in their relations have source $A^{\widehat\otimes^* k}$ and target $A^{\widehat\otimes^!l}$ for some $k,l\geq 0$. The framework of PROPs set up so far fails to describe those algebras. Indeed for a graded Tate vector space $A$, the collection
\begin{align*}
    \End_A(X,Y):=\Hom(A^{\widehat\otimes^* X},A^{\widehat\otimes^!Y})
\end{align*}
for $X,Y$ finite sets does not have a canonical structure of a PROP. This is because in general there is no canonical map $A^{\widehat\otimes^!X}\to A^{\widehat\otimes^* X} $ and thus no composition.\\
Our remedy for this problem is to work with dioperads and the following proposition:
\begin{Prop}\label{prop:tate dioperad}
    Let $A$ be a graded Tate vector space. The collection $\End_A(X,Y)$ has a canonical structure of a dioperad in $\Zgrmod_0$.
\end{Prop}
\begin{proof}
    We recall from Subsection \ref{subsec:dioperad} that a dioperad structure consists of a $(S_X,S_Y)$-bimodule structure, a unit and composition maps.\\
    The objects $\End_A(X,Y)$ have a canonical $(S_X,S_Y)$-bimodule structure given by permuting the coordinates.\\
    The unit is given by $\id_A\in \End_{A}(*,*)=\Hom(A,A)$.\\
    The composition map uses the map
    \begin{align*}
        \alpha \colon A\widehat{\otimes}^*(B\widehat\otimes^!C)\to (A\widehat\otimes^* B)\widehat\otimes^! C
    \end{align*}
    for graded Tate vector spaces $A,B,C$ as described in \cite[Proposition 5.5]{cieliebak2024rabinowitz}.\\
    We fix some sets $W,X,Y,Z$ and elements $x\in X$, $y\in Y$ and an element $f\in \End_A(W,X)$ and $g\in \End_A(Y,Z)$. We define $g {{}_y\circ_x} f\in \End_A((Y\setminus \{y\})\sqcup W, Z\sqcup (X\setminus \{x\}))$ as the composition
    \begin{align*}
        A^{\widehat{\otimes}^* (Y\setminus \{y\})\sqcup W}\cong&  (A^{\widehat\otimes^* Y\setminus \{y\}}) \widehat\otimes^*(A^{\widehat\otimes^* W}) \overset{f}{\to} (A^{\starotimes Y\setminus \{y\}})\widehat\otimes^*(A^{\widehat\otimes^! X})\\
        \overset{\alpha}{\to}& (A^{\starotimes Y})\widehat\otimes^! (A^{\widehat\otimes^! X\setminus \{x\}}) \overset{g}{\to}(A^{\shriekotimes Z})\shriekotimes (A^{\widehat\otimes^! X\setminus \{x\}})\cong A^{\shriekotimes Z\sqcup (X\setminus \{x\})}.
    \end{align*}
    Checking that this data satisfies the unitality, associativity and equivariance relations is analogous to checking that the endomorphism dioperad of an object in any symmetric monoidal category is a dioperad.
\end{proof}

This allows us to define a Tate vector space having the structure of an algebra over a dioperad:
\begin{Def}
    Let $\mathcal{D}$ a dioperad in $\Zgrmod_0$ and let $A$ be a Tate vector space. A \textit{$\mathcal{D}$-algebra structure} on $A$ is a morphism of dioperads
    \begin{align*}
        \mathcal{D}\to \End_A
    \end{align*}
    enriched over $\Zgrmod_0$, i.e.\ a morphism of dioperads which is linear and preserves degree.
\end{Def}

We recall that Proposition \ref{prop:subdioperad} shows that if $c,d$ are of different parity working with the forest dioperad does not lose operations compared to the PROP in characteristic $\neq 2$. The main example of Rabinowitz loop homology is such a case. We can thus state the closed part of \cite[Corollary 9.3 (c)]{cieliebak2022cofrobenius} in the following example.

\begin{Exmp}\label{exmp:LM w good coprod}
    Let $M$ be an oriented $d$-dimensional manifold. The Rabinowitz loop homology $\widehat{H}_*(LM)$ has the structure of a $\GrCob_{d,d+1}^{{\mathrm{forest}}}$-algebra.
\end{Exmp}

\begin{Rem}[Graded Open-Closed TQFTs]
    Cieliebak, Hingston and Oancea also use graded open-closed TQFTs. See for example \cite[Theorem 1.3.]{cieliebak2020poincar}. As far as we are aware, there does not exist a graph description of graded open-closed TQFTs which describes \cite[Theorem 1.3.]{cieliebak2020poincar}.\\
    One structure that describes some graded open-closed TQFTs is given in \cite[Subsection 6.3.]{wahlwesterland}. Taking the bottom homology of the PROP $\mathcal{OC}_d$ describes graded open-closed TQFTs. However, it has the two following restrictions: 
    \begin{itemize}
        \item $H_{bot}(\mathcal{OC}_d)$ does not describe a closed unit;
        \item the degree of the open multiplication is zero, while the degree of the open comultiplication, the closed multiplication and the closed comultiplication are all equal to a fixed $d\in \Z$. 
    \end{itemize}
    Especially, the second restriction means that this PROP does not describe \cite[Theorem 1.3.]{cieliebak2020poincar} even after possibly shifting the algebras.\\
    This leaves the question open whether there is a graph cobordism category describing a PROP or dioperad encoding the graded open-closed TQFTs as in \cite[Theorem 1.3.]{cieliebak2020poincar}. One promising outline of a construction may be the following. One constructs a dioperad analogous to $\GrCob_{c,d}^{\mathrm{forest}}$ but with two types of vertices. Vertices of \textit{open type}, that have cyclic ordering on $s^{-1}(v)$, and vertices of \textit{closed type}, that do not have a cyclic ordering. Moreover to make this well-defined, one stipulates that edges between vertices of different type may not be collapsed by graph morphisms. Such edges correspond to zippers and cozippers in \cite{lauda2008open, cieliebak2020poincar,cieliebak2022cofrobenius}.\\
    Then one can consider a twist by
    \begin{align*}
        \bolddet(G,\partial_{in}^o)^{c_o}\otimes \bolddet(G,\partial_{in}^c)^{c_c}\otimes \bolddet(G,\partial_{out}^o)^{d_o}\otimes \bolddet(G,\partial_{out}^c)^{d_c}
    \end{align*}
    where $\partial_{in}^o$ is the in-coming boundary of vertices of open type, $\partial_{in}^c$ is the in-coming boundary of vertices of closed type and so on. This is a dioperad with open and closed part whose multiplication and comultiplication may have any degree.\\
    However, an open-closed TQFT as in \cite{lauda2008open,cieliebak2020poincar,cieliebak2022cofrobenius} satisfies some relations involving the zipper and cozipper. For example, the zipper is an algebra morphism and the cozipper is a coalgebra morphism. These relations are not satisfied in the dioperad as we have described it. While it seems possible to force these relations in the dioperad, making this precise is beyond the scope of this article.
\end{Rem}

\section{Proof of Main Theorem}\label{sec:proof}
In this section, we introduce some computational tools which lets us more easily compute orientations in $\bolddet_{c,d}(G)$. With this machinery we can prove the main theorem in Subsection \ref{subsec:proof}.

\subsection{Computational Tools}\label{subsec:computational tools}
We recall that a generator $\omega(G)\in \bolddet_{c,d}(G)$ is called an orientation. For the proof, we need to be able to do three things:
\begin{enumerate}[(I)]
    \item\label{item:iso} compute how (fat) graph isomorphisms act on orientations;
    \item\label{item:edge collapse} compute how edge collapses act on orientations;
    \item\label{item:compo} compute the composition of two orientations.
\end{enumerate}
We recall that the composition in \eqref{item:compo} was defined in Lemma \ref{lem:def of compos} using a splitting of the short exact sequence:
\begin{center}
    \begin{tikzcd}
        0 \ar[r] & C_*(|G|,|\partial_{in}|) \ar[r] & C_*(|G'\circ G|,|\partial_{in}|)  \ar[r] & C_*(|G'|,|\partial_{in}|)   \ar[r] & 0.
    \end{tikzcd}
\end{center}
To explicitly compute \eqref{item:iso}, \eqref{item:edge collapse} and \eqref{item:compo}, we use Proposition \ref{prop:computedet} inspired by \cite[Proposition 4.14]{Getzler} and \cite[Lemma 14]{godin2007higher}. This gives us an expression of orientations in terms of top exterior powers of edges, half-edges and vertices:\\

For $X$ a finite set, we use the following notation
\begin{itemize}
    \item $\Or(X):=\Or(\Z^X)$ the top exterior concentrated in degree 0;
    \item $\det(X):=\det(\Z^X)$ the top exterior concentrated in degree $-|X|$.
\end{itemize}

\begin{Prop}\label{prop:computedet}
    Let $G$ be a graph. There is a natural isomorphism
    \begin{align*}
        \bolddet(G,\partial_{in}){\cong} &\det(E)^{-1}\otimes \Or(H)\otimes \det({V\setminus L_{in}}).
    \end{align*}
    Moreover, an analogous statement for $\bolddet(G,\partial_{out})$ holds after replacing $L_{in}$ by $L_{out}$.
\end{Prop}
\begin{proof}
    We set $L:=L_{in}$ and $\partial:=\partial_{in}$ or $L:=L_{out}$ and $\partial:=\partial_{out}$. Consider the cellular chain complex $C_*(|G|,|L|)$ for $|G|$ relative $|L|$:
    \begin{center}
        \begin{tikzcd}
            0 \ar[r] & \bigoplus_{e\in E} \Or(\{h_e,\sigma h_e\}) \ar[r] & \Z^V/\Z^{L} \ar[r] & 0,\\
             & h_e\wedge \sigma h_e\ar[r, mapsto] & s(\sigma h)-s(h).
        \end{tikzcd}
    \end{center}
    We have that the homology of this complex is $H_*(|G|,|L|)$. By Lemma \ref{lem:det}, we have a sequence of natural isomorphisms
    \begin{align*}
        \bolddet(G,\partial)\overset{\eqref{det:5}}{\cong} &\bolddet(C_*(G,L))\\
        \overset{\eqref{det:2}}{\cong}&\bolddet \left(\Sigma \bigoplus_{e\in E} \Or(\{h_e,\sigma h_e\})\right)\otimes \det(\Z^V/\Z^L)\\
        {\cong}\ &\bolddet \left(\Sigma \bigoplus_{e\in E} \Or(\{h_e,\sigma h_e\})\right)\otimes \det({V\setminus L})\\
        \overset{\eqref{det:3}}{\cong }&\det \left( \bigoplus_{e\in E} \Or(\{h_e,\sigma h_e\})\right)^{-1}\otimes \det({V\setminus L})\\
        \overset{\eqref{det:2}}{\cong} &\bigotimes_{e\in E}\det(\Or(\{h_e,\sigma h_e\}))^{-1}\otimes  \det({V\setminus L})
    \end{align*}
    We note that $\det(\Or(\{h_e,\sigma h_e\}))^{-1}$ is a free group of rank 1 concentrated in degree 1. Choosing a generator of this group is equivalent to choosing a generator in the degree 0 group $\Or(\{h_e,\sigma h_e\})$. We thus find 
    \begin{align*}
        \det(\Or(\{h_e,\sigma h_e\})^{-1}\cong \det(\{e\})^{-1}\otimes \Or(\{h_e,\sigma h_e\})
    \end{align*}
    where tensoring with $\det(\{e\})^{-1}$ ensures that the right hand side is in degree 1. We thus find
    \begin{align*}
              \bolddet(G,\partial) \cong &\ \bigotimes_{e\in E}( \det(\{e\})^{-1}\otimes \Or(\{h_e,\sigma h_e\}))\otimes \det({V\setminus L})\\
        \overset{\eqref{det:2}}{\cong} &\det(E)^{-1}\otimes \Or(H)\otimes  \det({V\setminus L})
    \end{align*}
    where the last isomorphism comes from the fact that $H=\sqcup_{e\in E} \{h_e,\sigma h_e\}$.
\end{proof}

\begin{Rem}
    Explicitly, for $E(G)=\{e_1,\dots,e_\alpha\}$, $e_i=\{h_i,\sigma h_i\}$ and $V(G)\setminus L_{in}(G)=\{v_1,\dots,v_\beta\}$ a generator of $\bolddet(G,\partial_{in})$ can be written as
    \begin{align}\label{eq:genofdet}
        (e_1\wedge \dots \wedge e_\alpha)^{{-1}} \otimes h_1\wedge \sigma h_1\wedge \dots \wedge \sigma h_\alpha\otimes v_1\wedge \dots \wedge v_\beta.
    \end{align}
\end{Rem}

A graph isomorphism induces a map of edges, half-edges and vertices. This lets us compute isomorphism on orientations which gives \eqref{item:iso} on our wish list. For computing edge collapses as in \eqref{item:edge collapse}, we use Lemma \ref{lem:edge collapse} which lets us compute an edge collapse on an orientation expressed in terms of edges, half-edges and vertices. For computing compositions as in \eqref{item:compo}, we use Lemma \ref{lem:compute compo}. The proofs of both lemmas are quite technical. But once they are proven doing \eqref{item:iso}, \eqref{item:edge collapse} and \eqref{item:compo} are straight-forward computations.

\begin{Lem}\label{lem:edge collapse}
    Let $G$ be a graph and $E=\{e_0,e_1,\dots, e_\alpha\}$, $e_i=\{h_i,\sigma h_i\}$ and $V=\{v_0,v_1,\dots,v_\beta\}$. Assume that $e_0$ is not a tadpole and $s(\sigma h_0)=v_0\notin L_{in}(G)$.\\
    The canonical map $c\colon |G|\to |G|/e_0$ does the following to generators as in \eqref{eq:genofdet}:
    \begin{align*}
        \bolddet(c)\colon \det({E(G)})^{-1}&\otimes \Or(H(G))\otimes \det({V(G)\setminus L_{in}(G)})\\
        &\to \det({E(G/e_0)})^{-1}\otimes \Or(H(G/e_0))\otimes \det({V(G/e_0)\setminus L_{in}(G/e_{0})})\\
        (e_0 \wedge e_1\wedge \dots \wedge e_\alpha)^{-1} &\otimes h_0\wedge \sigma h_0 \wedge h_1\wedge \dots \wedge \sigma h_\alpha\otimes v_0\wedge v_1\wedge \dots \wedge v_\beta\\
        &\mapsto (e_1\wedge \dots \wedge e_\alpha)^{-1} \otimes h_1\wedge \dots \wedge \sigma h_\alpha\otimes v_1\wedge \dots \wedge v_\beta.
    \end{align*} 
    Here we identify the newly created vertex in $G/e_0$ with $s(h_0)$ which gives an identification
    \begin{align*}
        V(G)\setminus \{v_0\}\cong V(G/e_0)
    \end{align*}
    together with the canonical identification
    \begin{align*}
        E(G)\setminus \{e_0\} &\cong E(G/e_0), &  H(G)\setminus \{h_0,\sigma h_0\} &\cong H(G/e_0).
    \end{align*}
    Moreover, an analogous statement for $\bolddet(G,\partial_{out})$ holds after replacing $L_{in}$ by $L_{out}$.
\end{Lem}
\begin{proof}
    We set $L:=L_{in}$ or $L:=L_{out}$ and $\partial:=\partial_{in}$ or $\partial:=\partial_{out}$. The map $c$ induces a map of CW-complexes. This in turn induces a map of cellular chain complexes
    \begin{center}
        \begin{tikzcd}
            0\ar[r] & \bigoplus_{e\in E(G)}\Or(e) \ar[r]\ar[d, "c|_E"] & \Z^{V(G)}/\Z^L \ar[r]\ar[d, "c|_V"]  & 0\\
            0\ar[r] & \bigoplus_{e\in E(G/e_0)}\Or(e) \ar[r] & \Z^{V(G/e_0)}/\Z^L\ar[r]   &0
        \end{tikzcd}
    \end{center}
    where $c|_E$ sends $\Or(e_0)$ to zero and all other edges to themselves and $c|_V$ sends $v_0$ to $v:=s(h_0)$ and all other vertices to themselves. This is a surjective map of chain complexes with kernel:
    \begin{center}
         \begin{tikzcd}
             0 \ar[r] & \Or(\{h_0,\sigma h_0\})\ar[r] & \langle v_0-v\rangle \ar[r] &0 \\
              & h_0\wedge \sigma h_0 \ar[r, mapsto] &v_0-v
         \end{tikzcd}
    \end{center}
    we denote this complex by $C$. Because $e_0$ is not a tadpole $v_0-v\neq 0$ and thus $C$ has trivial homology. Lemma \ref{lem:det} \eqref{det:5} gives the isomorphism 
    \begin{align*}
        \mathbbm{1}\cong \bolddet(C)\cong \det(\{e_0\})^{-1} \otimes\Or(e_0)\otimes \det(\langle v_0-v\rangle)
    \end{align*}
    given by 
    \begin{align}\label{eq:trivialization of edge complex}
        1\mapsto e_0^{-1}\otimes h_0\wedge \sigma h_0\otimes (v_0-v).
    \end{align}
    The geometric map $\bolddet(c)\colon \bolddet(G,\partial)\to \bolddet(G/e_{0},\partial)$ is given by the composition
    \begin{align*}
        \bolddet(H_*(G,\partial))\cong \bolddet(C)\otimes \bolddet(H_*(G/e_0,\partial))\cong \mathbbm1\otimes \bolddet(H_*(G/e_0,\partial))\cong \bolddet(H_*(G/e_0,\partial)).
    \end{align*}
    The isomorphism $\bolddet(H_*(G,\partial))\cong \bolddet(C)\otimes \bolddet(H_*(G/e_0,\partial))$ induces
    \begin{align*}
        \det&({E(G)})^{-1}\otimes \Or(H(G))\otimes \det(V(G)\setminus L) \\
        \to&\det({E(G/e_0)})^{-1}\otimes \det(\{ e_0\})^{-1}\otimes\Or(e_0)\otimes\Or(H(G/e_0))\otimes \det(\{v_0-v\})\otimes \det({V(G/e_0)\setminus L})\\
        \to &\det(\{e_0\})^{-1}\otimes \Or(e_0)\otimes \det(\{v_0-v\}) \otimes \det({E(G/e_0)})^{-1}\otimes \Or(H(G/e_0))\otimes \det({V(G/e_0)\setminus L}).
    \end{align*}
    This isomorphism is given by the identifications 
    \begin{align*}
        E(G)=\{e_0\}\sqcup E(G/e_0),\: H(G)= e_0\sqcup H(G/e_0),\: \Z^{V(G)\setminus L}\cong \Z^{v_0-v}\oplus \Z^{V(G/e_0)\setminus L}
    \end{align*}
    and then permuting the factors. We recall that in our conventions for graded abelian groups, the canonical isomorphism between the dual of the tensor product and the tensor product of the dual is order reversing. Thus the map on edges is order reversing:
    \begin{align*}
        \det(E(G))^{-1}&\to \det(E(G/e_0))\otimes \det(\{e_0\})^{-1},\\
        (e_0\wedge e_1 \wedge \dots \wedge e_\alpha)^{-1} &\mapsto (e_1\wedge \dots \wedge e_\alpha)^{-1} \otimes e_0^{-1}.
    \end{align*}
    Concretely on generators, we thus get the map
    \begin{align*}
        &(e_0 \wedge e_1\wedge \dots \wedge e_\alpha)^{-1} \otimes h_0\wedge \sigma h_0 \wedge h_1\wedge \dots \wedge \sigma h_\alpha\otimes v_0\wedge v_1\wedge \dots \wedge v_\beta\\
        =&\ (e_0 \wedge e_1\wedge \dots \wedge e_\alpha)^{-1} \otimes h_0\wedge \sigma h_0 \wedge h_1\wedge \dots \wedge \sigma h_\alpha\otimes (v_0-v)\wedge v_1\wedge \dots \wedge v_\beta\\
        \mapsto&\ (e_1\wedge \dots \wedge e_\alpha)^{-1}\otimes e_0^{-1} \otimes h_0\wedge \sigma h_0 \otimes h_1\wedge \dots \wedge \sigma h_\alpha\otimes (v_0-v)\otimes v_1\wedge \dots \wedge v_\beta\\
        \mapsto&\  e_0^{-1}\otimes h_0\wedge \sigma h_0\otimes (v_0-v) \otimes (e_1\wedge \dots \wedge e_\alpha)^{-1} \otimes h_1\wedge \dots \wedge \sigma h_\alpha\otimes  v_1\wedge \dots \wedge v_\beta. 
    \end{align*}
    The isomorphism $\bolddet(C)\otimes \bolddet(H_*(G/e_0,\partial))\cong \bolddet(H_*(G/e_0,\partial))$ is now given by applying \eqref{eq:trivialization of edge complex}. Concretely on generators, we compute
    \begin{align*}
        e_0^{-1}&\otimes h_0\wedge \sigma h_0\otimes (v_0-v) \otimes (e_1\wedge \dots \wedge e_\alpha)^{-1} \otimes h_1\wedge \dots \wedge \sigma h_\alpha\otimes  v_1\wedge \dots \wedge v_\beta\\
        &\mapsto (e_1\wedge \dots \wedge e_\alpha)^{-1} \otimes h_1\wedge \dots \wedge \sigma h_\alpha\otimes  v_1\wedge \dots \wedge v_\beta.
    \end{align*}
    Composing these two maps thus gives the equation we wanted to compute.
\end{proof}

\begin{Lem}\label{lem:compute compo}
    Let $G,G'$ be composable graphs in $\underline{\GrCob}$. The composition map
    \begin{align*}
        \mathrm{comp}\colon\bolddet(G',\partial_{in})\otimes \bolddet(G,\partial_{in})\cong& \bolddet(G'\circ G,\partial_{in})
    \end{align*}
    given in Lemma \ref{lem:def of compos} is computed via
    \begin{align*}
        \det(E(G'))^{-1}\otimes \Or(H(G'))&\otimes \det({V(G')\setminus L_{in}(G')})\\
        \otimes \det(E(G))^{-1}&\otimes \Or(H(G))\otimes \det({V(G)\setminus L_{in}(G)})\\
        &\cong \det(E(G'\circ G))^{-1}\otimes \Or(H(G'\circ G))\otimes \det({V(G'\circ G)\setminus L_{in}(G'\circ G)})
    \end{align*}
    where the map is given by the canonical identifications
    \begin{align}\label{eq:set ident}
    \begin{aligned}
        E(G')\sqcup E(G)&\cong E(G'\circ G),\\
        H(G')\sqcup H(G)&\cong H(G'\circ G),\\
        (V(G')\setminus L_{in}(G'))\sqcup (V(G)\setminus L_{in}(G))&\cong V(G'\circ G)\setminus L_{in}(G'\circ G).
    \end{aligned}
    \end{align}
    Moreover, an analogous statement holds for $\bolddet(G,\partial_{out})$ after replacing $L_{in}$ by $L_{out}$.
\end{Lem}

\begin{proof}
    We set $L:=L_{in}$ or $L:=L_{out}$ and $\partial:=\partial_{in}$ or $\partial:=\partial_{out}$.\\
    We consider the cellular chain complexes $C_*(G,L),C_*(G',L)$ and $C_*(G'\circ G,L)$. The map $\mathrm{comp}$ is given by Lemma \ref{lem:det} \eqref{det:4} applied to the short exact sequence of complexes
    \begin{align*}
        0\to C_*(G,L)\to C_*(G'\circ G,L)\to C_*(G',L)\to 0.
    \end{align*}
    One splitting of this short exact sequence is given by the maps on chains coming from the identification of the generating sets as in \eqref{eq:set ident}.
\end{proof}

\begin{Rem}\label{rem:sign compo}
    We write $E(G)=\{e_1,\dots,e_k\}$, $e_i=\{h_i,\sigma h_i\}$, $V(G)\setminus L_{in}(G)=\{v_1,\dots,v_l\}$, $E(G')=\{e_1',\dots,e_m'\}$, $e_i'=\{h_i',\sigma h_i'\}$ and $V(G')\setminus L_{in}(G')=\{v_1',\dots,v_n'\}$. Then we compute the following Koszul signs:
    \begin{align*}
        \mathrm{comp}(&(e_1'\wedge \dots \wedge e_m')^{{-1}} \otimes h_1'\wedge \sigma h_1'\wedge \dots \wedge \sigma h_m'\otimes v_1'\wedge \dots \wedge v_n'\\
        &\otimes (e_1\wedge \dots \wedge e_k)^{{-1}} \otimes h_1\wedge \sigma h_1\wedge \dots \wedge \sigma h_k\otimes v_1\wedge \dots \wedge v_l)\\
        =&\ (-1)^{(m+n)k}(e_1'\wedge\dots\wedge e_m'\wedge e_1\wedge \dots \wedge e_k)^{-1} \otimes h_1'\wedge\dots \wedge  \sigma h_m'\wedge h_1\wedge \dots \wedge \sigma h_k\\
        &\otimes v_1'\wedge \dots \wedge v_n' \wedge v_1\wedge \dots \wedge v_l
    \end{align*}
    in $\det(E(G'\circ G))^{-1}\otimes \Or(H(G'\circ G))\otimes \det({V(G'\circ G)\setminus L_{in}(G'\circ G)})$.
\end{Rem}

\subsection{Proof}\label{subsec:proof}
With the calculation tools in Proposition \ref{prop:computedet}, Lemma \ref{lem:edge collapse} and Lemma \ref{lem:compute compo}, we are in a position to prove Theorem \ref{thm:gen and rel}.

\begin{proof}[Proof of Theorem \ref{thm:gen and rel}]
    We recall that Theorem \ref{thm:gen and rel} gives a description of $(c,d)$-graded (commutative or symmetric) Frobenius algebras in terms of maps and relations that those maps satisfy. We defined a $(c,d)$-graded (commutative or symmetric) Frobenius algebra as a (symmetric) monoidal functor preserving degrees out of $\GrCob_{c,d}$ (or $\fatGrCob_{c,d}$ or $\pGrCob_{c,d}$).\\
    For the parts of the proof that apply to all three categories, we denote $\mathcal{G}_{c,d}$ for the category $\GrCob_{c,d}$, $\fatGrCob_{c,d}$ or $\pGrCob$. This is a category enriched in graded abelian groups. We denote by $\underline{\mathcal{G}}$ the corresponding $2$-category of graph cobordisms $\underline{\GrCob}$, $\underline{\fatGrCob}$ or $\underline{\pGrCob}$. Finally, we denote by $\mathcal{G}$ the associated $1$-categories $\GrCob$, $\fatGrCob$ or $\pGrCob$ where morphisms are equivalence classes of graph cobordisms as in Definition \ref{def:ass 1-cat}.\\
    
    We recall that Proposition \ref{prop:hom sets direct sum} gives isomorphisms
    \begin{align}
    \begin{split}\label{eq:maps as direct sum}
        {\mathcal{G}_{c,d}}(X,Y)&\cong\bigoplus_{[G]\in \mathcal{G}(X, Y)} {\bolddet}_{c,d}(G)_{\pi_1(|\underline{\mathcal{G}}(X,Y)|,G)}.
    \end{split}
    \end{align}
    An element $[G]\in \mathcal{G}(X,Y)$ is an equivalence class of a graph cobordism up to zig-zag of graph morphisms. Such an equivalence class can be represented by a graph that can be written as a composition of the four elementary graphs $\multi,\unit,\comulti,\counit$ (and possibly the twist map if $\mathcal{G}$ is $\GrCob$ or $\fatGrCob$). Indeed, expand the graph such that it has only graphs of valence $\leq 3$ and no two valence $3$-vertex are next to each other. Direct the edges such that no vertex of valence $2$ or $3$ has only incoming or only outgoing edges. To make this possible, one has to possibly introduce some units or counits. See Figure \ref{fig:graph decomp}.
    \begin{figure}[H]
        \centering
        \includegraphics[width=0.5\linewidth]{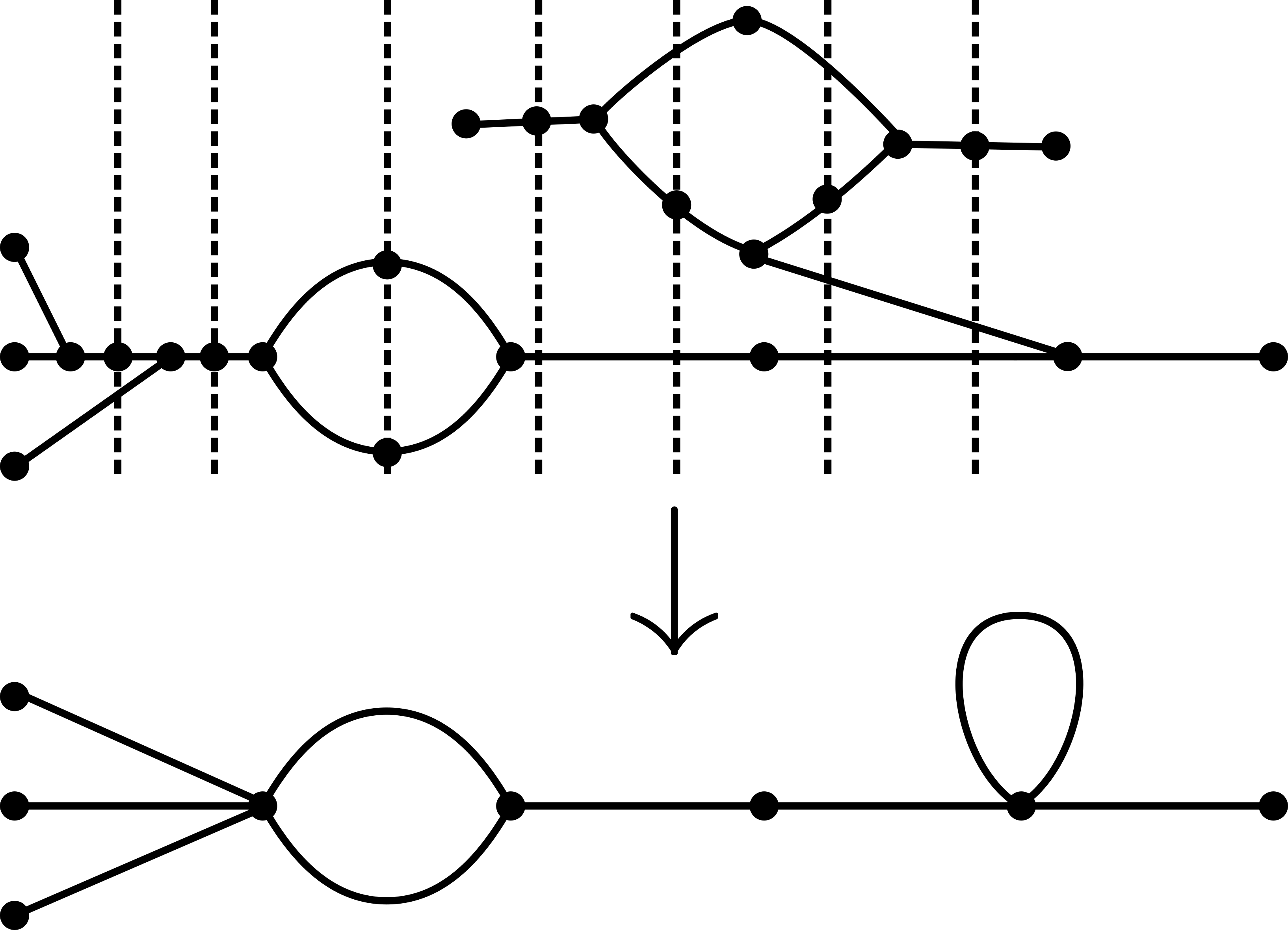}
        \caption{The graph on the bottom can be lifted to a graph which can be decomposed into the four elementary graphs. This graph can for example be written as $(\counit\otimes \multi) \circ (\multi \otimes \id \otimes \id)\circ  (\id\otimes \comulti\otimes \id)\circ (\comulti \otimes \id)\circ (\unit\otimes \multi)\circ \comulti\circ \multi\circ (\multi \otimes \id)$. (Note the order-reversing convention for the composition $\circ$.)}
        \label{fig:graph decomp}
    \end{figure}
    This shows that $\pGrCob_{c,d}$ is generated under composition and monoidal product by chosen generators
    \begin{align*}
        \omega(\multi)&\in {\bolddet}_{c,d}(\multi)_{\pi_1(|\underline{\mathcal{G}}(*\sqcup*,*)|,\multi)},& \omega(\unit)&\in  {\bolddet}_{c,d}(\unit)_{\pi_1(|\underline{\mathcal{G}}(\varnothing,*)|,\unit)},\\
        \omega (\comulti)&\in {\bolddet}_{c,d}(\comulti)_{\pi_1(|\underline{\mathcal{G}}(*,*\sqcup *)|,\comulti)},&\omega(\counit)&\in {\bolddet}_{c,d}(\counit)_{\pi_1(|\underline{\mathcal{G}}(*,\varnothing)|,\counit)}.
    \end{align*}
    The categories $\GrCob_{c,d}$ and $\fatGrCob_{c,d}$ are also generated by those four elements but as a \textit{symmetric} monoidal category rather than just as a monoidal category.\\
    We note that those four elements are non-trivial for all $c,d\in \Z$. This can be seen as there exists examples of $(c,d)$-graded commutative Frobenius algebras with non-trivial $\mu,\eta,\nu$ and $\varepsilon$ (see Example \ref{exmp:simplest algebra}). However, our proof does not use this fact.\\
    Any (symmetric) monoidal functor $F\colon \mathcal{G}_{c,d}\to \mathcal{C}$ is thus uniquely determined up to isomorphism by the images $\mu= F(\omega(\multi))$, $\eta= F(\omega(\unit))$, $\nu=F(\omega(\comulti))$ and $\varepsilon=F(\omega(\counit))$.\\

    We now fix explicit orientations $\omega(\multi)$, $\omega(\unit)$, $\omega(\comulti)$ and $\omega(\counit)$: we label the graphs $\multi$, $\unit$, $\comulti$ and $\counit$ as in Figures \ref{fig:gen of multi}, \ref{fig:gen of comulti}, \ref{fig:gen of unit} and \ref{fig:gen of counit}. We choose the convention that all half-edges named $\sigma h_i$ are glued to the internal vertex. This is convenient when we later apply Lemma \ref{lem:edge collapse} to compute edge collapses.
    \begin{figure}[H]
        \centering
        \begin{minipage}{0.45\textwidth}
            \centering
            \includegraphics[width=0.45\textwidth]{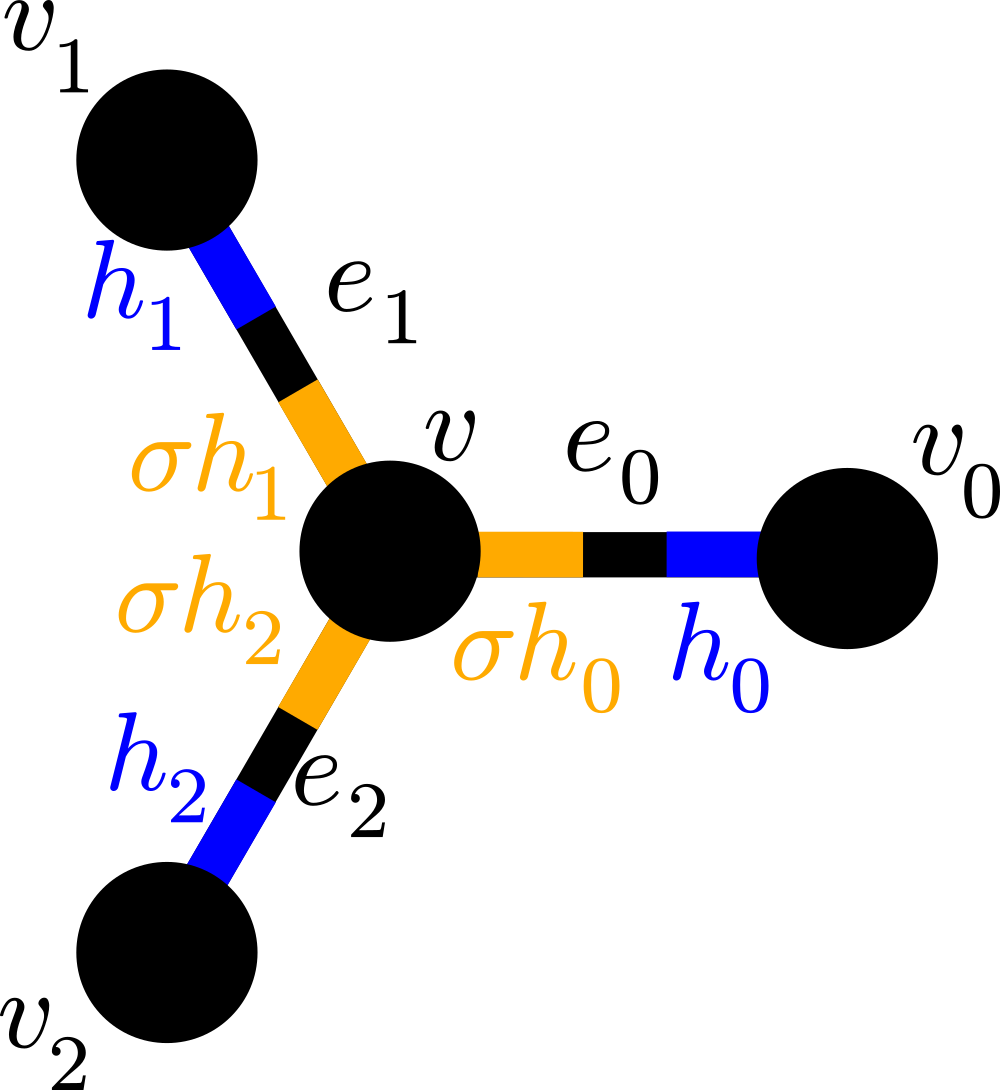} 
            \caption{The graph $\multi$ with $v$ the lone internal vertex, $v_0$ the lone outgoing vertex and $v_1$ and $v_2$ the incoming vertices. The vertex $v_1$ corresponding to the first coordinate of the input and $v_2$ to the second coordinate.}\label{fig:gen of multi}
        \end{minipage}\hfill
        \begin{minipage}{0.45\textwidth}
            \centering
            \includegraphics[width=0.45\textwidth]{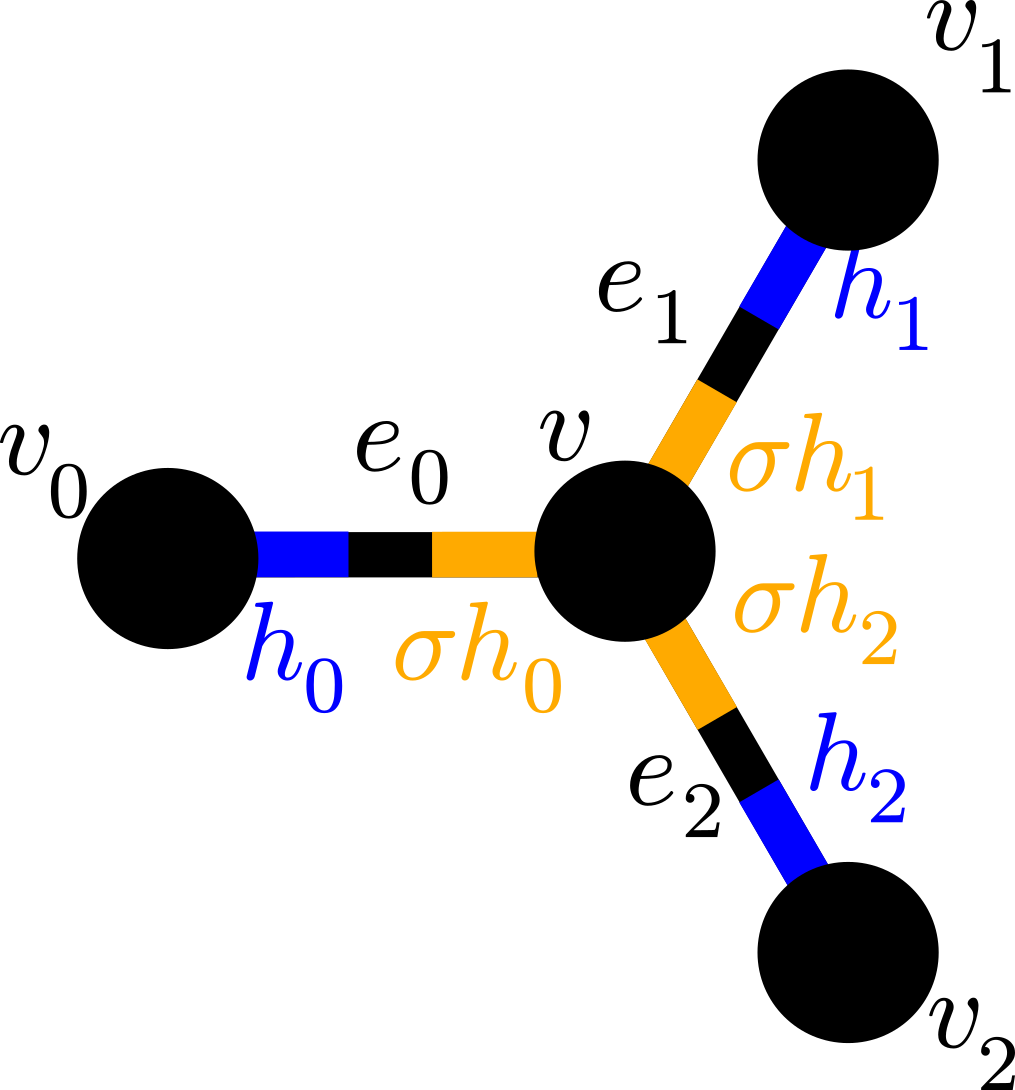} 
            \caption{The graph $\comulti$ with $v$ the lone internal vertex, $v_0$ the lone incoming vertex and $v_1$ and $v_2$ the outgoing vertices. The vertex $v_1$ corresponding to the first coordinate of the output and $v_2$ to the second coordinate.}\label{fig:gen of comulti}
        \end{minipage}
    \end{figure}
    \begin{figure}[H]
        \begin{minipage}{0.45\textwidth}
            \centering
            \includegraphics[width=0.3\textwidth]{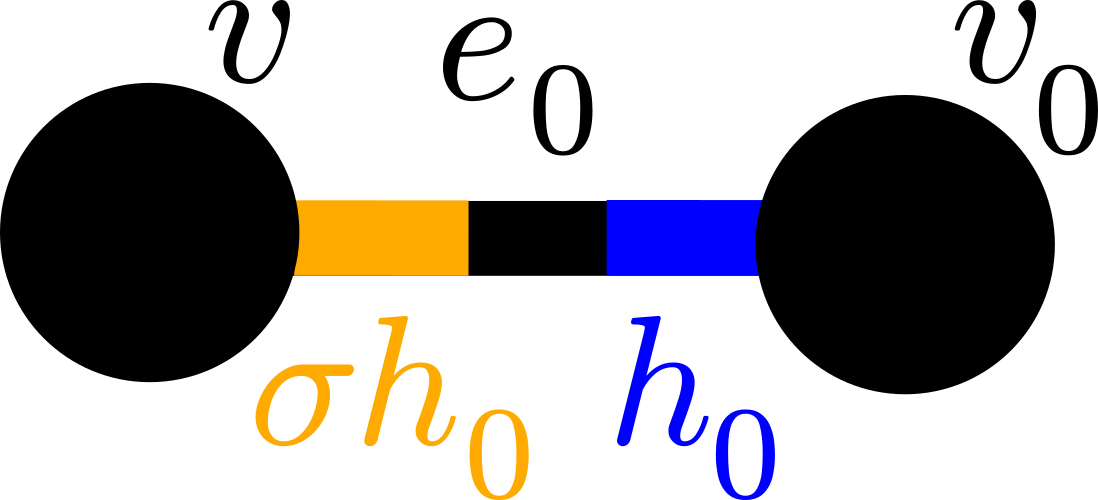}
            \caption{The graph $\unit$ with $v$ the lone internal vertex and $v_0$ the lone outgoing vertex.}\label{fig:gen of unit}
        \end{minipage}\hfill \begin{minipage}{0.45\textwidth}
            \centering
            \includegraphics[width=0.3\textwidth]{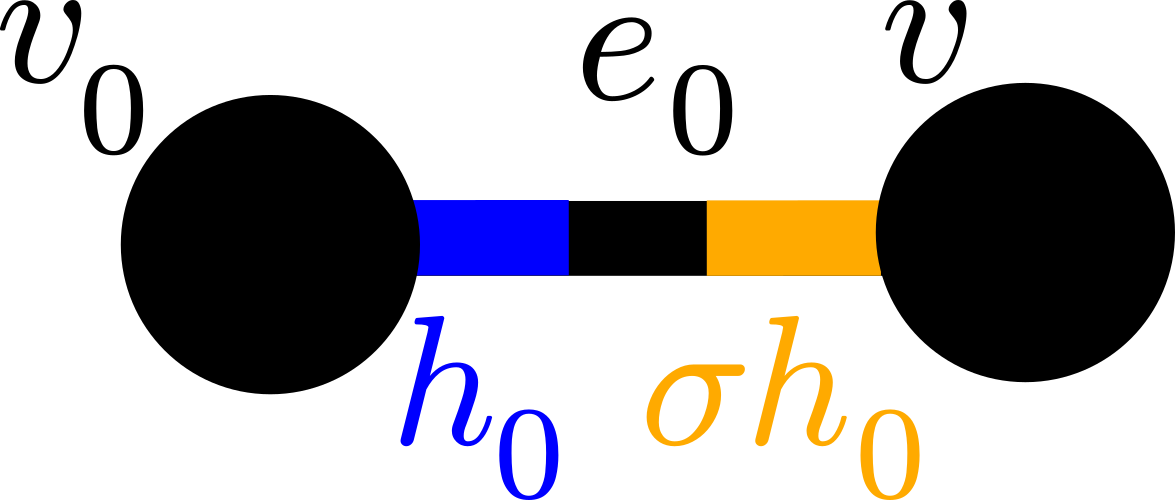}
            \caption{The graph $\counit$ with $v$ the lone internal vertex and $v_0$ the lone incoming vertex.}\label{fig:gen of counit}
        \end{minipage}
    \end{figure}
    Using the identification as in Proposition \ref{prop:computedet}, we define generators
    \begin{align*}
        \omega_{in}(\multi) &:=(e_2\wedge e_1\wedge e_0)^{-1}\otimes h_0\wedge \sigma h_0 \wedge h_1\wedge \sigma h_1 \wedge h_2\wedge \sigma h_2\otimes v\wedge v_0\in \bolddet(\multi,\partial_{in}),\\
        \omega_{out}(\multi)&:= 1\in \bolddet(\multi,\partial_{out})\cong \mathbbm{1},\\
        \omega_{in}(\unit) &:=e_0^{-1}\otimes h_0\wedge \sigma h_0 \otimes v\wedge v_0\in \bolddet(\unit,\partial_{in}),\\
        \omega_{out}(\unit)&:=1\in \bolddet(\unit,\partial_{out})\cong \mathbbm{1},\\
        \omega_{in}(\comulti) &:=1\in \bolddet(\comulti,\partial_{in})\cong \mathbbm{1},\\
        \omega_{out}(\comulti)&:=(e_2\wedge e_1\wedge e_0)^{-1}\otimes h_0\wedge \sigma h_0 \wedge h_1\wedge \sigma h_1 \wedge h_2\wedge \sigma h_2\otimes v\wedge v_0\in \bolddet(\comulti,\partial_{out}),\\
        \omega_{in}(\counit) &:=1\in \bolddet(\counit,\partial_{in})\cong \mathbbm{1},\\
        \omega_{out}(\counit)&:=e_0^{-1}\otimes h_0\wedge \sigma h_0 \otimes v\wedge v_0\in \bolddet(\counit,\partial_{out}).
    \end{align*}
    The generator $\omega_{in}(\unit)$ is the canonical generators coming from the isomorphism
    \begingroup
    \setlength{\tabcolsep}{1pt}
    \begin{center}
        \begin{tabular}{ccccc}
            $\bolddet(H_*(\unit,\partial_{in}))$&$\overset{\ref{lem:det}\eqref{det:5}}{\to}$& $\det(C_*(\unit,\partial_{in}))$&$\overset{\ref{prop:computedet}}{\to}$&$\det(E(\unit))^{-1}\otimes \Or(H(\unit))\otimes \det(V(\unit)\setminus L_{in}(\unit))$,\\
            $v_0$&$\mapsto$ & $(h_0\land \sigma h_0)^{-1}\otimes v\wedge v_0$&$\mapsto$& $e_0^{-1}\otimes h_0\wedge \sigma h_0 \otimes v\wedge v_0$.
        \end{tabular}
    \end{center}
    \endgroup
    \noindent An analogous statement holds for the chosen generator $\omega_{out}(\counit)$.\\
    We note that $\omega_{out}(\multi),\omega_{out}(\unit),\omega_{in}(\comulti)$ and $\omega_{in}(\counit)$ all have degree zero. The generators $\omega_{in}(\multi)$ and $\omega_{out}(\comulti)$ are in degree 1 while $\omega_{in}(\unit)$ and $\omega_{out}(\counit)$ are in degree $-1$.
    Then we set for fixed $c,d\in \Z$ the orientations
    \begin{align*}
        \omega(\multi)&:=\omega_{in}(\multi)^{\otimes c}\otimes \omega_{out}(\multi)^{\otimes d}\in {\bolddet}_{c,d}(\multi), &\omega(\unit)&:=\omega_{in}(\unit)^{\otimes c}\otimes \omega_{out}(\unit)^{\otimes d}\in {\bolddet}_{c,d}(\unit),\\
        \omega(\comulti)&:=\omega_{in}(\comulti)^{\otimes c}\otimes \omega_{out}(\comulti)^{\otimes d}\in {\bolddet}_{c,d}(\comulti), &\omega(\counit)&:=\omega_{in}(\counit)^{\otimes c}\otimes \omega_{out}(\counit)^{\otimes d}\in {\bolddet}_{c,d}(\counit).
    \end{align*}
    We note that $\omega(\multi)$ is in degree $c$ and $\omega(\unit)$ is in degree $-c$ while $\omega(\comulti)$ is in degree $d$ and $\omega(\counit)$ is in degree $-d$. As $F$ is enriched over $\Zgrmod_0$, it preserves degrees of maps. We thus have $|\mu|=c$, $|\eta|=-c$, $|\nu|=d$ and $|\varepsilon|=-d$.\\
    We have proven so far that a (symmetric) monoidal functor $F\colon \mathcal{G}_{c,d}\to \mathcal{C}$ is uniquely determined up to isomorphism by the data
    \begin{align*}
        A&=F(*)\in \mathcal{C},\\
        \mu&=F(\omega(\multi))\in \mathcal{C}_{c}(A\otimes A,A), & \eta&=F(\omega(\unit))\in \mathcal{C}_{-c}(\mathbbm{1},A),\\
        \nu&=F(\omega(\comulti))\in C_d(A,A\otimes A), & \varepsilon&=F(\omega(\counit))\in\mathcal{C}_{-d}(A,\mathbbm{1}).
    \end{align*}

    The theorem we are proving states that this data defines a functor \textit{if and only if} certain relations are satisfied. We first argue the only if part of the statement assuming the following lemma:

    \begin{Lem}\label{lem:mainproof}
        In the category $\mathcal{G}_{c,d}$, the following identities hold:
        \begin{align*}
            \omega(\multi)\circ (\omega(\multi)\otimes 1)&=(-1)^c\omega(\multi)\circ (1\otimes \omega(\multi)) \in {\bolddet}_{c,d}(\assgraph),\\
            (-1)^c\omega(\multi)\circ (\omega(\unit)\otimes 1)&=(-1)^{\frac{c(c-1)}{2}}=\omega(\multi)\circ (1\otimes \omega(\unit))\in {\bolddet}_{c,d}(\idgraph)\cong \mathbbm{1},\\
            (\omega(\comulti)\otimes 1)\circ \omega(\comulti)&=(-1)^d(1 \otimes \omega(\comulti))\circ \omega(\comulti)\in {\bolddet}_{c,d}(\coassgraph),\\
            (-1)^d(1\otimes \omega(\counit))\circ \omega(\comulti)&=(-1)^{\frac{d(d-1)}{2}}=(\omega(\counit)\otimes 1)\circ \omega(\comulti)\in {\bolddet}_{c,d}(\idgraph)\cong \mathbbm{1},\\
            (\omega(\multi)\otimes 1)\circ (1\otimes \omega(\comulti))&=(-1)^{cd}\omega(\comulti)\circ \omega(\multi)=(1\otimes \omega(\multi))\circ (\omega(\comulti)\otimes 1)\in {\bolddet}_{c,d}(\frobgraph).
        \end{align*}
        Moreover in the category $\GrCob_{c,d}$, the automorphism of the graph $\multi$ which switches the two incoming edges acts by $(-1)^c$ on $\bolddet_{c,d}(\multi)$ and, in the category $\fatGrCob_{c,d}$, the fat graph automorphism of the fat graph $\symgraph$ which switches the two incoming edges acts by $(-1)^c$ on $\bolddet_{c,d}(\symgraph)$.
    \end{Lem}

    These identities are preserved under $F$ because $F$ is linear and monoidal. The lemma thus shows that the relations are satisfied by $\mu$, $\eta$, $\nu$ and $\varepsilon$ if they come from a functor $F$. This thus concludes the only if part.
    
    Before proving the lemma, we show the if part of the statement that the data $(A,\mu,\eta,\nu,\varepsilon)$ satisfying the relations defines a functor:\\
    Consider a morphism $f\in \mathcal{G}(X,Y)$. As we have
    \begin{align*}
        {\mathcal{G}_{c,d}}(X,Y)&\cong\bigoplus_{[G]\in \mathcal{G}(X, Y)} {\bolddet}_{c,d}(G)_{\pi_1(|\underline{\mathcal{G}}(X,Y)|,G)},
    \end{align*}
    we may assume that $f$ is a homogeneous generator, i.e., it is given by the $\pi_1(|\underline{\mathcal{G}}(X,Y)|,G)$-orbit of $ \omega(G)$ for some generator $\omega(G)\in \bolddet_{c,d}(G)$. We write $[G]$ as a composition of the graphs $\multi$, $\unit$, $\comulti$ and $\counit$ as in Figure \ref{fig:graph decomp}. We then define $F(\omega(G))$ as the corresponding composition of $\mu$, $\eta$, $\nu$ and $\varepsilon$.\\
    The map $F$ is by definition linear and degree-preserving, i.e., it is enriched over $\Zgrmod_0$. Moreover up to choosing the right decomposition of $G$, this is functorial and (symmetric) monoidal. However, our definition of $F(f)$ depends on two choices: the decomposition of $G$ and the choice of representative of the $\pi_1(|\underline{\mathcal{G}}(X,Y)|,G)$-orbit.\\    
    We first argue why the definition is independent of the choice of decomposition: in \cite{kock2004frobenius,lauda2008open}, Kock and Lauda and Pfeiffer describe a normal form for a decomposition of a closed or an open surface. Moreover, they show that any decomposition can be brought into this normal form using the six generating relations: associativity, unitality, coassociativity, counitality, the Frobenius relation and possibly commutativity or symmetry. Thus any two decompositions are related to each other by a sequence of those six relations.\\
    An analogous argument for (fat) graphs shows that any two decompositions into the four graphs $\multi$, $\unit$, $\comulti$ and $\counit$ are related by a sequence of those six relations.\\
    Therefore, any two decompositions are related by applying the relations in Lemma \ref{lem:mainproof}. If the morphisms $\mu$, $\eta$, $\nu$ and $\varepsilon$ satisfy the corresponding relation in $\mathcal{C}$, the definition of $F(\omega(G))$ as above is independent of the choice of decomposition.\\
    It remains to show that our definition of $F$ is independent of choice of representative in the $\pi_1(|\underline{\mathcal{G}}(X,Y)|,G)$-orbit. This is equivalent to showing that our definition of $F(\omega(G))$ is invariant under the $\pi_1(|\underline{\mathcal{G}}(X,Y)|,G)$-action:\\
    Any element in $\pi_1(|\underline{\mathcal{G}}(X,Y)|,G)$ is represented by a zig-zag of edge collapses and (fat) graph isomorphisms:
    \begin{center}
        \begin{tikzcd}
            G=G_0 \ar[r] & G_1 & G_2  \ar[r] \ar[l] & G_3 &  \dots \ar[r] \ar[l] & G_{n-1} & G_n=G. \ar[l] 
        \end{tikzcd}
    \end{center}
    Applying one of the relations in Lemma \ref{lem:mainproof} to a decomposition $G$ produces a new decomposition $G'$. The relation then specifies a zig-zag 
    \begin{center}
        \begin{tikzcd}
            G\ar[r] & G''  & G' \ar[l]
        \end{tikzcd}
    \end{center}
    where $G$ and $G'$ are given by a different decomposition. We show that all zig-zags can be lifted to a sequence of zig-zags coming from the relations in Lemma \ref{lem:mainproof}.\\
    We can lift each graph $G_{2i}$ to a decomposition into the four elementary graphs: $\widetilde{G}_{2i}\to G_{2i}$ where $\widetilde{G}_{2i}$ has a decomposition into the four elementary graphs. We fix such a decomposition for all $\widetilde{G}_{2i}$.\\
    We thus have the following diagram:
    \begin{center}
        \begin{tikzcd}
            \widetilde{G}_0 \ar[d, equal] & & \widetilde{G}_2 \ar[d] & & \dots & & \widetilde{G}_n \ar[d, equal] \\
            G\ar[r] & G_1 & G_2  \ar[r] \ar[l] & G_3 &  \dots \ar[r] \ar[l] & G_{n-1} & G. \ar[l] 
        \end{tikzcd}
    \end{center}
    Each $\widetilde{G}_{2i}\to G_{2i}\to G_{2i\pm 1}$ gives a decomposition of $G_{2i\pm 1}$. We thus can rewrite the above diagram suggestively as
    \begin{center}
        \begin{tikzcd}
            \widetilde{G}_0 \ar[d, equal] \ar[r, equal] &\widetilde{G}_0 \ar[d] & \widetilde{G}_2  \ar[d] \ar[r,equal]& \widetilde{G}_2 \ar[d]\ar[r, equal] & \widetilde{G}_2 \ar[d] & \dots & \widetilde{G}_n\ar[d] \ar[r, equal]& \widetilde{G}_n \ar[d, equal] \ar[l, equal]\\
            G \ar[r] & G_1 & G_1 \ar[l, equal]& G_2 \ar[r] \ar[l] & G_3 & \dots \ar[r] \ar[l] & G_{n-1}  & G\ar[l].
        \end{tikzcd}
    \end{center}
    The decompositions $\widetilde{G}_{2i}\to G_{2i+1}$ and $\widetilde{G}_{2i+2}\to G_{2i+1}$ are two decompositions of $G_{2i+1}$. Therefore $\widetilde{G}_{2i}$ and $\widetilde{G}_{2i+2}$ are related by a sequence of the relations in Lemma \ref{lem:mainproof}. Thus there exists a zig-zag coming from those relations fitting into the lifting problem 
    \begin{center}
        \begin{tikzcd}
            \widetilde{G}_{2i}\ar[d] \ar[r, dashed] & \widetilde{G}_{2i+2}\ar[d]\\
            G_{2i+1}\ar[r, equal] & G_{2i+1}.
        \end{tikzcd}
    \end{center}
    This lifts the whole zig-zag to a zig-zag coming from the relations in Lemma \ref{lem:mainproof}.\\
    The action of the element of $\pi_1(|\underline{\mathcal{C}}(X,Y)|,G)$ can thus be computed by applying the relations in Lemma \ref{lem:mainproof}.\\
    Consequently, if $\mu$, $\eta$, $\nu$ and $\varepsilon$ satisfy the corresponding relations in $\mathcal{C}$, our definition of $F(\omega(G))$ does not depend on the representative of the $\pi_1(\underline{\mathcal{G}}(X,Y),G)$-orbit as it is invariant under the $\pi_1(\underline{\mathcal{G}}(X,Y),G)$-action.\\
    This shows that the data $\mu$, $\eta$, $\nu$ and $\varepsilon$ satisfying the relations lets us define a functor $F$ and thus concludes the only if part of the statement. The only thing left to prove is Lemma \ref{lem:mainproof}. The proof makes up the rest of the paper.   
\begin{proof}[Proof of Lemma \ref{lem:mainproof}]
    We prove the seven identities using Lemmas \ref{lem:edge collapse} and \ref{lem:compute compo}.
    \subsection*{Graded associativity.} We show that
    \begin{align}\label{eq:ass}
        \omega(\multi)\circ (\omega(\multi)\otimes 1)=(-1)^c\omega(\multi)\circ (1\otimes \omega(\multi)) \in {\bolddet}_{c,d}(\assgraph)
    \end{align}
    holds in $\mathcal{G}_{c,d}$.\\
    We start computing
    \begin{align*}
        \omega(\multi)\circ (\omega(\multi)\otimes 1)=&\ (\omega_{in}(\multi)^{\otimes c}\otimes \omega_{out}(\multi)^{\otimes d})\circ (\omega_{in}(\multi)^{\otimes c}\otimes \omega_{out}(\multi)^{\otimes d}\otimes 1)\\
        =&\  \comp_1(\omega_{in}(\multi)^{\otimes c}\otimes \omega_{in}(\multi)^{\otimes c})\otimes \comp_1(\omega_{out}(\multi)^{\otimes d}\otimes \omega_{out}(\multi)^{\otimes d})\\
        =&\  (-1)^\frac{c(c-1)}{2}\comp_1(\omega_{in}(\multi)\otimes \omega_{in}(\multi))^{\otimes c}\otimes \comp_1(\omega_{out}(\multi)\otimes \omega_{out}(\multi))^{\otimes d}.
    \end{align*}
    Here $\comp_1$ is the composition we get by gluing one copy of $\multi$ at $v_0$ to $v_1$ of the other copy. The sign $(-1)^{\frac{c(c-1)}{2}}$ comes from the fact that we have to reshuffle $\omega_{in}(\multi)^{\otimes c}\otimes \omega_{in}(\multi)^{\otimes c}$ into $c$ pairs of $\omega_{in}(\multi)\otimes \omega_{in}(\multi)$ (via a riffle shuffle). For $\omega_{out}(\multi)$ no such sign appears as they are in degree 0.\\
    Analogously, we compute
    \begin{align*}
        \omega(\multi)\circ (1\otimes \omega(\multi))=(-1)^\frac{c(c-1)}{2}\comp_2(\omega_{in}(\multi)\otimes \omega_{in}(\multi))^{\otimes c}\otimes \comp_2(\omega_{out}(\multi)\otimes \omega_{out}(\multi))^{\otimes d}.
    \end{align*}
    Here $\comp_2$ is the composition we get by gluing one copy of $\multi$ at $v_0$ to $v_2$ of the other copy because $v_2$ is the vertex corresponding to the second input.\\
    It remains to compare $\comp_1(\omega_{in}(\multi)\otimes \omega_{in}(\multi))$ to $\comp_2(\omega_{in}(\multi)\otimes \omega_{in}(\multi))$ and $\comp_1(\omega_{out}(\multi)\otimes \omega_{out}(\multi))$ to $\comp_2(\omega_{out}(\multi)\otimes \omega_{out}(\multi))$.\\
    We note that because of triviality of the $\omega_{out}(\multi)$ orientations, we have
    \begin{align*}
        \comp_1(\omega_{out}(\multi)\otimes \omega_{out}(\multi))=\comp_2(\omega_{out}(\multi)\otimes \omega_{out}(\multi)).
    \end{align*}
    Equation \eqref{eq:ass} then follows from the following identity
    \begin{align}\label{eq:ass2}
        \comp_1(\omega_{in}(\multi)\otimes \omega_{in}(\multi))=-\comp_2(\omega_{in}(\multi)\otimes \omega_{in}(\multi)).
    \end{align}
    To show this identity, we start with $\comp_1(\omega_{in}(\multi)\otimes \omega_{in}(\multi))$. For this, we denote the edges, half-edges and vertices of the left copy of $\multi$ by $e,e_0,e_1,e_2$ and so on, as before. This represents the multiplication that is applied second. We denote the edges, half-edges and vertices of the right copy of $\multi$ by $e',e_0',e_1',e_2'$ and so on. This represents the multiplication that is applied first.\\
    We compute using Lemma \ref{lem:compute compo} and Remark \ref{rem:sign compo}:
    \begin{align*}
        \comp_1(\omega_{in}(\multi)\otimes \omega_{in}(\multi))=&\ \comp_1((e_2\wedge e_1\wedge e_0)^{-1}\otimes h_0\wedge \sigma h_0 \wedge h_1\wedge \sigma h_1 \wedge h_2\wedge \sigma h_2\otimes v\wedge v_0\\
        &\otimes (e_2'\wedge e_1'\wedge e_0')^{-1}\otimes h_0'\wedge \sigma h_0' \wedge h_1'\wedge \sigma h_1' \wedge h_2'\wedge \sigma h_2'\otimes v'\wedge v_0')\\
        =&\ -(e_2\wedge e_1\wedge e_0 \wedge e_2' \wedge e_1' \wedge e_0')^{-1}\\
        &\otimes h_0\wedge \sigma h_0 \wedge h_1\wedge \sigma h_1 \wedge h_2\wedge \sigma h_2\wedge h_0'\wedge \sigma h_0' \wedge h_1'\wedge \sigma h_1' \wedge h_2'\wedge \sigma h_2'\\
        &\otimes v\wedge v_0 \wedge v' \wedge v_0'
    \end{align*}
    \begin{figure}[H]
        \centering
        \includegraphics[width=0.57\linewidth]{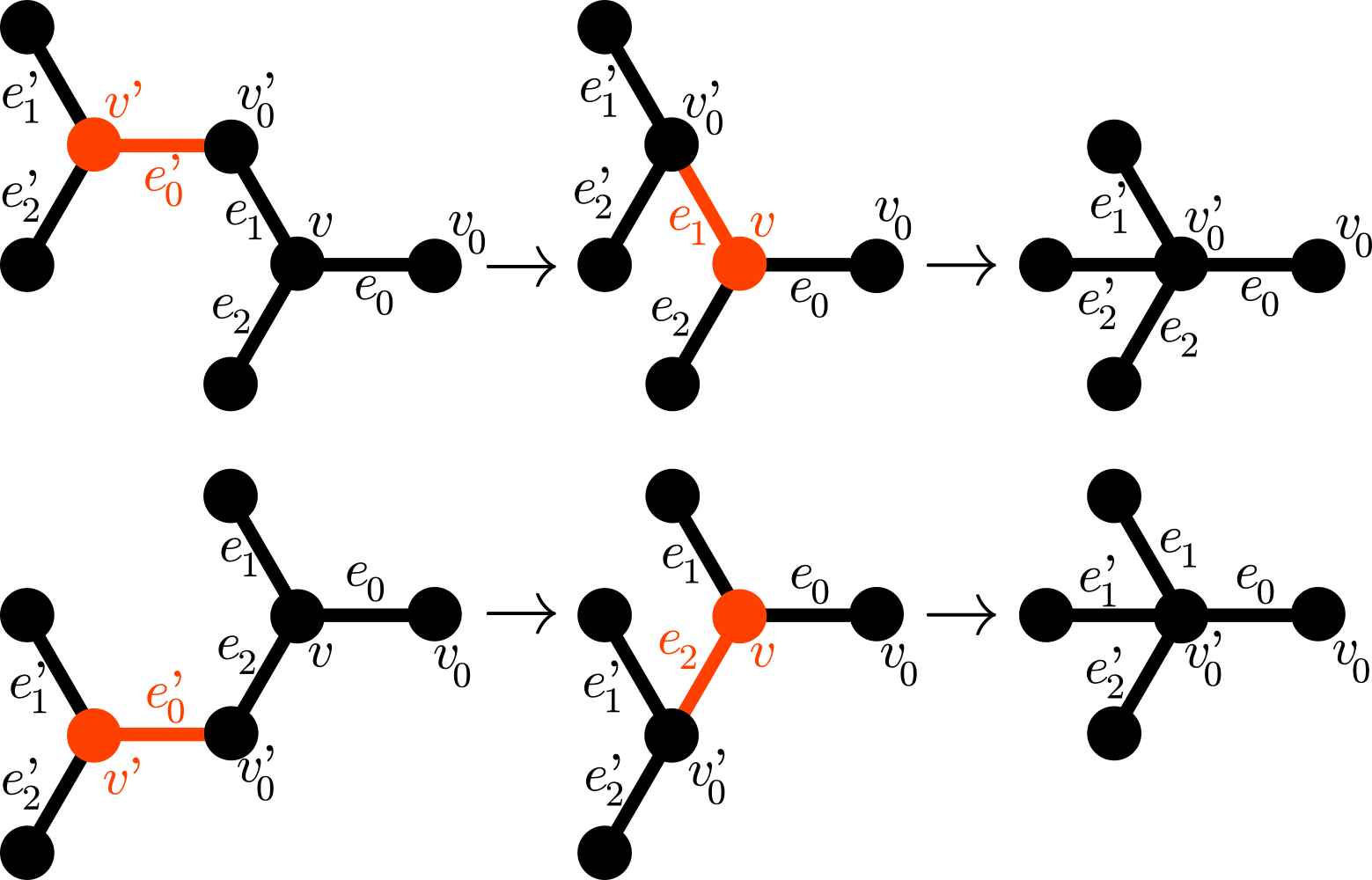}
        \caption{The two composition of two multiplications give the same graph after collapsing two edges.}
        \label{fig:ass strat}
    \end{figure}
    As we glued $v_0'$ to $v_1$, we can collapse $e_0'$ and $e_1$ to obtain the graph $\assgraph$ (see top of Figure \ref{fig:ass strat}). We recall the naming convention that $s(\sigma h_i)$ is an internal vertex. Therefore when collapsing $e_0'$, we can apply Lemma \ref{lem:edge collapse} using that $s(\sigma h_0')=v'$ is the vertex we want to delete from the orientation. Thereafter deleting $e_1$, we note that $s(\sigma h_1)=v$ is the vertex we delete in the process. Importantly to apply Lemma \ref{lem:edge collapse}, we must reshuffle the edge and vertex, we want to delete to the first position. This gives some signs. We do the reshuffle explicitly in the first computation and thereafter only denote the corresponding sign after the edge collapse. We find:
    \begin{align*}
        -(e_2\wedge e_1&\wedge e_0 \wedge e_2' \wedge e_1' \wedge e_0')^{-1}\\
        \otimes&\ h_0\wedge \sigma h_0 \wedge h_1\wedge \sigma h_1 \wedge h_2\wedge \sigma h_2\wedge h_0'\wedge \sigma h_0' \wedge h_1'\wedge \sigma h_1' \wedge h_2'\wedge \sigma h_2'\otimes v\wedge v_0 \wedge v' \wedge v_0'\\
        =&\ (e_0'\wedge e_2\wedge e_1\wedge e_0 \wedge e_2' \wedge e_1' )^{-1}\\
        &\otimes h_0'\wedge \sigma h_0'\wedge h_0\wedge \sigma h_0 \wedge h_1\wedge \sigma h_1 \wedge h_2\wedge \sigma h_2 \wedge h_1'\wedge \sigma h_1' \wedge h_2'\wedge \sigma h_2' \otimes v'\wedge v\wedge v_0 \wedge v_0'\\
        \mapsto&\ (e_2\wedge e_1\wedge e_0 \wedge e_2' \wedge e_1' )^{-1}\\
        &\otimes h_0\wedge \sigma h_0 \wedge h_1\wedge \sigma h_1 \wedge h_2\wedge \sigma h_2 \wedge h_1'\wedge \sigma h_1' \wedge h_2'\wedge \sigma h_2' \otimes v\wedge v_0 \wedge v_0'\\
        \mapsto&\ -(e_2\wedge e_0 \wedge e_2' \wedge e_1' )^{-1}\otimes h_0\wedge \sigma h_0 \wedge h_2\wedge \sigma h_2 \wedge h_1'\wedge \sigma h_1' \wedge h_2'\wedge \sigma h_2'\otimes v_0 \wedge v_0'.
    \end{align*}
    The edges of $\assgraph$ are here labelled (from top to bottom) by $e_1'$, $e_2'$ and $e_2$ as in the top right of Figure \ref{fig:ass strat}.\\
    Analogously, we compute
    \begin{align*}
        \comp_2(\omega_{in}(\multi)\otimes \omega_{in}(\multi))= -(e_2\wedge 
        \comp_2(\omega_{in}(\multi)\otimes \omega_{in}(\multi))=&  -(e_2\wedge e_1\wedge e_0 \wedge e_2' \wedge e_1' \wedge e_0')^{-1}\\
        &\otimes h_0\wedge \sigma h_0 \wedge h_1\wedge \sigma h_1 \wedge h_2\wedge \sigma h_2\wedge h_0'\wedge \sigma h_0' \wedge h_1'\wedge \sigma h_1' \wedge h_2'\wedge \sigma h_2'\\
        &\otimes v\wedge v_0 \wedge v' \wedge v_0'.
    \end{align*}
    But as we glued $v_0'$ to $v_2$, we can now collapse $e_0'$ and $e_2$ (see bottom of Figure \ref{fig:ass strat}). In this process, we first delete the vertex $s(\sigma h_0')=v'$ and then the vertex $s(\sigma h_2)=v$. We thus find
    \begin{align*}
        -(e_2\wedge e_1&\wedge e_0 \wedge e_2' \wedge e_1' \wedge e_0')^{-1}\\
        \otimes &h_0\wedge \sigma h_0 \wedge h_1\wedge \sigma h_1 \wedge h_2\wedge \sigma h_2\wedge h_0'\wedge \sigma h_0' \wedge h_1'\wedge \sigma h_1' \wedge h_2'\wedge \sigma h_2' \otimes v\wedge v_0 \wedge v' \wedge v_0'\\
        \mapsto &(e_2\wedge e_1\wedge e_0 \wedge e_2' \wedge e_1')^{-1}\\
        &\otimes h_0\wedge \sigma h_0 \wedge h_1\wedge \sigma h_1 \wedge h_2\wedge \sigma h_2\wedge h_1'\wedge \sigma h_1' \wedge h_2'\wedge \sigma h_2' \otimes v\wedge v_0 \wedge v_0'\\
        \mapsto &(e_1\wedge e_0 \wedge e_2' \wedge e_1')^{-1}\otimes h_0\wedge \sigma h_0 \wedge h_1\wedge \sigma h_1 \wedge h_1'\wedge \sigma h_1' \wedge h_2'\wedge \sigma h_2' \otimes v_0 \wedge v_0'.
    \end{align*}
    The edges of $\assgraph$ are here labelled (from top to bottom) by $e_1$, $e_1'$ and $e_2'$ as in the bottom right of Figure \ref{fig:ass strat}. This differs from the naming scheme of the other composition by an even permutation. This shows that the two orientations differ by the sign $-1$. We have proved \eqref{eq:ass2} and thus \eqref{eq:ass}.
    \subsection*{Graded unitality.} We want to compute
    \begin{align}\label{eq:unit}
        (-1)^c\omega(\multi)\circ (\omega(\unit)\otimes 1)=(-1)^{\frac{c(c-1)}{2}}=\omega(\multi)\circ (1\otimes \omega(\unit))\in {\bolddet}_{c,d}(\idgraph)\cong \mathbbm{1}.
    \end{align}
    As for graded associativity, we compute the signs from the riffle shuffle
    \begin{align*}
        \omega(\multi)\circ (\omega(\unit)\otimes 1)=(-1)^\frac{c(c-1)}{2}\comp_1(\omega_{in}(\multi)\otimes \omega_{in}(\unit))^{\otimes c}\otimes \comp_1( \omega_{out}(\multi)\otimes \omega_{out}(\unit))^{\otimes d}
    \end{align*}
    and 
    \begin{align*}
        \omega(\multi)\circ (1\otimes \omega(\unit))=(-1)^\frac{c(c-1)}{2}\comp_2( \omega_{in}(\multi)\otimes \omega_{in}(\unit))^{\otimes c}\otimes \comp_2( \omega_{out}(\multi)\otimes \omega_{out}(\unit))^{\otimes d}.
    \end{align*}
    We note that, by triviality of the $\omega_{out}$-orientations, we have
    \begin{align*}
        \comp_1( \omega_{out}(\multi)\otimes \omega_{out}(\unit))=1=\comp_2( \omega_{out}(\multi)\otimes \omega_{out}(\unit))\in\bolddet(\idgraph,\partial_{out}).
    \end{align*}
    The equation \eqref{eq:unit} then follows from the following identity
    \begin{align}\label{eq:unit2}
        -\comp_1( \omega_{in}(\multi)\otimes \omega_{in}(\unit))=1=\comp_2( \omega_{in}(\multi)\otimes \omega_{in}(\unit))\in \bolddet(\idgraph,\partial_{in})\cong \mathbbm{1}.
    \end{align}
    To show this, we denote the edges, half-edges and vertices of $\unit$ by $e_0$ and so on. On the other hand, we denote the edges, half-edges and vertices of $\multi$ by $e_0',e_1',e_2'$ and so on. We thus compute by Lemma \ref{lem:compute compo} and Remark \ref{rem:sign compo}:
    \begin{align*}
        \comp_1(\omega_{in}(\multi)\otimes \omega_{in}(\unit))=&\ \comp_1((e_2'\wedge e_1'\wedge e_0')^{-1}\otimes h_0'\wedge \sigma h_0' \wedge h_1'\wedge \sigma h_1' \wedge h_2'\wedge \sigma h_2'\otimes v'\wedge v_0'\\
        &\otimes (e_0)^{-1}\otimes h_0\wedge \sigma h_0\otimes v\wedge v_0)\\
        =&\ -( e_2'\wedge e_1'\wedge e_0'\wedge e_0)^{-1}\\
        &\otimes h_0'\wedge \sigma h_0' \wedge h_1'\wedge \sigma h_1' \wedge h_2'\wedge \sigma h_2'\wedge h_0\wedge \sigma h_0\otimes  v'\wedge v_0'\wedge v\wedge v_0.
    \end{align*}
    \begin{figure}[H]
        \centering
        \includegraphics[width=0.75\linewidth]{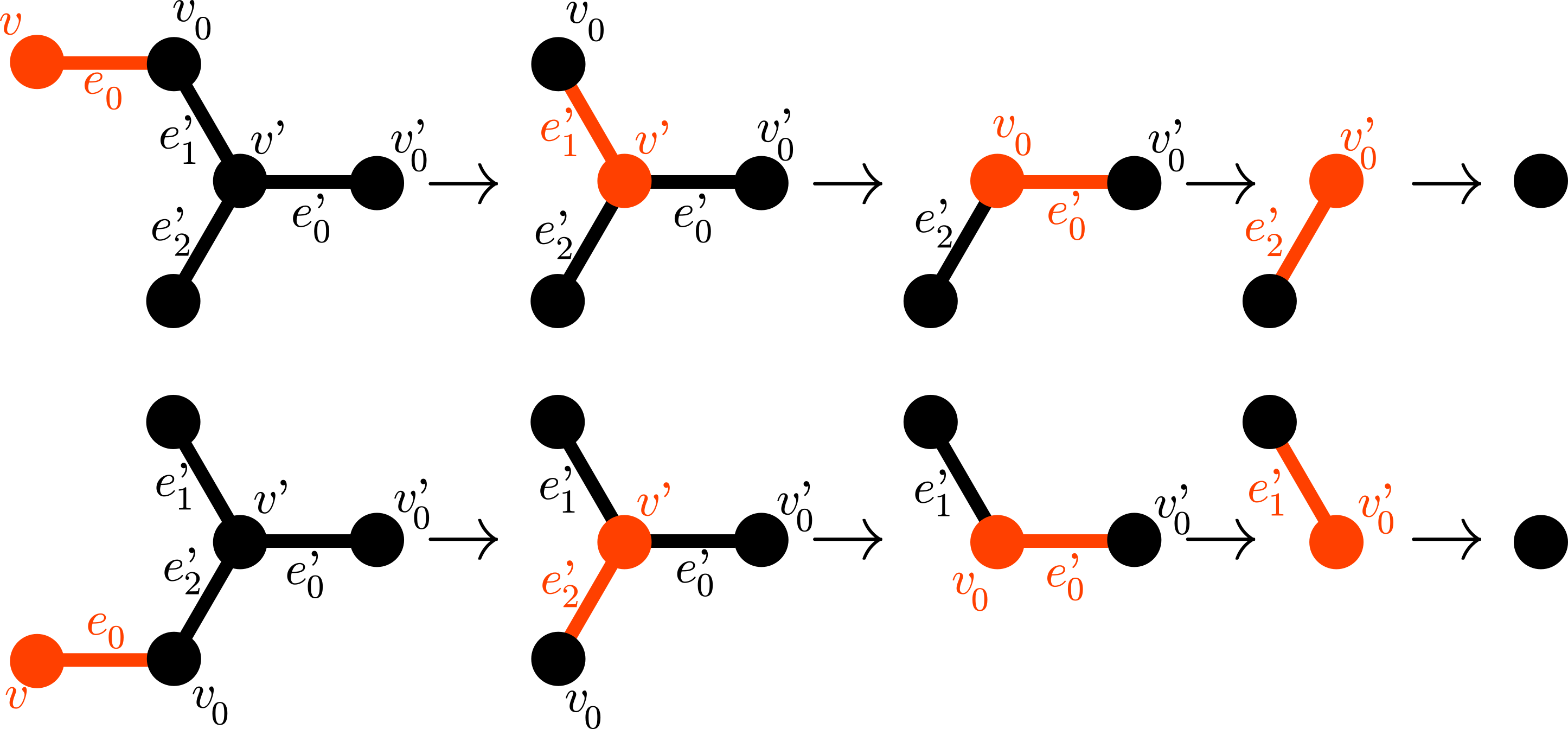}
        \caption{The two composition of a unit and a multiplications give the graph with one vertex representing $\id$ after collapsing four edges.}
        \label{fig:unit strat}
    \end{figure}
    Now we can collapse $e_0$ and $e_1'$ as in the top of Figure \ref{fig:unit strat}. We use Lemma \ref{lem:edge collapse} and note that when we collapse $e_0$, we have to delete $s(\sigma h_0)=v$. If we then collapse $e_1'$, we note that $s(\sigma h_1')=v'$. We compute:
    \begin{align*}
        -(e_2'\wedge e_1'&\wedge e_0'\wedge e_0)^{-1}\otimes h_0'\wedge \sigma h_0' \wedge h_1'\wedge \sigma h_1' \wedge h_2'\wedge \sigma h_2'\wedge h_0\wedge \sigma h_0\otimes  v'\wedge v_0'\wedge v\wedge v_0\\
        \mapsto&\ (e_2'\wedge e_1'\wedge e_0')^{-1}\otimes h_0'\wedge \sigma h_0' \wedge h_1'\wedge \sigma h_1' \wedge h_2'\wedge \sigma h_2'\otimes  v'\wedge v_0'\wedge v_0\\
        \mapsto& -(e_2'\wedge e_0')^{-1}\otimes h_0'\wedge \sigma h_0' \wedge h_2'\wedge \sigma h_2'\otimes   v_0'\wedge v_0.
    \end{align*}
    Next we collapse $e_0'$ which deletes $s(\sigma h_0')=v_0$:
    \begin{align*}
         -(e_2'\wedge e_0')^{-1}&\otimes h_0'\wedge \sigma h_0' \wedge h_2'\wedge \sigma h_2'\otimes   v_0'\wedge v_0\\
         \mapsto & -(e_2')^{-1}\otimes  h_2'\wedge \sigma h_2'\otimes   v_0'.
    \end{align*}
    Finally collapsing $e_2'$, we delete $s(\sigma h_2')=v_0'$ and find:
    \begin{align*}
        -(e_2')^{-1}\otimes  h_2'\wedge \sigma h_2'\otimes   v_0'\mapsto -1.
    \end{align*}
    This shows the left equality in \eqref{eq:unit2}.\\
    For the other equation, we compute:
    \begin{align*}
        \comp_2( \omega_{in}(\multi)\otimes \omega_{in}(\unit))=&-( e_2'\wedge e_1'\wedge e_0'\wedge e_0)^{-1}\\
        &\otimes h_0'\wedge \sigma h_0' \wedge h_1'\wedge \sigma h_1' \wedge h_2'\wedge \sigma h_2'\wedge h_0\wedge \sigma h_0\otimes  v'\wedge v_0'\wedge v\wedge v_0.
    \end{align*}
    But now we collapse first $e_0$ and at the same time deleting $s(\sigma h_0)=v$. Then we collapse $e_2'$ and delete $s(\sigma h_2')=v'$. Then we delete $e_0'$ and delete $s(\sigma h_0')=v_0$. Finally, we collapse $e_1'$ and delete $s(\sigma h_1')=v_0'$ as in the bottom of Figure \ref{fig:unit strat}. We thus find
    \begin{align*}
        -(e_0\wedge e_2'\wedge e_1'\wedge e_0')^{-1}\otimes h_0\wedge \sigma h_0\wedge h_0'\wedge \sigma h_0' \wedge h_1'\wedge \sigma h_1' \wedge h_2'\wedge \sigma h_2'\otimes v\wedge v_0\wedge v'\wedge v_0' \mapsto 1.
    \end{align*}
    One can see that by either doing the calculation explicitly. Or one can notice that we are deleting the vertices in the same order as before but switched the order of the two edges $e_1'$ and $e_2'$ in which we delete them (see bottom of Figure \ref{fig:unit strat}). This shows \eqref{eq:unit2} and thus \eqref{eq:unit}.
    \subsection*{Graded coassociativity.} This computation is dual to the computation of graded associativity, we thus obtain: 
    \begin{align*}
        (\omega(\comulti)\otimes 1)\circ \omega(\comulti)=(-1)^d(1 \otimes \omega(\comulti))\circ \omega(\comulti)\in {\bolddet}_{c,d}(\coassgraph).
    \end{align*}
    \subsection*{Graded counitality.} This computation is dual to the computation of graded unitality, we thus obtain: 
    \begin{align*}
         (-1)^d(1\otimes \omega(\counit))\circ \omega(\comulti)=(-1)^{\frac{d(d-1)}{2}}=(\omega(\counit)\otimes 1)\circ \omega(\comulti)\in {\bolddet}_{c,d}(\idgraph)\cong \mathbbm{1}.
    \end{align*}
    We note that the dual term to $\omega(\multi)\circ (\omega(\unit)\otimes 1)$ is $(1\otimes \omega(\counit))\circ \omega(\comulti)$.
    \subsection*{The graded Frobenius relation.} We want to show that
    \begin{align}\label{eq:frob}
        (\omega(\mu)\otimes 1)\circ (1\otimes \omega(\nu))=(-1)^{cd}\omega(\nu)\circ \omega(\mu)=(1\otimes \omega(\mu))\circ (\omega(\nu)\otimes 1)\in {\bolddet}_{c,d}(\frobgraph).
    \end{align}
    We start computing
    \begin{align*}
        (\omega(\mu)\otimes 1)\circ (1\otimes \omega(\nu))&=(\omega_{in}(\mu)^{\otimes c}\otimes \omega_{out}(\mu)^{\otimes d}\otimes 1)\circ (1\otimes \omega_{in}(\nu)^{\otimes c}\otimes \omega_{out}(\nu)^{\otimes d})\\
        &=(\omega_{in}(\mu)^{\otimes c}\otimes 1)\circ (1\otimes \omega_{in}(\nu)^{\otimes c})\otimes (\omega_{out}(\mu)^{\otimes d}\otimes 1)\circ (1\otimes \omega_{out}(\nu)^{\otimes d})\\
        &=\comp_{1,2}(\omega_{in}(\mu)\otimes \omega_{in}(\nu))^{\otimes c}\otimes \comp_{1,2}(\omega_{out}(\mu)\otimes \omega_{out}(\nu))^{\otimes d}.
    \end{align*}
    Here $\comp_{1,2}$ is the composition given by gluing together the first outgoing vertex of $\comulti$ to the second incoming vertex of $\multi$.\\
    Analogously, we compute
    \begin{align*}
        (1\otimes \omega(\mu))\circ (\omega(\nu)\otimes 1)=\comp_{2,1}(\omega_{in}(\mu)\otimes \omega_{in}(\nu))^{\otimes c}\otimes \comp_{2,1}(\omega_{out}(\mu)\otimes \omega_{out}(\nu))^{\otimes d}.
    \end{align*}
    Here $\comp_{2,1}$ is the composition given by gluing together the second outgoing vertex of $\comulti$ to the first incoming vertex of $\multi$.\\
    In contrast, we compute
    \begin{align*}
        \omega(\nu)\circ \omega(\mu)=&\ (\omega_{in}(\nu)^{\otimes c}\otimes \omega_{out}(\nu)^{\otimes d})\circ (\omega_{in}(\mu)^{\otimes c}\otimes \omega_{out}(\mu)^{\otimes d})\\
        = &\  (-1)^{cd}(\omega_{in}(\nu)^{\otimes c}\circ \omega_{in}(\mu)^{\otimes c})\otimes (\omega_{out}(\nu)^{\otimes d}\circ \omega_{out}(\mu)^{\otimes d}).
    \end{align*}
    The sign $(-1)^{cd}$ comes from the fact that we have to move $\omega_{out}(\nu)^{\otimes d}$ which is in degree $d$ past $\omega_{in}(\mu)^{\otimes c}$ which is in degree $c$. By definition, we have:
    \begin{align*}
        (-1)^{cd}(\omega_{in}(\nu)^{\otimes c}&\circ \omega_{in}(\mu)^{\otimes c})\otimes (\omega_{out}(\nu)^{\otimes d}\circ \omega_{out}(\mu)^{\otimes d})\\
        =&\ (-1)^{cd}\comp(\omega_{in}(\nu)\otimes \omega_{in}(\mu))^{\otimes c}\otimes \comp(\omega_{out}(\nu)\otimes \omega_{out}(\mu))^{\otimes d}.
    \end{align*}
    The equation \eqref{eq:frob} now follows from the following identities
    It holds that
    \begin{align*}
        \comp_{1,2}(\omega_{in}(\mu)\otimes \omega_{in}(\nu))=\comp(\omega_{in}(\nu)\otimes \omega_{in}(\mu))=\comp_{2,1}(\omega_{in}(\mu)\otimes \omega_{in}(\nu))
    \end{align*}
    and
    \begin{align*}
        \comp_{1,2}(\omega_{out}(\mu)\otimes \omega_{out}(\nu))=\comp(\omega_{out}(\nu)\otimes \omega_{out}(\mu))=\comp_{2,1}(\omega_{out}(\mu)\otimes \omega_{out}(\nu)).
    \end{align*}
    This can be seen either by doing the hands-on computation using Lemmas \ref{lem:edge collapse} and \ref{lem:compute compo} or by noting that $\omega_{in}(\nu)$ and $ \omega_{out}(\mu)$ are trivial. Composition with them thus gives the trivial map. 
    
    \subsection*{Graded Commutativity.} We want to show that the automorphism of the graph $\multi$ which switches the two incoming edges acts by $(-1)^c$ on $\omega(\multi)$. For this, we note that this automorphism acts trivially on $\bolddet(\multi,\partial_{out})\cong \mathbbm{1}$. On the other hand, we compute the action $\bolddet(\multi,\partial_{in})$ as
    \begin{align*}
        \bolddet(\multi,\partial_{in}) &\to \bolddet(\multi,\partial_{in}),\\
        (e_2\wedge e_1\wedge e_0)^{-1}\otimes h_0\wedge \sigma h_0 \wedge h_1\wedge &\sigma h_1 \wedge h_2\wedge \sigma h_2\otimes v\wedge v_0\\
        &\mapsto (e_1\wedge e_2\wedge e_0)^{-1}\otimes h_0\wedge \sigma h_0 \wedge h_2\wedge \sigma h_2 \wedge h_1\wedge \sigma h_1\otimes v\wedge v_0.
    \end{align*}
    It only changes the order of $e_1$ and $e_2$ and their half-edges. In other words, it acts by $-1$ on $\bolddet(\multi,\partial_{in})$. This combines to a $(-1)^c$-action on $\bolddet_{c,d}(\multi)$.
    \subsection*{Graded Symmetry.} We note that flipping the two incoming vertices of $\multi$ is not a fat graph automorphism. As before, we call the edges, half-edges and vertices of $\multi$ $e_0,e_1,e_2$ and so on. On the other hand, we call the edges, half-edges and vertices of $\counit$ $e_0'$, $h_0, \sigma h_0$ and $v_0',v'$. We glue $v_0$ to $v_0'$, then we can collapse the edge $e_0$ and $e_0'$. On the resulting graph $\symgraph$, flipping $v_1$ with $v_2$ is a fat graph automorphism. This acts on $\bolddet(\symgraph,\partial_{in})$ by $-1$ and on $\bolddet(\symgraph,\partial_{out})$ trivially. This combines to a $(-1)^c$-action on $\bolddet_{c,d}(\symgraph)$.\\
    Therefore, we have computed all signs appearing in the identities in Lemma \ref{lem:mainproof}. This concludes the proof of Lemma \ref{lem:mainproof}, the proof of Theorem \ref{thm:gen and rel} and thus the paper.
\end{proof}

\end{proof}

\normalem
\bibliographystyle{alpha}
\bibliography{main.bib}

\end{document}